\documentclass[10pt]{amsart}
\usepackage[dvipdfmx]{graphicx}
\usepackage{amsmath,amscd,amsthm,amsxtra,yhmath,stackrel}
\usepackage{epsfig,graphics,color,colortbl}
\usepackage{amssymb,latexsym}
\usepackage{mathrsfs}
\usepackage{bm}
\usepackage[poly,all]{xy}
\usepackage{hyperref}
\usepackage{tikz-cd}
\usepackage[draft]{todonotes}

\newtheorem{thm}{\bf Theorem}[section]
\newtheorem{df}[thm]{\bf Definition}
\newtheorem{prop}[thm]{\bf Proposition}
\newtheorem{cor}[thm]{\bf Corollary}
\newtheorem{lem}[thm]{\bf Lemma}
\newtheorem{rem}[thm]{\bf Remark}

\numberwithin{equation}{section}

\newcommand{\mc}{\mathcal}
\newcommand{\mf}{\mathfrak}
\newcommand{\ms}{\mathscr}
\newcommand{\mb}{\bm}

\newcommand{\pf}{\noindent{\bfseries Proof. }}
\newcommand{\ov}{\overline}
\newcommand{\un}{\underline}

\newcommand{\cP}{\mathscr{P}}

\newcommand{\be}{{\bf e}}

\newcommand{\e}{\epsilon}
\newcommand{\de}{\delta}

\newcommand{\td}{\widetilde}

\newcommand{\La}{\Lambda}

\newcommand{\la}{\lambda}

\newcommand{\hf}{\frac{1}{2}}

\newcommand{\si}{(-1)^{\e_i}}

\newcommand{\ot}{\otimes}

\newcommand{\bq}{{\bf q}}

\usepackage{stmaryrd}
\usepackage{verbatim}
\usepackage{mathtools}
\usepackage[title]{appendix}

\theoremstyle{definition}
\newtheorem{defn}[thm]{\bf Definition}
\theoremstyle{remark}

\newcommand{\aha}{H^{\rm aff}_{\ell}(q^{2})}
\newcommand{\fha}{{H}_{\ell}(q^{2})}
\newcommand{\ket}[1]{\left|#1\right\rangle}

\begin{document}
\title
{Super duality for quantum affine algebras of type $A$}

\author{JAE-HOON KWON}

\address{Department of Mathematical Sciences and RIM, Seoul National University, Seoul 08826, Korea}
\email{jaehoonkw@snu.ac.kr}

\author{SIN-MYUNG LEE}

\address{Department of Mathematical Sciences, Seoul National University, Seoul 08826, Korea}
\email{luckydark@snu.ac.kr}

\thanks{This work is supported by the National Research Foundation of Korea(NRF) grant funded by the Korea government(MSIT) (No.\,2019R1A2C108483311 and 2020R1A5A1016126).}

\begin{abstract}
We introduce a new approach to the study of finite-dimensional representations of the quantum group of the affine Lie superalgebra ${\rm L}{\mf{gl}}_{M|N}=\mathbb{C}[t,t^{-1}]\ot\mf{gl}_{M|N}$ ($M\neq N$). We explain how the representations of the quantum group of ${\rm L}{\mf{gl}}_{M|N}$ are directly related to those of the quantum affine algebra of type $A$, using an exact monoidal functor called truncation. This can be viewed as an affine analogue of super duality of type $A$.
\end{abstract}

\maketitle
\setcounter{tocdepth}{1}
\tableofcontents

\section{Introduction}

\subsection{Lie superalgebra and super duality}
Let $\mf{gl}_{M|N}$ be the general linear Lie superalgebra over $\mathbb{C}$. The character of a finite-dimensional irreducible representation of $\mf{gl}_{M|N}$ has been obtained in \cite{B,BS,CL,Ser} by various independent methods.

One of these methods is the notion of {\em super duality}, which was conjectured in \cite{CWZ}, and proved in \cite{CL} (see also \cite{BS} for another proof). Roughly speaking, the super duality provides a concrete connection between the category of finite dimensional $\mf{gl}_{M|N}$-modules and a maximal parabolic BGG category of $\mf{gl}_{M+N}$, which is induced from an equivalence between these two categories at infinite ranks. This  naturally explains why the theory of finite-dimensional $\mf{gl}_{M|N}$-module is not parallel to that of $\mf{gl}_{M+N}$ at all, but how it is still determined by the Kazhdan-Lusztig theory of $\mf{gl}_{M+N}$. We also refer the reader to \cite{BLW,CLW15} for its generalization to the BGG category of $\mf{gl}_{M|N}$.

\subsection{Quantum affine superalgebra}
Let ${\rm L}{\mf{gl}}_{M|N}=\mathbb{C}[t,t^{-1}]\ot\mf{gl}_{M|N}$ be the affine Lie superalgebra associated to $\mf{gl}_{M|N}$.
In this paper, we introduce a new approach to understanding finite-dimensional representations of the quantum group of ${\rm L}{\mf{gl}}_{M|N}$. In other words, we explain how they are directly related to or determined by finite-dimensional representations of the quantum affine algebra of type $A$, which can be viewed as an affine analogue of super duality of type $A$ due to Cheng-Lam \cite{CL}.
For this purpose, we adopt the notion of a generalized quantum group introduced by Kuniba-Okado-Sergeev \cite{KOS} as the quantum group for ${\rm L}{\mf{gl}}_{M|N}$. 

Let $q$ be an indeterminate. The generalized quantum group $\mathcal{U}(\e)$ of type $A$ is a Hopf algebra over $\Bbbk=\mathbb{Q}(q)$ associated to a sequence $\epsilon=(\e_1,\ldots,\e_n)$  with $n\ge 4$ and $\e_i\in\{0,1\}$, which originates in the study of the three-dimensional Yang-Baxter equation (see \cite{KOS} and references therein). 
It is given by an analogue of Drinfeld-Jimbo's presentation including quantum Serre relation for the odd generator corresponding to $(\e_i,\e_{i+1})=(0,1)$ or $(1,0)$ with respect to cyclic order.
The sequence $\e$ represents the type of non-conjugate Borel subalgebra of $\mf{gl}_{M|N}$ under the Weyl group action of $\mf{gl}_{M|N}$, where $M$ and $N$ are the numbers of $0$ and $1$ in $\e$, respectively.

We assume that $M\neq N$ in this paper.
Note that $\mathcal{U}(\e)$ is different from the previously known quantum affine superalgebras introduced by Yamane \cite{Ya99} (say $U(\e)$ with a Drinfeld-Jimbo type presentation) and by Zhang \cite{Z14} (as a realization by Drinfeld type generators). However, the representations of $\mathcal{U}(\e)$ are still closely connected to those of $U(\e)$. We explain this by constructing an isomorphism of algebras (but not preserving comultiplications) from $\mathcal{U}(\e)$ to $U(\e)$ after taking a semi-direct product with the group algebra of $(\mathbb{Z}/2\mathbb{Z})^n$. Thus we obtain a representation of the one from the other in a natural way. 
We also show that $\mathcal{U}(\e)$ has a non-degenerate bilinear form, which is an analogue of the bilinear form on the usual quantum group, which enables us to construct a universal $R$ matrix (cf.~\cite{Lu93}).

We remark that a systematic study of finite-dimensional representations of $U(\e)$ has been recently done by Zhang \cite{Z14,Z16,Z17}, where he gives
a classification of irreducible modules with respect to Drinfeld realization (when $\e$ is standard) and studies certain properties of tensor product of fundamental representations under $R$ matrix including its irreducibility. However, no direct connection with the representation theory of the quantum affine algebra of type $A$ is given as far as we understand.

\subsection{Fusion construction of irreducible modules} 
Let $\mc{C}(\e)$ be the category of finite-dimensional $\mathcal{U}(\e)$-modules with polynomial weights. 
When $\e$ is homogeneous, that is, $\e_{n|0}=(0,\dots,0)$ or $\e_{0|n}=(1,\dots,1)$, it is the category of finite-dimensional modules over the quantum group of type $A_{n-1}^{(1)}$ with polynomial weights. As a module over the subalgebra $\mathring{\mathcal{U}}(\epsilon)$ corresponding to $\mf{gl}_{M|N}$, any object in $\mc{C}(\e)$ is completely reducible and decomposes into a direct sum of irreducible polynomial representations $V_\e(\la)$ parametrized by $(M|N)$-hook partitions $\lambda$.

For $l\geq 1$, let $\mathcal{W}_{l,\e}(x)$ denote a $q$-analogue of the $l$-supersymmetric power of the natural representation of $\mathcal{U}(\e)$ with spectral parameter $x\in \Bbbk^\times$, which is the $l$-th fundamental representation in case of $\e_{0|n}$. We denote by $\mathcal{W}_{l,\e}(x)_{\rm aff}$ the affinization of $\mathcal{W}_{l,\e}(x)$.

We construct a large family of irreducible $\mathcal{U}(\e)$-modules in $\mc{C}(\e)$ by fusion construction.
More precisely, we show that there exists a well-defined renormalized $R$ matrix in the sense of \cite{KKKO15} on arbitrary tensor product of fundamental representations
\begin{equation}\label{eq:renormal R}
\xymatrixcolsep{2pc}\xymatrixrowsep{3pc}\xymatrix{
{\bf r}_{\mb{l},\e}(\mb{c}) : 
\mathcal{W}_{l_1,\e}(c_1)\ot \cdots \ot\mathcal{W}_{l_t,\e}(c_t) \ar@{->}[r] \ &
\mathcal{W}_{l_t,\e}(c_t)\ot \cdots \ot\mathcal{W}_{l_1,\e}(c_1)},
\end{equation}
for ${\mb l}=(l_1,\dots,l_t)\in \mathbb{Z}_+^t$, and $\mb{c}=(c_1,\dots,c_t)\in (\Bbbk^\times)^t$,
and that the image of ${\bf r}_{\mb{l},\e}(\mb{c})$, say $\mathcal{W}_{\e}(\mb{l},\mb{c})$, is irreducible if it is not zero. The proof is based on the properties of the $R$ matrix on real simple modules developed in  \cite{KKKO15} together with the universal $R$ matrix for $\mathcal{U}(\e)$ and the spectral decomposition of its normalized form on $\mathcal{W}_{l,\e}(x)_{\rm aff} \ot \mathcal{W}_{m,\e}(y)_{\rm aff}$.

\subsection{Truncation functor}

Suppose that $\e'$ is a subsequence of $\e$, say $\e'<\e$. 
Motivated by the super duality for general linear and ortho-symplectic Lie superalgebra due to Cheng-Lam-Wang \cite{CL,CLW,CW}, the first author introduced a functor $\mf{tr}^\e_{\e'} : \mc{C}(\e) \longrightarrow \mc{C}(\e')$, which is induced from a homomorphism from $\mathcal{U}(\e')$ to $\mathcal{U}(\e)$ \cite{KY}. It sends $\mathcal{W}_{l,\e}(x)$ and $V_{\e}(\la)$ to $\mathcal{W}_{l,\e'}(x)$ and $V_{\e'}(\la)$, respectively. Above all, its most important feature is that it is exact and preserves tensor products.

We show that the fusion construction in \eqref{eq:renormal R} is compatible with the truncation, that is, $\mf{tr}^\e_{\e'}({\bf r}_{\mb{l},\e}(\mb{c}))={\bf r}_{\mb{l},\e'}(\mb{c})$ and hence $\mf{tr}^\e_{\e'}\left(\mathcal{W}_{\e}(\mb{l},\mb{c})\right)=\mathcal{W}_{\e'}(\mb{l},\mb{c})$. 
Hence if we take $\e''$ such that $\e_{0|n}<\e''$ and $\e<\e''$, then we have the following:
\begin{equation}\label{eq:triple of truncation}
\xymatrixcolsep{3pc}\xymatrixrowsep{1pc}\xymatrix{
& \mc{C}(\e'')  \ar@{->}_{\mf{tr}^{\e''}_{\e}}[dl]\ar@{->}^{\mf{tr}^{\e''}_{\e_{0|n}}}[dr] &  \\
\mc{C}(\e) & &  \mc{C}(\e_{0|n}) 
}
\end{equation}

In particular, this implies that the composition multiplicity of $\mathcal{W}_{\e}(\mb{l},\mb{c})$ in a standard module, and the branching multiplicity of a polynomial $\ring{\mathcal{U}}(\e)$-module in $\mathcal{W}_{\e}(\mb{l},\mb{c})$ are equal to those in $\mc{C}(\e_{0|n})$ whenever they are non-zero. In this sense, the structure of irreducible $\mathcal{W}_{\e}(\mb{l},\mb{c})$ is completely determined by its counterpart in $\mc{C}(\e_{0|n})$. 

\subsection{Duality}
In order to have a more concrete connection, we consider the full subcategory $\mc{C}_J(\e)$ of $\mc{C}(\e)$, whose composition factors are the irreducible $\mathcal{U}(\e)$-modules appearing as composition factors in a tensor product of $\mathcal{W}_{1,\e}(q^{2j})$ for $j\in \mathbb{Z}$.
Then we construct and apply the generalized quantum affine Schur-Weyl duality functor by Kang-Kashiwara-Kim \cite{KKK} on $\mc{C}_J(\e)$. 

Let $R^{J}\text{-gmod}$ denote the category of finite-dimensional graded modules over the quiver Hecke algebra $R^{J}$ associated to a Dynkin quiver of type $A_\infty$.
Here $J$ denotes the set of vertices of the quiver which is equal to $\mathbb{Z}$ as a set.
As was observed in \cite{KY} (also in \cite{KOS} when $\e$ is standard), the spectral decomposition of the normalized $R$ matrix on $\mathcal{W}_{l,\e}(x)_{\rm aff} \ot \mathcal{W}_{m,\e}(y)_{\rm aff}$ does not depend on the choice of $\e$ up to decomposition into irreducible $\ring{\mathcal{U}}(\e)$-modules. Thus for any $\e$, we can construct an analogue of generalized quantum affine Schur-Weyl duality functor $\mc{F}_\e : R^{J}\text{-gmod} \longrightarrow \mc{C}_J(\e)$ as in \cite{KKK}, which is associated to $S : J \longrightarrow \{\mathcal{W}_{1,\e}(1)\} $ and $X: J \longrightarrow \Bbbk $ with $X(j)=q^{-2j}$. 
The functor $\mc{F}_\e$ is exact and preserves tensor products.
Furthermore, it is compatible with the truncation so that we have  $\mf{tr}^{\e}_{\e'}\circ\mc{F}_\e \cong \mc{F}_{\e'}$ for $\e'<\e$.
We remark that $\mc{F}_{\e}$ can be also defined for other choice of $S$ and $J$ as in \cite{KKK}.

Let $\mc{C}_J^\ell(\e)$ be the full subcategory consisting of objects which are direct sums of polynomial $\ring{\mathcal{U}}(\e)$-modules of degree $\ell$. For the subalgebra $R^J(\ell)$ of $R^J$ of degree $\ell$, let  $R^{J}(\ell)\text{-mod}_0$ be the category of finite-dimensional $R^J(\ell)$-modules on which $x_{k}$ acts nilpotently.
We show that the functor $\mc{F}_\e$ gives an equivalence from $R^{J}(\ell)\text{-mod}_0$ to $\mc{C}_J^\ell(\e)$ if $n > \ell$ by adapting similar arguments for the affine Schur-Weyl duality due to Chari-Pressley \cite{CP} and then using the isomorphism between the quiver Hecke algebra and affine Hecke algebra of type $A$ (after a suitable completion) by Brundan-Kleshchev  \cite{BK} and Rouquier \cite{R}.

If we take $\e''$ and $\e'$ such that $\e_{0|n}<\e'<\e''$ and $\e<\e''$ in \eqref{eq:triple of truncation} and the lengths of $\e', \e''$ are greater than $\ell$, then we obtain
\begin{equation*}\label{eq:quadruple of truncation}
\xymatrixcolsep{3pc}\xymatrixrowsep{1pc}\xymatrix{
& \mc{C}^\ell_J(\e'')  \ar@{->}_{\mf{tr}^{\e''}_{\e}}[dl] \ar@{->}^{\mf{tr}^{\e''}_{\e'}}_{\cong}[r]
& \mc{C}^\ell_J(\e')    \ar@{->}^{\mf{tr}^{\e'}_{\e_{0|n}}}[dr] 
&  \\
\mc{C}^\ell_J(\e) &  &  &  \mc{C}^\ell_J(\e_{0|n})
},
\end{equation*}
where  $\mf{tr}^{\e''}_{\e'}$ is an equivalence.
Now the above relation induces an equivalence between the inverse limit categories by taking limits of the sequences $\e''$ and $\e'$(cf.~\cite{EA}). Let $\e^\infty=(\e_i)_{i\geq 1}$ be an infinite sequence with infinitely many $0$'s and $1$'s such that $\e$ is a subsequence of $\e^\infty$, and let $\ov{\e}^\infty=(1,1,1,\dots)$. Then we have an equivalence of monoidal categories
\begin{equation}\label{eq:quadruple of truncation-2}
\xymatrixcolsep{3pc}\xymatrixrowsep{1pc}\xymatrix{
& \mc{C}_J(\e^\infty)  \ar@{->}_{\mf{tr}_{\e}}[dl] \ar@{->}_{\cong}[r]
& \mc{C}_J(\ov{\e}^\infty)    \ar@{->}^{\mf{tr}_{\e_{0|n}}}[dr] 
&  \\
\mc{C}_J(\e) &  &  &  \mc{C}_J(\e_{0|n})
},
\end{equation}
where $\mc{C}_J(\e^\infty)$ (resp. $\mc{C}_J(\ov{\e}^\infty)$) is the restricted inverse limit category associated to $\e^\infty$ (resp. $\ov{\e}^\infty$), and $\mf{tr}_{\e}$ (resp. $\mf{tr}_{\e_{0|n}}$) denotes the truncation associated to the corresponding subsequence.

The equivalence between $\mc{C}_J(\ov{\e}^\infty)$ and $\mc{C}_J(\un{\e}^\infty)$ induced from \eqref{eq:quadruple of truncation-2} with $\un{\e}^\infty=(0,0,0,\dots)$ gives a kind of duality on the finite-dimensional representations of the quantum affine algebra of type $A_{n-1}^{(1)}$, which looks interesting in itself. This correspondence at the level of Grothendieck rings
can be viewed as an affine analogue of the involution on the ring of symmetric functions sending a Schur function to another one of conjugate shape.

\subsection{Applications}
We present a couple of applications. First, we immediately have that the Grothendieck ring of $\mc{C}_J(\e)$ is a quotient of the polynomial ring generated by the variables corresponding to the fundamental representations in $\mc{C}_J(\e)$. In particular, the polynomial, which yields the isomorphism class of a given irreducible module in $\mc{C}_J(\e)$, is the same as the one in $\mc{C}_J(\e_{0|n})$ for a sufficiently large $n$ (up to truncation of variables).
  
Next, we consider a family of irreducible modules in $\mc{C}_J(\e)$, which corresponds to the Kirillov-Reshetikhin modules over the usual quantum affine algebras. 
Recall that there exists a short exact sequence consisting of tensor products of Kirillov-Reshetikhin modules \cite{Her06,Na03,Na04}, which produces an equation called $T$-system \cite{KNS}. It plays an important role in the study of $q$-characters \cite{FR} and the monoidal categorification of a cluster algebra \cite{HL10,HL16}. Now the diagram \eqref{eq:quadruple of truncation-2} also implies the existence of such a short exact sequence in $\mc{C}_J(\e)$, which yields the same $T$-system.

We expect that our approach \eqref{eq:quadruple of truncation-2} can be used to apply various other properties of $\mc{C}(\e_{0|n})$ to finite-dimensional representations in $\mc{C}(\e)$.

\subsection{Organization}
The paper is organized as follows. 
In Section \ref{sec:GQG}, we recall the definition of $\mathcal{U}(\e)$ and show that the positive (or negative) half of it has a non-degenerate bilinear form. We also construct an isomorphism from $\mathcal{U}(\e)$ to the quantum affine superalgebra in \cite{Ya99} after a suitable extension of both of algebras. 
In Section \ref{sec:R matrix}, we define a renormalized $R$ matrix on a tensor product of fundamental representations and construct a family of irreducible objects in $\mc{C}(\e)$ by fusion construction \eqref{eq:renormal R}.
In Section \ref{sec:irr char}, we show that the fusion construction is compatible with the truncation so that the truncation $\mf{tr}^\e_{\e'}$ sends $\mathcal{W}_{\e}(\mb{l},\mb{c})$ to $\mathcal{W}_{\e'}(\mb{l},\mb{c})$. To consider an object which is a limit of $\mathcal{W}_{\e}(\mb{l},\mb{c})$ with respect to truncation, we introduce the inverse limit category associated to $\{\,\mc{C}(\e^{(k)})\,\}_{k\geq 1}$ with $\e^{(k)}<\e^{(k+1)}$.  
In Section \ref{sec:QSWD}, we verify that an analogue of the generalized quantum affine Schur-Weyl duality functor is well-defined for any $\e$.
In Section \ref{sec:duality}, we prove the equivalence of $\mc{F}_\e$ on $R^{J}(\ell)\text{-mod}_{0}$ for $\e$ with rank $n> \ell$, which implies the equivalence of $\mf{tr}^\e_{\e'}$ for $\e'$ with rank $n'> \ell$.

\vskip 2mm

{\bf Acknowledgement} 
The authors would like to thank M.-H. Kim and E. Park for their kind explanation on the works \cite{KKK,KKKO15}, and J.H.K. also would like to thank M. Okado for helpful discussion on generalized quantum groups and related topics. They also would like to thank the referees for very careful reading of the manuscript and helpful comments.

\section{Generalized quantum groups of type $A$}\label{sec:GQG}

\subsection{Generalized quantum group $\mathcal{U}(\e)$}
In this paper, we assume that $n$ is a positive integer such that $n\geq 4$ and $q$ is an indeterminate. Let $\mathbb{Z}_+$ denote the set of non-negative integers.
We put 
\begin{equation*}
[m]=\frac{q^m-q^{-m}}{q-q^{-1}}\quad (m\in \mathbb{Z}_+).
\end{equation*}
We assume the following notations: 
\begin{itemize}
\item[$\bullet$] $\Bbbk=\mathbb{Q}(q)$, $\Bbbk^{\times}=\Bbbk \setminus \left\{0\right\}$,

\item[$\bullet$] $\e=(\e_1,\cdots,\e_n)$ : a sequence with $\e_i\in \{0,1\}$ ($i=1,\dots, n$),

\item[$\bullet$] $M=|\{\,i\,|\,\e_i=0\,\}|$ and $N=|\{\,i\,|\,\e_i=1\,\}|$,

\item[$\bullet$] $\e_{M|N}$ : a sequence with $\e_1=\dots=\e_M=0$, $\e_{M+1}=\dots=\e_{M+N}=1$ ($M+N=n$),
 
\item[$\bullet$] $\mathbb{I}= \{1<2<\cdots <n\}$ : a linearly ordered set with $\mathbb{Z}_2$-grading $\mathbb{I}=\mathbb{I}_0\cup\mathbb{I}_1$ such that $$\mathbb{I}_0=\{\,i\,|\,\e_i=0\,\},\quad \mathbb{I}_1=\{\,i\,|\,\e_i=1\,\},$$

\item[$\bullet$] $P$ : the free abelian group with a basis $\{\,\delta_i\,|\,i\in \mathbb{I}\,\}$,

\item[$\bullet$] $(\,\cdot\,|\,\cdot\,)$ : a bilinear form on $P$ such that $(\de_i|\de_j)=(-1)^{\e_i}\de_{ij}$ $(i,j\in \mathbb{I})$,

\item[$\bullet$] $P^\vee={\rm Hom}_\mathbb{Z}(P,\mathbb{Z})$ with a basis $\{\,\de^\vee_i\,|\,i\in \mathbb{I}\,\}$ such that $\langle \de_i, \de^\vee_j \rangle =\de_{ij}$ $(i,j\in \mathbb{I})$,

\item[$\bullet$] $I=\{\,0,1,\ldots,n-1\,\}$,

\item[$\bullet$] $\alpha_i=\de_i-\de_{i+1}\in P$ $(i\in I)$,

\item[$\bullet$] $\alpha_i^\vee = \de^\vee_i-(-1)^{\e_i+\e_{i+1}}\de^\vee_{i+1}\in P^\vee$  ($i\in I$),

\item[$\bullet$] $I_{\rm even}=\{\,i\in I\,|\,(\alpha_i|\alpha_i)=\pm 2\,\}$, 
$I_{\rm odd}=\{\,i\in I\,|\,(\alpha_i|\alpha_i)=0\,\}$,

\item[$\bullet$] $q_i=\si q^{\si}$ $(i\in \mathbb{I})$, that is,
\begin{equation*}
q_i=
\begin{cases}
q & \text{if $\e_i=0$},\\
-q^{-1} & \text{if $\e_i=1$},\\
\end{cases} \quad (i\in \mathbb{I}),
\end{equation*}

\item[$\bullet$] ${\bq}(\,\cdot\,,\,\cdot\,)$ : a symmetric biadditive function from $P\times P$ to $\Bbbk^{\times}$ given by
\begin{equation*}
\bq(\mu,\nu) = \prod_{i\in \mathbb{I}}q_i^{\langle\mu,\delta^\vee_i\rangle \langle\nu,\delta^\vee_i\rangle},
\end{equation*}


\item[$\bullet$] $P_{\rm af}=\mathbb{Z}\mb{\delta}_1 \oplus \cdots \oplus \mathbb{Z}\mb{\delta}_{n} \oplus\mathbb{Z}\mb{\delta}$ : a free abelian group of rank $n+1$,

\item[$\bullet$]  $\mb{\alpha}_i=\mb{\delta}_i-\mb{\delta}_{i+1} +\delta_{i0}\mb{\delta}\in P_{\rm af}$ ($i\in I$),

\item[$\bullet$] $Q=\bigoplus_{i\in I}\mathbb{Z}\mb{\alpha}_i\subset P_{\rm af}$, 

\item[$\bullet$] $Q_+=\sum_{i\in I}\mathbb{Z}_+\mb{\alpha}_i$, $Q_-=-Q_+$,

\item[$\bullet$] ${\rm cl} : P_{\rm af} \longrightarrow P$ : the linear map given by ${\rm cl}(\mb{\delta}_i)=\delta_i$ for $i\in \mathbb{I}$ and ${\rm cl}(\mb{\delta})=0$,

\item[$\bullet$] $\iota : P \longrightarrow P_{\rm af}$ : a section of ${\rm cl}$ given by $\iota(\delta_i)=\mb{\delta}_i$ for $i\in \mathbb{I}$.

\end{itemize}

Throughout the paper, we understand the subscript $i\in I$ modulo $n$.
\begin{df}\label{def:U(e)}
{\rm
We define ${\mathcal{U}}(\e)$ to be the associative $\Bbbk$-algebra with $1$ 
generated by $k_\mu, e_i, f_i$ for $\mu\in P$ and $i\in I$ 
satisfying
{\allowdisplaybreaks
\begin{gather}
k_0=1, \quad k_{\mu +\mu'}=k_{\mu}k_{\mu'} \quad (\mu, \mu' \in P),\label{eq:Weyl-rel-1} \\ 
k_\mu e_i k_{-\mu}=\bq(\mu,\alpha_i)e_i,\quad 
k_\mu f_i k_{-\mu}=\bq(\mu,\alpha_i)^{-1}f_i\quad (i\in I, \mu\in P), \label{eq:Weyl-rel-2} \\ 
e_if_j - f_je_i =\delta_{ij}\frac{k_{\alpha_i} - k_{-\alpha_i}}{q-q^{-1}}\quad (i,j\in I),\label{eq:Weyl-rel-3}\\
e_i^2= f_i^2 =0 \quad (i\in I_{\rm odd}),\label{eq:Weyl-rel-4}
\end{gather}
and 
\begin{equation}\label{eq:Serre-rel-1}
\begin{split}
&\ \, e_i e_j -  e_j e_i = f_i f_j -  f_j f_i =0
 \quad \text{($i-j\not\equiv \pm 1\!\!\pmod n$)},\\ 
&
\begin{array}{ll}
e_i^2 e_j- (-1)^{\e_i}[2] e_i e_j e_i + e_j e_i^2= 0\\ 
f_i^2 f_j- (-1)^{\e_i}[2] f_i f_j f_i+f_j f_i^2= 0
\end{array}
\quad \text{($i\in I_{\rm even}$ and $i-j\equiv \pm 1\!\!\pmod n$)}, 
\end{split}
\end{equation}
\begin{equation}\label{eq:Serre-rel-2}
\begin{array}{ll}
  e_{i}e_{i-1}e_{i}e_{i+1}  
- e_{i}e_{i+1}e_{i}e_{i-1} 
+ e_{i+1}e_{i}e_{i-1}e_{i} \\  
\hskip 2cm - e_{i-1}e_{i}e_{i+1}e_{i} 
+ (-1)^{\e_i}[2]e_{i}e_{i-1}e_{i+1}e_{i} =0, \\ 
  f_{i}f_{i-1}f_{i}f_{i+1}  
- f_{i}f_{i+1}f_{i}f_{i-1} 
+ f_{i+1}f_{i}f_{i-1}f_{i}  \\  
\hskip 2cm - f_{i-1}f_{i}f_{i+1}f_{i} 
+ (-1)^{\e_i}[2]f_{i}f_{i-1}f_{i+1}f_{i} =0,
\end{array}\quad \text{($i\in I_{\rm odd}$)}.
\end{equation}}
We call $\mathcal{U}(\e)$ the {\em generalized quantum group of affine type $A$ associated to $\e$} (see \cite{KOS}).
}
\end{df}

Let $\mathcal{U}(\e)^+$ (resp. $\mathcal{U}(\e)^-$) be the subalgebra of $\mathcal{U}(\e)$ generated by $e_i$ (resp. $f_i$) for $i\in I$, and $\mathcal{U}(\e)^0$ be the one generated by $k_\mu$ for $\mu\in P$.
Note that $\mathcal{U}(\e)^\pm$ is naturally graded by $Q_\pm$.

For simplicity, we put $k_i=k_{\alpha_i}$ for $i\in I$. We have 
\begin{gather*}
k_ie_jk_i^{-1}=\bq(\alpha_i,\alpha_j)e_j,\quad 
k_if_jk_i^{-1}=\bq(\alpha_i,\alpha_j)^{-1}f_j,\\
e_if_j - f_je_i =\delta_{ij}\frac{k_i - k_i^{-1}}{q-q^{-1}},
\end{gather*}
for $i,j\in I$.
Note that 
\begin{equation*}
\bq(\alpha_i,\alpha_i)=
\begin{cases}
q^2  & \text{if $i\in I_{\rm even}$ and $(\e_i,\e_{i+1})=(0,0)$},\\
q^{-2}  & \text{if $i\in I_{\rm even}$ and $(\e_i,\e_{i+1})=(1,1)$},\\
-1  & \text{if $i\in I_{\rm odd}$},
\end{cases}
\end{equation*}
for $i\in I$.
If $\e_i=0$ (resp. $\e_i=1$) for all $i\in I$, then the subalgebra of $\mathcal{U}(\e)$ generated by $e_i, f_i, k_i$ for $i\in I$ is isomorphic to $U^{\prime}_q\left(A_{n-1}^{(1)}\right)$ (resp. $U^{\prime}_{-q^{-1}}\left(A_{n-1}^{(1)}\right)$), the quantum affine algebra of type $A_{n-1}^{(1)}$ without derivation, more precisely its quotient by a central element.

There is a Hopf algebra structure on $\mathcal{U}(\e)$, where the comultiplication $\Delta$ is given by 
\begin{equation*}\label{eq:comult-1}
\begin{split}
\Delta(k_\mu)&=k_\mu\otimes k_\mu, \\ 
\Delta(e_i)&= 1\ot e_i + e_i\ot k_i^{-1}, \\
\Delta(f_i)&= f_i\ot 1 + k_i\ot f_i , \\  
\end{split}
\end{equation*}
for $\mu\in P$ and  $i\in I$.

Let $- : \mathcal{U}(\e)\longrightarrow \mathcal{U}(\e)$ be the involution of $\mathbb{Q}$-algebras given by 
\begin{equation*}
\ov{q}=q^{-1}, \quad \ov{e_i}=e_i,\quad \ov{f_i}=f_i,\quad \ov{k_\mu}=k_{\mu}^{-1},
\end{equation*}
for $i\in I$ and $\mu\in P$.
Let ${\e}^c=({\e}^c_1,\dots,{\e}^c_n)$ be given by ${\e}^c_i=1-\e_i$ for $1\leq i\leq n$. Then there is an isomorphism of $\mathbb{Q}$-algebras $\ \sim \ : \mathcal{U}(\e) \longrightarrow \mathcal{U}({\e}^c)$ given by 
\begin{equation}\label{eq:iso tilde}
\td{q}=-q^{-1}, \quad \td{e_i}=e_i,\quad \td{f_i}=f_i,\quad \td{k}_\mu=k_\mu,
\end{equation}
for $i\in I$ and $\mu\in P$.
%

\subsection{Quantum group $\mathbf{U}(\e)$}
Let $'{\bf{f}}(\e)$ be the free associative $\Bbbk$-algebra with $1$ generated by $\theta_i$ for $i\in I$. Then ${'\bf f}(\e)$ is naturally graded by $Q_+$, that is, $'{\bf f}(\e)=\bigoplus_{\beta\in Q_+}{'\bf f}(\e)_\beta$ with $\theta_i\in {'{\bf f}}(\e)_{\mb{\alpha}_i}$ for $i\in I$. We let $|x|={\rm cl}(\beta)$ for $x\in {'\bf f}(\e)_\beta$. 

We regard ${'\bf f}(\e)\ot {'\bf f}(\e)$ as an associative $\Bbbk$-algebra with multiplication 
\begin{equation*}
(x_1\ot x_2)(y_1\ot y_2) = \bq(|x_2|,|y_1|)^{-1}(x_1y_1)\ot (x_2y_2)
\end{equation*}
for homogeneous $x_i, y_i\in {'\bf f}(\e)$ ($i=1,2$).
Let $r: {'\bf f}(\e)\longrightarrow {'\bf f}(\e)\ot {'\bf f}(\e)$ be the unique $\Bbbk$-algebra homomorphism given by $r(\theta_i)=\theta_i\ot 1 + 1\ot \theta_i$ for $i\in I$.

\begin{prop}\label{prop:bilinear form on f}
There exists a unique symmetric bilinear form $(\,,\,)$ on $'\bf f(\e)$ with values in $\Bbbk$ such that 
\begin{itemize}
\item[(1)] $(1,1)=1$,

\item[(2)] $(\theta_i,\theta_j)=\delta_{ij}$ for $i,j\in I$,

\item[(3)] $(x,y'y'') = (r(x),y'\ot y'')$ for $x, y',y''\in {'\bf f}(\e)$,

\item[(4)] $(x'x'',y) = (x'\ot x'',r(y))$ for $x',x'', y\in {'\bf f}(\e)$.
\end{itemize}
\end{prop}
\pf The proof is the same as in \cite[Proposition 1.2.3]{Lu93}.
\qed\vskip 2mm

Let $\mc{I}$ be the radical of $(\,,\,)$ in Proposition \ref{prop:bilinear form on f}, and define 
\begin{equation*}
{\bf f}(\e) = {'\bf f}(\e)/\mc{I}.
\end{equation*}
Note that ${\bf f}(\e)$ is also $Q_+$-graded since $\mc{I}$ is $Q_+$-graded.
For $i\in I$,
let ${}_{i}r : {'\bf f}(\e) \longrightarrow {'\bf f}(\e)$ be the unique $\Bbbk$-linear map such that 
\begin{itemize}
\item[(1)] ${}_{i}r(1)=0$ and ${}_{i}r(\theta_j)=\delta_{ij}$ for $j\in I$,

\item[(2)] ${}_{i}r(xy) = {}_{i}r(x)y + \bq(|x|,\alpha_i)^{-1}x{}_{i}r(y)$ for homogeneous $x,y\in {'\bf f}(\e)$.
\end{itemize}
Similarly, let $r_i : {'\bf f}(\e) \longrightarrow {'\bf f}(\e)$ be the $\Bbbk$-linear map such that 
\begin{itemize}
\item[(1)] $r_i(1)=0$ and $r_i(\theta_j)=\delta_{ij}$  for $j\in I$,

\item[(2)] $r_i(xy) = xr_i(y) + \bq(|y|,\alpha_i)^{-1}r_i(x)y$ for homogeneous $x,y\in {'\bf f}(\e)$.
\end{itemize}
By the same arguments as in \cite[1.2.13]{Lu93}, we have
\begin{equation}\label{eq:adjoint of multiplication}
(\theta_iy,x) = (y, {}_ir(x)),\quad (y\theta_i,x) = (y, r_i(x)),
\end{equation}
for $x,y\in {'\bf f}(\e)$.
Then we can check that  $\mc{I}$ is invariant under ${}_ir$ and $r_i$. We also denote by ${}_ir$ and $r_i$ the induced maps on ${\bf f}(\e)$.

\begin{lem}\label{lem:e^2=0}
For $i\in I_{\rm odd}$, we have ${}_ir(\theta_i^2)=0$ and $\theta_i^2\in \mc{I}$.
\end{lem}
\pf 
We have $\bq(\alpha_i,\alpha_i)=q_iq_{i+1}=-1$ since $\e_i\neq \e_{i+1}$, and hence
\begin{equation*}
{}_ir(\theta_i^2) = \theta_i + \bq(\alpha_i,\alpha_i)\theta_i = \theta_i - \theta_i=0.
\end{equation*}
By \eqref{eq:adjoint of multiplication}, we have
$(\theta_i^2,\theta_i^2)=(\theta_i,{}_ir(\theta_i^2))=0$,
which implies that $\theta_i^2\in \mc{I}$ since $({'\bf f}(\e)_\mu,{'\bf f}(\e)_\nu)=0$ for $\mu\neq\nu$. 
\qed

\begin{prop}\label{prop:Serre for f(e)}
The generators $\theta_i$ $(i\in I)$ of ${\bf f}(\e)$ satisfy the relations \eqref{eq:Serre-rel-1} and \eqref{eq:Serre-rel-2}.
\end{prop}
\pf The proof of \eqref{eq:Serre-rel-1} is the same as in \cite[Proposition 1.4.3]{Lu93}.

Let us prove the relation \eqref{eq:Serre-rel-2} for $\theta_1, \theta_2, \theta_3$ when $\e=(\e_1,\e_2,\e_3,\e_4)$ with $\e_2\neq \e_3$. 
We consider only the case when $\e=(0,0,1,1)$ since the other cases can be proved similarly.

Let
\begin{equation*}
\begin{array}{ll}
S=  
\theta_{2}\theta_{1}\theta_{2}\theta_{3}  
- \theta_{2}\theta_{3}\theta_{2}\theta_{1} 
+ \theta_{3}\theta_{2}\theta_{1}\theta_{2}
 - \theta_{1}\theta_{2}\theta_{3}\theta_{2} 
+ [2]\theta_{2}\theta_{1}\theta_{3}\theta_{2}.
\end{array}
\end{equation*}
We claim that $(S,\theta_{i_1}\theta_{i_2}\theta_{i_3}\theta_{i_4})=0$ for all $(i_1,i_2,i_3,i_4)$ such that $|\theta_{i_1}\theta_{i_2}\theta_{i_3}\theta_{i_4}| = \alpha_1 + 2\alpha_2+\alpha_3$.
Note that
\begin{equation*}
\left( \bq(\alpha_i,\alpha_j) \right)_{1\leq i,j\leq 3} =
\begin{pmatrix}
q^2 & q^{-1} & 1 \\
q^{-1} & -1 & -q \\
1 & -q & q^{-2}  \\
\end{pmatrix}.
\end{equation*}
For example, if $(i_1,i_2,i_3,i_4)=(2,1,2,3)$, then we have
\begin{equation*}
\begin{split}
(\theta_2\theta_1\theta_2\theta_3,\theta_2\theta_1\theta_2\theta_3)&= 1-q^{2},\\
(\theta_2\theta_3\theta_2\theta_1,\theta_2\theta_1\theta_2\theta_3)&= 0,\\
(\theta_3\theta_2\theta_1\theta_2,\theta_2\theta_1\theta_2\theta_3)&= q^{-2}(1-q^{2}),\\
(\theta_1\theta_2\theta_3\theta_2,\theta_2\theta_1\theta_2\theta_3)&= 0,\\
(\theta_2\theta_1\theta_3\theta_2,\theta_2\theta_1\theta_2\theta_3)&= q-q^{-1},\\
\end{split}
\end{equation*}
and hence 
\begin{equation*}
\begin{split}
(S,\theta_2\theta_1\theta_2\theta_3)= (1-q^2) + q^{-2}(1-q^{2}) + (q+q^{-1})(q-q^{-1})=0.
\end{split}
\end{equation*}
The other cases can be also proved in a straightforward manner. Therefore, we conclude that $S \in \mc{I}$.
\qed\vskip 2mm

Let ${'\bf U}(\e)$ be the associative $\Bbbk$-algebra with $1$ generated by ${\rm K}_\mu$, ${\rm E}_i$, ${\rm F}_i$ for $\mu\in P$ and $i\in I$ subject to the following relations:
\begin{gather*}
{\rm K}_0=1, \quad {\rm K}_{\mu}{\rm K}_{\mu'}={\rm K}_{\mu+\mu'} \quad  (\mu, \mu' \in P),  \\ 
{\rm K}_\mu {\rm E}_i {\rm K}_{-\mu}=\bq(\mu,\alpha_i){\rm E}_i,\quad 
  {\rm K}_\mu {\rm F}_i {\rm K}_{-\mu}=\bq(\mu,\alpha_i)^{-1}{\rm F}_i \quad (i\in I, \mu\in P),  \\ 
{\rm E}_i{\rm F}_j - {\rm F}_j{\rm E}_i =\delta_{ij}\frac{{\rm K}_{\alpha_i} - {\rm K}_{-\alpha_i}}{q-q^{-1}}\quad (i,j\in I).
\end{gather*}
Then we define ${\bf U}(\e)$ to be the quotient of ${'\bf U}(\e)$ by the two-sided ideal generated by
$h({\rm E}_0,\dots,{\rm E}_n)$ and $h({\rm F}_0,\dots,{\rm F}_n)$ 
for any $h(\theta_0,\dots,\theta_n)\in\mathcal{I}$.

Let ${\bf U}(\e)^+$ (resp. ${\bf U}(\e)^-$) be the subalgebra of ${\bf U}(\e)$ generated by ${\rm E}_i$ (resp. ${\rm F}_i$) for $i\in I$, and ${\bf U}(\e)^0$ be the one generated by ${\rm K}_\mu$ for $\mu\in P$. By similar arguments as in \cite[3.2]{Lu93}, it is not difficult to see that we have an isomorphism of $\Bbbk$-vector spaces
\begin{equation}\label{0eq:triangular decomp over A}
{\bf U}(\e)^+ \ot {\bf U}(\e)^0 \ot {\bf U}(\e)^- \longrightarrow {\bf U}(\e),
\end{equation}
sending $u\ot {\rm K}_\mu \ot v$ to $u{\rm K}_\mu v$.

Let $\pm : {\bf f}(\e) \longrightarrow {\bf U}(\e)^\pm$ be the homomorphism of $\Bbbk$-algebras given by $\theta_i^+ = {\rm E}_i$ and $\theta_i^-= {\rm F}_i$ for $i\in I$, respectively. By Lemma \ref{lem:e^2=0} and Proposition \ref{prop:Serre for f(e)}, there exists a surjective homomorphism of $\Bbbk$-algebras 
\begin{equation}\label{eq:projection to radical quotient}
\pi : \mathcal{U}(\e) \longrightarrow {\bf U}(\e),
\end{equation}
such that $\pi(k_\mu)={\rm K}_{\mu}$, $\pi(e_i)={\rm E}_i$, $\pi(f_i)={\rm F}_i$ for $\mu \in P$ and $i\in I$.

\subsection{Quantum affine superalgebra}\label{subsec:Yamane quantum group}
Let us first recall the quantized enveloping algebra associated to the affine Lie superalgebra corresponding to the Cartan matrix $A=(\langle \alpha_j, \alpha_i^\vee \rangle)_{i,j\in I}$ \cite{Ya99}.

\begin{df}\label{df:Yamane quantum group}{\rm
Let $U(\e)$ be the associative $\Bbbk$-algebra with $1$ 
generated by $K_\mu, E_i, F_i$ for $\mu\in P$ and $i\in I$ 
satisfying
{\allowdisplaybreaks
\begin{gather}
K_0=1, \quad K_{\mu +\mu'}=K_{\mu}K_{\mu'} \quad  (\mu, \mu' \in P), \label{eq:Weyl-rel-1-Ya}\\ 
 K_\mu E_i K_\mu^{-1}=q^{(\mu|\alpha_i)}E_i,\quad 
 K_\mu F_i K_\mu^{-1}=q^{-(\mu|\alpha_i)}F_i\quad (i\in I, \mu\in P),\label{eq:Weyl-rel-2-Ya} \\ 
 E_iF_j - (-1)^{p(i)p(j)} F_jE_i =
{(-1)^{\e_i}}\delta_{ij}\frac{K_{{\alpha}_i} - K_{-{\alpha}_i}}{q-q^{-1}} \quad (i,j\in I),\label{eq:Weyl-rel-3-Ya}\\
 E_i^2 = F_i^2=0 \quad (i\in I_{\rm odd}), \label{eq:Serre-rel-1-Ya} \\
 E_i E_j - (-1)^{p(i)p(j)} E_j E_i = F_i F_j - (-1)^{p(i)p(j)}  F_j F_i =0,
 \quad (\text{$i,j\in I$ and $i-j\not\equiv 1\!\!\pmod{n}$}), \label{eq:Serre-rel-2-Ya} \\ 
\!\!\!
\begin{array}{ll}
E_i^2 E_j- [2] E_i E_j E_i + E_j E_i^2= 0\\ F_i^2 F_j-[2] F_i F_j F_i+F_j F_i^2= 0
\end{array}
\quad (\text{$i\in I_{\rm even}$ and $i-j\equiv 1\!\!\pmod{n}$}), \label{eq:Serre-rel-3-Ya}  \\
\!\!\!
\begin{array}{ll}
\left[E_i,\left[\left[E_{i-1},E_i\right]_{(-1)^{p(i-1)}q},E_{i+1}\right]_{(-1)^{(p(i-1)+p(i))p(i+1)}q^{-1}}\right]_{(-1)^{p(i-1)+p(i)+p(i+1)}}=0\\
\left[F_i,\left[\left[F_{i-1},F_i\right]_{(-1)^{p(i)}q},F_{i+1}\right]_{(-1)^{(p(i-1)+p(i))p(i+1)}q^{-1}}\right]_{(-1)^{p(i-1)+p(i)+p(i+1)}}=0 \end{array}\ (i\in I_{\rm odd}),\label{eq:Serre-rel-4-Ya} 
\end{gather}
\noindent where
$p(i)=\e_i+\e_{i+1}$ $(i\in I)$, and $[X,Y]_t = XY-tYX$ for $t\in \Bbbk$.}
}
\end{df}

Let $\Sigma$ be the bialgebra over $\Bbbk$ generated by $\sigma_j$ for $j\in \mathbb{I}$, which commute with each other and satisfy $\sigma_j^2=1$.
Here the comultiplication is given by $\Delta(\sigma_j)=\sigma_j\otimes \sigma_j$ for $j\in \mathbb{I}$.
Then $U(\e)$ is a $\Sigma$-module algebra where $\Sigma$ acts on $U(\e)$ by 
\begin{equation}\label{eq:sigma-rel}
\begin{split}
&\sigma_j K_\mu =K_\mu,\quad
\sigma_jE_i=(-1)^{\e_j(\delta_j|\alpha_i)}E_i,\quad 
\sigma_jF_i=(-1)^{\e_j(\delta_j|\alpha_i)}F_i,
\end{split}
\end{equation}
for $j\in \mathbb{I}$, $\mu\in P$ and $i\in I$.
Let $U(\e)[\sigma]$ be the semidirect product of $U(\e)$ and $\Sigma$.

Let 
\begin{equation}\label{eq:extended weight lattice}
P^e_{\rm af}= P_{\rm af} \oplus \mathbb{Z}\La_0 = \bigoplus_{j\in \mathbb{I}}\mathbb{Z}\mb{\delta}_i \oplus \mathbb{Z}\mb{\delta} \oplus \mathbb{Z} \Lambda_0,
\end{equation}
with a symmetric bilinear form 
\begin{equation}\label{eq:extended bilinear form}
(\mb{\delta}_i|\mb{\delta}_j)=(-1)^{\e_i}\delta_{ij},\quad (\mb{\delta}_i|\mb{\delta})=(\mb{\delta}_i|\Lambda_0)=(\Lambda_0|\Lambda_0)=(\mb{\delta}|\mb{\delta})=0,\quad (\Lambda_0|\mb{\delta})=1.
\end{equation}
We define $U(\e)^e_{\rm af}$ to be the $\Bbbk$-algebra defined in the same way as in Definition \ref{df:Yamane quantum group}, with $P$, $(\,\cdot\,|\,\cdot\,)$ and $\alpha_i$ replaced by \eqref{eq:extended weight lattice}, \eqref{eq:extended bilinear form}, and $\mb{\alpha}_i$, respectively.
Let $U(\e)^e_{\rm af}[\sigma]$ be the semidirect product of $U(\e)^e_{\rm af}$ and $\Sigma$ with respect to \eqref{eq:sigma-rel}.

\begin{rem}\label{rem:connection with Yamane's definition}
{\rm
Let $\sigma=\sigma_1\cdots\sigma_n$. 
Then the subalgebra of $U(\e)^e_{\rm af}[\sigma]$ generated by $K_{\mu}$ ($\mu\in Q + \mathbb{Z}\Lambda_0\subset P^e_{\rm af}$), $E_i$, $F_i$ ($i\in I$) and $\sigma$ 
is equal to the quantized enveloping algebra $U_q^\sigma$ \cite[Section 6.4]{Ya99} of the affine Lie superalgebra associated to $A=(\langle \alpha_j, \alpha_i^\vee \rangle)_{i,j\in I}$. 
Note that
\begin{equation}\label{eq:commutation with sigma}
\sigma K_\mu = K_\mu \sigma,\quad \sigma E_i = (-1)^{p(i)} E_i\sigma,\quad \sigma F_i = (-1)^{p(i)} F_i\sigma,
\end{equation}
for $\mu$ and $i\in I$.
Indeed, the algebra $U_q^{\sigma}$ is defined in \cite[Section 6.4]{Ya99} as a quotient of the algebra $\td{U}^\sigma_q$ satisfying the relations \eqref{eq:Weyl-rel-1-Ya}-\eqref{eq:Weyl-rel-3-Ya} and \eqref{eq:commutation with sigma} by the 2-sided ideal generated by the radical of a symmetric bilinear form. 
Then it is shown \cite[Theorem 6.8.2]{Ya99} that $U_q^\sigma$ is isomorphic to the one defined by the relations in Definition \ref{df:Yamane quantum group} and \eqref{eq:commutation with sigma} when $M\neq N$.
}
\end{rem}

We also define $\mathcal{U}(\e)^e_{\rm af}$ to be the $\Bbbk$-algebra defined in the same way as in Definition \ref{def:U(e)} where $P$ and $\bq(\,\cdot\,,\,\cdot\,)$ replaced by $P^e_{\rm af}$ and 
\begin{equation*}
\bq(\mu,\nu)=\bq(\mu_0,\nu_0)q^{ab'+a'b},
\end{equation*}
for $\mu=\mu_0+a\Lambda_0+b\mb{\delta}$, $\nu=\nu_0+a'\Lambda_0+b'\mb{\delta}\in P^e_{\rm af}$ with $\mu_0, \nu_0\in P$.
Then we define $\mathcal{U}(\e)^e_{\rm af}[\sigma]$ to be the semidirect product of $\mathcal{U}(\e)^e_{\rm af}$ and $\Sigma$, where $\Sigma$ acts on $\mathcal{U}(\e)^e_{\rm af}$ by
\begin{equation*}\label{eq:sigma-rel-2}
\begin{split}
&\sigma_j k_\mu =k_\mu,\quad
\sigma_je_i=(-1)^{\e_j(\delta_j|\alpha_i)}e_i,\quad 
\sigma_jf_i=(-1)^{\e_j(\delta_j|\alpha_i)}f_i,
\end{split}
\end{equation*}
for $j\in \mathbb{I}$, $\mu\in P^e_{\rm af}$ and $i\in I$.  

We may assume that $M,N\neq 0$. 
By \eqref{eq:iso tilde}, we may also assume that $\e_1=0$. Let $1=i_0< i_1<\dots<i_l< i_{l+1}=n$ be the unique sequence such that 
\begin{itemize}
\item[(1)] $\e_{i_k}\neq \e_{i_{k+1}}$ for $1\leq k\leq l$,

\item[(2)] $\e_a$ is constant on $\mathbb{I}^{(k)}:=\{\,a\,|\,i_{k-1}+1-\delta_{1k}\leq a\leq i_k\,\}$ for each $1\leq k\leq l+1$.
\end{itemize}
Note that $\mathbb{I}=\mathbb{I}^{(1)}\sqcup \dots \sqcup \mathbb{I}^{(l+1)}$. 
Put  
\begin{gather*}
\sigma_{\leq j} = \sigma_1\sigma_2\cdots\sigma_j\quad 
(j\in \mathbb{I}), \quad  \varsigma_i=\sigma_i\sigma_{i+1}\quad (i\in I).
\end{gather*}
Now, let us define a map $\tau$ on $E_i$, $F_i$, $K_i:=K_{\mb{\alpha}_i}$ for $i\in I$ as follows;

\begin{itemize}
\item[(1)] If  $i\in I_{\rm even}$ with $(\e_i,\e_{i+1})=(0,0)$, then we define
\begin{equation*}
\tau(E_i)=e_i,\quad \tau(F_i)=f_i,\quad \tau(K_i)=k_i.
\end{equation*}

\item[(2)] If $i\in I_{\rm odd}$ with $(\e_i,\e_{i+1})=(0,1)$, then we define
\begin{equation*}
\tau(E_i)=e_i\sigma_{\leq i},\quad \tau(F_i)=f_i\sigma_{\leq i}\varsigma_i,\quad \tau(K_i)=k_i\varsigma_i.
\end{equation*}

\item[(3)] If $i\in I_{\rm even}$ with $(\e_i,\e_{i+1})=(1,1)$, then we have $\{\,i,i+1\,\}\subset \mathbb{I}^{(k)}\cap \mathbb{I}_{1}$ for some $1\leq k\leq l+1$, and define
\begin{equation*}
\tau(E_i)=e_i\varsigma_i^{i-i_{k-1}},\quad \tau(F_i)=f_i\varsigma_i^{i-i_{k-1}-1},\quad \tau(K_i)=-k_i\varsigma_i.
\end{equation*}

\item[(4)] If $i\in I_{\rm odd}\setminus\{0\}$ with $(\e_i,\e_{i+1})=(1,0)$, then we have $i=i_k$ for some $1< k\leq l$, and define
\begin{equation*}
\tau(E_i)=e_i\sigma_{\leq i}\varsigma_i^{i_k-i_{k-1}},\quad \tau(F_i)=f_i\sigma_{\leq i}\varsigma_i^{i_k-i_{k-1}-1},\quad \tau(K_i)=(-1)^{i_k-i_{k-1}}k_i\varsigma_i.
\end{equation*}

\item[(5)] If $i=0$ and $i\in I_{\rm odd}$, that is, $(\e_n,\e_1)=(1,0)$, 
then we define
\begin{equation*}
\tau(E_0)=e_0\sigma_{\leq n}\varsigma_n^{n-i_{l}},\quad \tau(F_0)=f_0\sigma_{\leq n}\varsigma_n^{n-i_{l}-1},\quad \tau(K_0)=(-1)^{n-i_{l}}k_0\varsigma_n.
\end{equation*}
\end{itemize}

\begin{thm}\label{thm:tau iso}
Suppose that $M,N\neq 0$. Let $\ms{X}$ be the $\Bbbk$-subalgebra of $U(\e)^e_{\rm af}[\sigma]$ generated by $E_i, F_i, K_i, K_i^{-1}$ $(i\in I)$, $K_{\pm\Lambda_0}$, and $\sigma_j$ $(j\in \mathbb{I})$. 
Let $\ms{Y}$ be the $\Bbbk$-subalgebra of $\mathcal{U}(\e)^e_{\rm af}[\sigma]$ generated by $e_i, f_i, k_i, k^{-1}_i$ $(i\in I)$, $k_{\pm\Lambda_0}$, and $\sigma_j$ $(j\in \mathbb{I})$. 
Then $\tau$ extends to an isomorphism of $\Bbbk$-algebras
$\tau : \ms{X} \longrightarrow \ms{Y}$
satisfying
\begin{gather*}
\tau(\sigma_j)=\sigma_j\quad (j\in \mathbb{I}),\quad \tau(K_{\Lambda_0})=k_{\Lambda_0}.
\end{gather*}
\end{thm}
\pf It is straightforward to check that $\tau$ yields a well-defined homomorphism of $\Bbbk$-algebras $\tau : \ms{X} \longrightarrow \ms{Y}$. The inverse map is defined in a similar way, and hence it is an isomorphism.
\qed

\begin{rem}{\rm
The isomorphism is a generalization of the one in \cite[Proposition 4.4]{KO} (see also \cite[Section 3.3]{KOS}), where the map is defined only on the subalgebra of finite type generated by $K_\mu$ ($\mu\in P$), $E_i,F_i$ ($i\in I\setminus\{0\}$) and $\sigma_j$ ($j\in \mathbb{I}$), and $\e$ is assumed to be standard, that is, $\e=\e_{M|N}$.
We also remark that the algebra $U(\e)^e_{\rm af}$ is a Hopf algebra \cite{Ya99}, but the isomorphism $\tau$ on $\ms{X}$ does not preserve the comultiplication.
}
\end{rem}

\subsection{Isomorphism from $\mathcal{U}(\e)$ to ${\bf U}(\e)$}

\begin{thm}\label{thm:isomorphism}
Suppose that $M\neq N$. Then the map $\pi : \mathcal{U}(\e) \longrightarrow {\bf U}(\e)$ in
\eqref{eq:projection to radical quotient} is an isomorphism of $\Bbbk$-algebras.
\end{thm}
\pf We may define ${\bf U}(\e)^e_{\rm af}$ in the same way as in $\mathcal{U}(\e)^e_{\rm af}$ (see Section \ref{subsec:Yamane quantum group}).
There is a well-defined action of $\Sigma$ on ${\bf U}(\e)^e_{\rm af}$ given by
\begin{equation*}\label{eq:sigma-rel-3}
\begin{split}
&\sigma_j {\rm K}_\mu ={\rm K}_\mu,\quad
\sigma_j{\rm E}_i=(-1)^{\e_j(\delta_j|\alpha_i)}{\rm E}_i,\quad 
\sigma_j{\rm F}_i=(-1)^{\e_j(\delta_j|\alpha_i)}{\rm F}_i,
\end{split}
\end{equation*}
for $j\in \mathbb{I}$, $\mu\in P^e_{\rm af}$ and $i\in I$, since the action of $\Sigma$ on ${'\bf f}(\e)$ is well-defined and it preserves $\mc I$. 
Then ${\bf U}(\e)^e_{\rm af}[\sigma]$ is defined in a similar way. The projection $\pi$ \eqref{eq:projection to radical quotient} can be extended naturally to
\begin{equation*}\label{eq:projection to radical quotient-2}
\pi : \mathcal{U}(\e)^e_{\rm af}[\sigma] \longrightarrow {\bf U}(\e)^e_{\rm af}[\sigma],
\end{equation*}
where $\pi(\sigma_j)=\sigma_j$ for $j\in \mathbb{I}$. 
Let $\ms{Z}$ be the $\Bbbk$-subalgebra of ${\bf U}(\e)^e_{\rm af}[\sigma]$ generated by ${\rm E}_i, {\rm F}_i, {\rm K}_i, {\rm K}^{-1}_i$ $(i\in I)$, ${\rm K}_{\pm\Lambda_0}$, and $\sigma_j$ $(j\in \mathbb{I})$, where ${\rm K}_i:={\rm K}_{\mb{\alpha}_i}$. 
Consider
\begin{equation*}
\xymatrixcolsep{2pc}\xymatrixrowsep{3pc}\xymatrix{
\ms{X} \ \ar@{->}^{\tau}[r] &\ \ms{Y} \ \ar@{->}^{\pi}[r] & \ms{Z}.}
\end{equation*}

Let $\ms{X}^+$ be the subalgebra of $\ms{X}$ generated by $E_i$ for $i\in I$ and 
let $\ms{X}^0$ be the subalgebra generated by $\sigma_j$ ($j\in \mathbb{I}$), $K_i, K_i^{-1}$ ($i\in I$) and $K_{\pm\Lambda_0}$. 
We define $\ms{Y}^+, \ms{Y}^0$ and $\ms{Z}^+, \ms{Z}^0$ similarly. Recall that $\ms{X}^+, \ms{Y}^+, \ms{Z}^+$ are $Q_+$-graded. For $\beta=\sum_{i\in I}c_i\mb{\alpha}_i\in Q_+$, let ${\rm ht}(\beta)=\sum_{i\in I}c_i$.

Suppose that $y \in {\rm Ker}\,\pi \cap \ms{Y}^+ \setminus \{0\}$ is given. 
Since $(\ms{Z}^+_\beta,\ms{Z}^+_{\beta'})=0$ for $\beta, \beta'\in Q_+$ ($\beta\neq \beta'$), we may assume that $y$ is homogeneous, that is, $y\in \ms{Y}^+_\beta$ for some $\beta\in Q_+$. We further assume that ${\rm ht}(\beta)$ is minimal such that ${\rm Ker}\pi\cap\mathscr{Y}^{+}_{\beta}\neq 0$.

By Theorem \ref{thm:tau iso}, there exist a unique non-zero $X\in \ms{X}^{+}_{\beta}$ and a monomial $\varsigma$ in $\sigma_i$ such that $\tau(X\varsigma)=y$. Since ${\rm Ker}\, \pi\circ \tau $ is a two-sided ideal and $\varsigma$ is invertible, we have $X\in {\rm Ker}\, \pi\circ \tau$ as well.

Consider $[X,F_i]:=XF_i - (-1)^{p(i)p(\beta)}F_iX$ for $i\in I$, where $p(\beta)=\sum_{i\in I_{\rm odd}}c_i$. Note that $[X,F_i]\in {\rm Ker}\,\pi\circ\tau$.
By applying \eqref{eq:Weyl-rel-3-Ya}, we have
\begin{equation}\label{eq:super q-commutator}
[X,F_i] = X'K_i + X''K_i^{-1} 
\end{equation}
for some $X', X''\in \ms{X}^+_{\beta-\mb{\alpha}_i}$. 
Write $\tau(X')=Y'\sigma'$ and $\tau(X'')=Y''\sigma''$ for some $Y', Y''\in \ms{Y}^+$ and $\sigma', \sigma''\in \Sigma$. 
Let $Z'=\pi(Y')$ and $Z''=\pi(Y'')$. 
Applying $\pi\circ \tau$ to \eqref{eq:super q-commutator}, we get
\begin{equation*} 
Z'{\rm K}_{i}\varsigma' + Z''{\rm K}_{i}^{-1}\varsigma''=0,
\end{equation*} 
for some $\varsigma', \varsigma''\in \Sigma$.
Since the vectors $Z'{\rm K}_{i}\varsigma', Z''{\rm K}_{i}^{-1}\varsigma''$ are linearly independent by \eqref{0eq:triangular decomp over A}, we have $Z'=Z''=0$. Hence $Y',\, Y''=0$ by the minimality of ${\rm ht}(\beta)$, which implies $X'=X''=0$ and $[X,F_i]=0$.
 
Finally  we have $X=0$ by \cite[Proposition 6.5.1]{Ya99}, and then $y=0$, which contradicts the fact that $y$ is non-zero.
Therefore we conclude that $\pi\vert_{\ms{Y}^+}$ is injective (and the same holds for $\ms{Y}^-$, where $\ms{Y}^-$ is the subalgebra generated by $f_i$ $(i\in I)$). Since $\pi\vert_{\ms{Y}^0}$ is an isomorphism, so is $\pi$ by \eqref{0eq:triangular decomp over A}.
\qed

\section{$R$ matrix}\label{sec:R matrix}

\subsection{Finite-dimensional modules and affinization}
Let $V$ be a $\mathcal{U}(\e)$-module, which is a $\Bbbk$-vector space. For $\la\in P$, let 
\begin{equation*}
V_\la 
= \{\,u\in V\,|\,k_{\mu} u= \bq(\la,\mu) u \ \ (\mu\in P) \,\}
\end{equation*}
be the $\la$-weight space of $V$. We have
\begin{equation*}\label{eq:action of e and f on weight space}
e_iV_\la \subset V_{\la+\alpha_i},\quad 
f_iV_\la \subset V_{\la-\alpha_i}\quad (i\in I).
\end{equation*}
Put ${\rm wt}(V)=\{\,\mu\in P\,|\,V_\mu\neq 0\,\}$. 
Let $P_{\geq0}=\sum_{i\in \mathbb{I}}\mathbb{Z}_+\delta _i$. For $\la\in P_{\ge 0}$ with $\la=\sum_{i}\la_i\de_i$, let ${\rm deg}(\la)=\sum_{i}\la_i$.
Let $\ring{\mathcal{U}}(\e)$ be the $\Bbbk$-subalgebra of $\mathcal{U}(\e)$ generated by $k_\mu$, $e_i$ and $f_i$ for $\mu\in P$ and $i\in I\setminus \{0\}$. 

\begin{df}{\rm \mbox{}
\begin{itemize}
\item[(1)]
Define $\mc{C}(\e)$ to be the category of finite-dimensional $\mathcal{U}(\e)$-modules $V$ such that $V$ has a weight space decomposition with ${\rm wt}(V)\subset P_{\ge 0}$, that is,
$V=\bigoplus_{\la\in P_{\ge 0}}V_\la$.

\item[(2)]
For $\ell\in \mathbb{Z}_+$, define $\mc{C}^\ell(\e)$ to be the full subcategory of $\mc{C}(\e)$ consisting of $V$ such that ${\rm deg}(\la)=\ell$ for all $\la\in {\rm wt}(V)$.

\item[(3)] Define $\ring{\mc{C}}(\e)$ to be the category of finite-dimensional $\ring{\mathcal{U}}(\e)$-modules $V$ such that $V$ has a weight space decomposition with ${\rm wt}(V)\subset P_{\ge 0}$. The subcategory $\ring{\mc{C}}^\ell(\e)$ is defined in the same way.

\end{itemize}
}
\end{df}

Note that $\mc{C}(\e)$ is closed under taking submodules, quotients and tensor products, while $\mc{C}^\ell(\e)$ is closed under taking submodules and quotients. We also have
\begin{equation*}
\mc{C}(\e)=\bigoplus_{\ell\in \mathbb{Z}_+}\mc{C}^\ell(\e).
\end{equation*}

\begin{rem}{\rm
Our notion of weight space is slightly different from that of usual quantum affine algebras. However, if we restrict the action of $\mathcal{U}(\e)^0$ to that of subalgebra corresponding to $\mathfrak{sl}(M|N)$, then the finite-dimensional representations of $\mathcal{U}(\e)$ and $U(\e)$ are directly related by the isomorphism $\tau$ in Theorem \ref{thm:tau iso} as follows. 

Let $X$ (resp. $Y$) be the subalgebra of $U(\epsilon)[\sigma]$ (resp. $\mathcal{U}(\epsilon)[\sigma]$) generated by $E_i,\,F_i,\,K_i, K^{-1}_i$ ($i\in I$) and $\sigma_j$ ($j\in \mathbb{I}$) (resp. $e_i,\,f_i,\,k_i, k^{-1}_i,\,\sigma_j$). Then the isomorphism in Theorem \ref{thm:tau iso} restricts to $\tau : X \rightarrow Y$
and induces an equivalence between two categories of modules over $X$ and $Y$. Let $V^\tau$ denote the $X$-module obtained from a $Y$-module $V$ by $\tau$.

Let $V$ be a $\mathcal{U}(\epsilon)$-module with weight space decomposition $V=\bigoplus_{\la\in P} V_\lambda$.
Then we have the weight space decomposition as a $Y$-module $V=\bigoplus_{\overline{\lambda}\in\overline{P}}V_{\overline{\lambda}}$, where
\begin{align*}
\overline{P} & \coloneqq P/\mathbb{Z}(\delta_1 +\delta_2 +\cdots +\delta_n ), \\
V_{\overline{\lambda}} & \coloneqq \bigoplus_{\overline{\mu}=\overline{\lambda}} V_\mu =
\left\{u\in V \,|\, k_i u=\mathbf{q}(\alpha_i,\lambda)u \ \   (i\in I) \right\} \quad (\la\in P).
\end{align*}
One can extend the $\mathcal{U}(\epsilon)$-action to $\mathcal{U}(\epsilon)[\sigma]$ by 
$\sigma_{i}u = (-1)^{\epsilon_i(\delta_i | \lambda)}u$ 
for $i\in\mathbb{I}$ and $u\in V_\lambda$. 
On the other hand, for a ${U}(\epsilon)$-module $W$ with the weight space decomposition $W=\bigoplus_{\la\in P} W_\lambda$, where
\[
W_\lambda = \left\{u\in W \,|\, K_\mu u = q^{(\mu|\lambda)}u  \ \ (\mu\in P) \right\},
\]
we have a similar weight space decomposition as an $X$-module, and extend $W$ to a $U(\epsilon)[\sigma]$-module in the same way.
Then the decomposition $V^\tau=\bigoplus_{\overline{\lambda}\in\overline{P}} V_{\overline{\lambda}}$ coincides with the weight space decomposition as an $X$-module, that is,
\[
V_{\ov{\lambda}} = \left\{u\in V^\tau \,|\, K_i u=q^{(\alpha_i|\lambda)}u \ \  (i\in I)\right\}.
\]
}
\end{rem}

Let us introduce an important family of $\mathcal{U}(\e)$-modules in $\mc{C}(\e)$.
Let
\begin{equation*}
\mathbb{Z}^n_+(\e)=\{\,{\bf m}=(m_1,\ldots,m_n)\,|\,\e_i=0 \Rightarrow m_i\in\mathbb{Z}_+,\  \e_i=1 \Rightarrow m_i\in \{0,1\}\,\}.
\end{equation*}
For ${\bf m}\in \mathbb{Z}_+^n(\e)$, let $|{\bf m}|=m_1+\dots+m_n$.
For $i\in \mathbb{I}$, put $\be_i=(0,\cdots, 1,\cdots, 0)$ where $1$ appears only in the $i$-th component. 
For $l\in \mathbb{Z}_{+}$, let 
\begin{equation*}
\mathcal{W}_{l,\e} = \bigoplus_{\substack{{\bf m}\in \mathbb{Z}^n_+(\e), |{\bf m}|=l}}\Bbbk|{\bf m}\rangle
\end{equation*}
be the $\Bbbk$-vector space spanned by $|{\bf m}\rangle$ for ${\bf m}\in \mathbb{Z}^n_+(\e)$ with $|{\bf m}|=l$.

For $x\in \Bbbk^\times$, we denote by $\mathcal{W}_{l,\e}(x)$ a $\mathcal{U}(\e)$-module $V$, where $V=\mathcal{W}_{l,\e}$ as a $\Bbbk$-space and the actions of $k_\mu, e_i, f_i$ are given by  
{\allowdisplaybreaks
\begin{align*}\label{eq:W_l}
k_\mu |{\bf m}\rangle &= \bq\left(\mu,\, \sum_{j\in \mathbb{I}}m_i\delta_j\right) |{\bf m} \rangle,\\
e_i |{\bf m}\rangle &= 
\begin{cases}
x^{\delta_{i,0}}[m_{i+1}]|{\bf m} + \be_{i} -\be_{i+1} \rangle  & \text{if ${\bf m} + \be_{i} -\be_{i+1}\in\mathbb{Z}_+^n(\e)$},\\
0& \text{otherwise},
\end{cases}
\\
f_i |{\bf m}\rangle &= 
\begin{cases}
x^{-\delta_{i,0}}[m_{i}]|{\bf m} - \be_{i} + \be_{i+1} \rangle  & \text{if ${\bf m} - \be_{i} + \be_{i+1}\in\mathbb{Z}_+^n(\e)$},\\
0 & \text{otherwise},
\end{cases}\\
\end{align*}}
for $\mu\in P$, $i\in I$ and ${\bf m}=(m_1,\ldots,m_n)\in \mathbb{Z}^n_+(\e)$. 
Here we understand $\be_0=\be_n$. It is clear that $\mathcal{W}_{l,\e}(x)\in\mc{C}(\e)$. 
We regard $\mathcal{W}_{l,\e}=\mathcal{W}_{l,\e}(1)$.
We call $\mathcal{W}_{l,\epsilon}(x) $ a \emph{fundamental representation} with spectral parameter $x$.

\begin{rem}\label{rem:homogeneous fundamental}
{\rm
If $\e=\e_{0|n}$ and $1\leq l<n$, then $\mathcal{W}_{l,\e}(x)$ is the $l$-th fundamental representation of $U_{-q^{-1}}^{\prime}\left(A_{n-1}^{(1)}\right)$.
If $\e=\e_{n|0}$ and $\ell\geq 1$, then $\mathcal{W}_{l,\e}(x)$ is the Kirillov-Reshetikhin module corresponding to the partition $(l)$.}
\end{rem}

\begin{rem}{\rm
The map $\phi(|{\bf m}\rangle)= q^{\hf\sum_i m_i^2 -\sum_{i<j}m_im_j}|{\bf m}\rangle$ gives an isomorphism of $\mathcal{U}(\e)$-modules from $\mathcal{W}_{l,\e}(x)$ to itself with another $\mathcal{U}(\e)$-action defined in \cite[(2.15)]{KO}, where a different comultiplication is used.
Note that the $A_0$-span of $|\bf{m}\rangle$ is a crystal lattice of $\mathcal{W}_{l,\e}$ \cite{KY}, where $A_0$ is the subring of $f(q)\in \Bbbk$ which are regular at $q=0$.}
\end{rem}\label{rem:fund-twist}

Let $\cP$ be the set of partitions. 
A partition $\la=(\la_i)_{i\ge 1}\in \cP$ is called an $(M|N)$-hook partition if $\la_{M+1}\leq N$. We denote the set of all $(M|N)$-hook partitions by $\cP_{M|N}$.
For $\la\in\cP_{M|N}$, let $V_\e(\la)$ be the irreducible highest weight $\ring{\mathcal{U}}(\e)$-module with the highest weight $\sum_{i\in \mathbb{I}} m_i\delta_i$. Here $m_i$ denotes the number of $i$'s in the tableau $H_{\lambda}$ of shape $\lambda$, which is defined inductively as follows (see also \cite[(2.37)]{CW}):
\begin{enumerate}
    \item If $\epsilon_1 =0$ (resp. $\e_1=1$), then fill the first row (resp. column) of $\lambda$ with $1$.
    \item After filling a subdiagram $\mu$ of $\lambda$ with $1,\dots,k$, fill the first row (resp. column) of $\lambda/\mu$ with $k+1$ if $\epsilon_{k+1} =0$ (resp. $\e_{k+1}=1$).
\end{enumerate}
If we put $V=V_\e((1))$, then $V^{\otimes \ell}$ is semisimple and its irreducible components are $V_\e(\la)$'s for $(M|N)$-hook partitions $\la$ of $\ell$ \cite{BKK,KY}. In particular, $V_\e(\la)\in \ring{\mc{C}}(\e)$.
Note that $\mathcal{W}_{l,\e}(x)\cong V_\e((l))$ for $l\geq 1$ as a $\ring{\mathcal{U}}(\e)$-module. 
Hence, any tensor product of a finite number of $\mathcal{W}_{l,\e}(x)$'s is completely reducible as a $\ring{\mathcal{U}}(\e)$-module, and decomposes into $V_\e(\la)$'s for $\la\in\cP_{M|N}$.

Let $V$ be a $\mathcal{U}(\e)$-module in $\mc{C}(\e)$ and $z$ an indeterminate. We define a $\mathcal{U}(\e)$-module $V_{\rm aff}$ by
\begin{equation*}\label{eq:affinization}
V_{\rm aff} = \Bbbk[z,z^{-1}]\ot V, 
\end{equation*}
where $k_\mu$, $e_i$, and $f_i$ acts as $1\ot k_\mu$, $z^{\delta_{i,0}}\ot e_i$, and $z^{-\delta_{i,0}}\ot f_i$, respectively.
For $\la\in P_{\rm af}$, we define the $\la$-weight space of $V_{\rm aff}$ by
\begin{equation*}
(V_{\rm aff})_{\la}= z^k\ot V_{{\rm cl}(\la)},
\end{equation*}
where $\la- \iota \circ {\rm cl}(\la)=k\mb{\delta}$. Then we have
\begin{equation*}
V_{\rm aff}=\bigoplus_{\la\in P_{\rm af}}(V_{\rm aff})_\la,
\end{equation*}
and $e_i(V_{\rm aff})_\la \subset (V_{\rm aff})_{\la+\mb{\alpha}_i}$, 
$f_i(V_{\rm aff})_\la \subset (V_{\rm aff})_{\la-\mb{\alpha}_i}$ for $i\in I$ and $\la\in P_{\rm af}$.
The map sending $g(z)\ot m$ to $zg(z)\ot m$ for $m\in V$ and $g(z)\in \Bbbk[z,z^{-1}]$ gives an isomorphism of $\mathcal{U}(\e)$-modules 
\begin{equation}\label{eq:automorphism z}
z : V_{\rm aff} \longrightarrow V_{\rm aff},
\end{equation}
such that $z(V_{\rm aff})_\la \subset (V_{\rm aff})_{\la+\mb{\delta}}$ for $\la\in P_{\rm af}$.  
For $x\in \Bbbk^\times$, we define
\begin{equation*}
V_x = V_{\rm aff}/(z-x)V_{\rm aff}.
\end{equation*}
Note that $\mathcal{W}_{l,\e}(x)\cong (\mathcal{W}_{l,\e}(1))_x\cong (\mathcal{W}_{l,\e}(x))_1$ for $x\in \Bbbk^\times$ as a $\mathcal{U}(\e)$-module.

\subsection{Universal $R$ matrix}
In this subsection, we assume that $M\neq N$.
Let $\ov{\Delta}$ be the comultiplication on $\mathcal{U}(\e)$ given by $\ov{\Delta}(u)=\ov{\Delta(\ov{u})}$ for $u\in \mathcal{U}(\e)$, that is,
\begin{equation*}\label{eq:comult-2}
\begin{split}
\ov{\Delta}(k_\mu)&=k_\mu\otimes k_\mu, \\ 
\ov{\Delta}(e_i)&= 1\ot e_i + e_i\ot k_i, \\
\ov{\Delta}(f_i)&= f_i\ot 1 + k_i^{-1}\ot f_i , \\  
\end{split}
\end{equation*}
for $\mu\in P$ and $i\in I$.
Let 
\begin{equation*}
\mathcal{U}(\e)^+ \widehat{\ot}\ \mathcal{U}(\e)^- = 
\bigoplus_{\xi\in Q}\prod_{\xi=\mu+\nu}
\left(\mathcal{U}(\e)^+_\mu \ot \mathcal{U}(\e)^-_\nu\right),
\end{equation*}
which is a ring under the multiplications induced from $\mathcal{U}(\e)^\pm$. 
Similarly, let $\mathcal{U}(\e)^- \td{\ot}\ \mathcal{U}(\e)^+$ be the one defined by replacing $x\ot y\in \mathcal{U}(\e)^+ \widehat{\ot}\ \mathcal{U}(\e)^-$ with $y\ot x$.

For $\beta\in Q_+$, let ${\bf B}_\beta$ be a basis of ${\bf f}(\epsilon)_\beta$ and let ${\bf B}^*_{\beta}=\{\,b^*\,|\,b\in {\bf B}_\beta\,\}$ be the dual basis of ${\bf B}_\beta$, that is, $(b^*,b')=\delta_{b,b'}$ for $b,b'\in {\bf B}_\beta$. We may identify $\mathcal{U}(\epsilon)$ with $\mathbf{U}(\epsilon)$ by Theorem~\ref{thm:isomorphism}.
By slightly modifying the arguments in the proof of \cite[Theorem 4.1.2]{Lu93}, we have the following analogue.
\begin{thm}\label{thm:quasi R matrix}\mbox{}
\begin{itemize}
\item[(1)]
There is a unique $\Theta_\beta\in \mathcal{U}(\e)^+_\beta\ot\mathcal{U}(\e)^-_{-\beta}$ for each $\beta\in Q_+$ such that $\Theta_0=1\ot 1$ and 
\begin{equation*}
\Theta\, \Delta(u) = \ov{\Delta}(u)\, \Theta, 
\end{equation*}
for $u\in \mathcal{U}(\e)$, where $\Theta = \sum_{\beta\in Q_+}\Theta_\beta\in \mathcal{U}(\e)^+\widehat{\ot}\ \mathcal{U}(\e)^-$.  
\item[(2)] For $\beta=\sum_{i\in I}a_i \mb{\alpha}_i\in Q_+$, we have
\begin{equation*}
\Theta_\beta = (q-q^{-1})^{{\rm ht}(\beta)}\sum_{b\in {\bf B}_\beta}b^+\ot (b^*)^-,
\end{equation*}
where ${\rm ht}(\beta)=\sum_{i\in I}a_i$.

\item[(3)] Let $\ov{\Theta}=\sum_{\beta\in Q_+}\ov{\Theta}_\beta$, where $\ov{\Theta}_\beta:= (- \ot -)(\Theta_\beta)$. Then 
\begin{equation*}
\ov{\Theta}{\Theta} ={\Theta}\ov{\Theta}=1.
\end{equation*}
\end{itemize}
\end{thm}

\begin{rem}\label{rem:properties of Theta}
{\rm
We can prove that $\Theta$ also satisfies the properties in  \cite[4.2]{Lu93} by slight modification of the arguments.
}
\end{rem}

Let $V, W$ be $\mathcal{U}(\e)$-modules in $\mc{C}(\e)$.
Let $z_1,z_2$ be indeterminates and let
\begin{equation*}
V_{\rm aff} = \Bbbk[z_1,z_1^{-1}]\ot V,\quad 
W_{\rm aff} = \Bbbk[z_2,z_2^{-1}]\ot W.
\end{equation*}
We consider two completions given in \cite[Section 7]{Kas02}:
\begin{equation*}
\begin{split}
V_{\rm aff}\, \widehat{\ot}\, W_{\rm aff} &= 
\sum_{\la,\mu\in P_{\rm af}} \prod_{\beta\in Q_+} 
(V_{\rm aff})_{\la+\beta}{\ot} (W_{\rm aff})_{\mu-\beta} ,\\
W_{\rm aff}\, \td{\ot}\, V_{\rm aff} &=
\sum_{\la,\mu\in P_{\rm af}} \prod_{\beta\in Q_+} 
(W_{\rm aff})_{\la-\beta}{\ot} (V_{\rm aff})_{\mu+\beta},
\end{split}
\end{equation*}
which are invariant under the action of $\mathcal{U}(\e)^+ \widehat{\ot}\ \mathcal{U}(\e)^-$ and $\mathcal{U}(\e)^- \td{\ot}\ \mathcal{U}(\e)^+$, respectively.
Let $\Pi_{\bq} : V_{\rm aff}\, \widehat{\ot}\, W_{\rm aff} \longrightarrow V_{\rm aff}\, \widehat{\ot}\, W_{\rm aff}$ be given by 
\begin{equation*}
\Pi_{\bq}(v\ot w)=\bq({\rm cl}(\mu),{\rm cl}(\nu))\, v\ot w,
\end{equation*}
for $v\in (V_{\rm aff})_\mu$ and $w\in (W_{\rm aff})_\nu$. 
Let $s :V_{\rm aff}\, \widehat{\ot}\, W_{\rm aff} \longrightarrow W_{\rm aff}\, \td{\ot}\, V_{\rm aff}$ be the map given by $s(v\ot w)=w\ot v$.
Then 
\begin{thm}\label{thm:univ R matrix}
We have an isomorphism of $\mathcal{U}(\e)$-modules
\begin{equation*}
\xymatrixcolsep{2pc}\xymatrixrowsep{3pc}\xymatrix{
\mc{R}^{\rm univ}:= \Theta \circ \Pi_{\bq} \circ s : 
V_{\rm aff}\, \widehat{\ot}\, W_{\rm aff}\ \ar@{->}[r] & \ W_{\rm aff}\, \td{\ot}\, V_{\rm aff}
}
\end{equation*}
\end{thm}
\pf The proof is almost the same as the one in \cite[Theorem 32.1.5]{Lu93}, where we replace ${}_f\Pi$ in \cite{Lu93} with $\Pi_{\bq}$. The inverse is given by $s\circ \Pi_{\bq}^{-1}\circ \ov{\Theta}$.
\qed\vskip 2mm

The operator $\mc{R}^{\rm univ}$ is called the {\em universal $R$ matrix}.
Since ${\rm wt}(V)$ and ${\rm wt}(W)$ are finite subsets of $P_{\geq 0}$, we have as $\mathcal{U}(\e)$-modules
\begin{equation*}
\begin{split}
V_{\rm aff}\, \widehat{\ot}\, W_{\rm aff}&= 
\Bbbk[[z_1/z_2]]\ot_{\Bbbk[z_1/z_2]}\left(V_{\rm aff}\,{\ot}\, W_{\rm aff}\right),\\
W_{\rm aff}\, \td{\ot}\, V_{\rm aff}&=
\Bbbk[[z_1/z_2]]\ot_{\Bbbk[z_1/z_2]}\left(W_{\rm aff}\,{\ot}\, V_{\rm aff}\right).
\end{split}
\end{equation*}
In particular, we have a $\mathcal{U}(\e)$-linear map
\begin{equation}\label{eq:R univ}
\xymatrixcolsep{2pc}\xymatrixrowsep{3pc}\xymatrix{
\mc{R}^{\rm univ}: 
V_{\rm aff} \ot W_{\rm aff}\ \ar@{->}[r] & \ 
\Bbbk[[z_1/z_2]]\ot_{\Bbbk[z_1/z_2]}\left(W_{\rm aff}\,{\ot}\, V_{\rm aff}\right). 
}
\end{equation}

\subsection{Truncation}\label{sec:truncation}

Let $\e'=(\e'_1,\dots,\e'_{n-1})$ be the sequence obtained from $\e$ by removing $\e_i$ for some $i\in \mathbb{I}$. We further assume that 
\begin{equation*}\label{eq:condition on e'}
\text{$\e'$ is homogeneous when it is of length $3$, that is, $\e'=(000)$ or $(111)$.}
\end{equation*}
Let $\mathbb{I}'=\{\,1,\dots,n-1\,\}$ with the $\mathbb{Z}_2$-grading given by $\mathbb{I}^\prime_\varepsilon=\{\,i\,|\,\epsilon^\prime_i =\varepsilon\,\}$ ($\varepsilon=0,1$).
We assume that the weight lattice for $\mathcal{U}(\e')$ is $P'=\bigoplus_{l\in \mathbb{I}'}\mathbb{Z}\delta'_l$.
Put $I'=\{0,1,\cdots,n-2\}.$ 
Let us denote by $k'_{\mu}$, $e'_j$, and $f'_j$ the generators of $\mathcal{U}(\e')$ for $\mu\in P'$ and $j\in I'$. Put $k'_j=k'_{\delta'_j-\delta'_{j+1}}$ for $j\in I'$.
Then we define $\hat{k}_{\delta'_l}, \hat{e}_j, \hat{f}_j\in \mathcal{U}(\e)$ for $l\in \mathbb{I}'$ and $j \in I'$ as follows: 
\begin{equation*}
\hat{k}_{\delta'_l}=
\begin{cases}
k_{\delta_l} & \text{for $1\leq l\leq i-1$}, \\
k_{\delta_{l+1}} & \text{for $i \leq l\leq n-1$.}
\end{cases}
\end{equation*}

{\allowdisplaybreaks
{\em Case 1}. If $2\leq i\leq n-1$, then
\begin{gather*}\label{eq:generators for e'}
(\hat{e}_j, \hat{f}_j) =
\begin{cases}
\left(e_j,f_j\right) & \text{for $j\leq i-2$},\\
\left([e_{i-1},e_{i}]_{\bq_{i-1,i}}, [f_{i},f_{i-1}]_{\bq_{i-1,i}^{-1}}\right) & \text{for $j= i-1$},\\
\left(e_{j+1},f_{j+1}\right) & \text{for $j\geq i$}.
\end{cases}
\end{gather*}

{\em Case 2}. If $i=n$, then  
\begin{gather*}\label{eq:generators for e'-2}
(\hat{e}_j, \hat{f}_j) =
\begin{cases}
(e_j,f_j) & \text{for $j\neq 0$},\\
\left([e_{n-1},e_{0}]_{\bq_{n-1,0}},[f_{0},f_{n-1}]_{\bq_{n-1,0}^{-1}}\right) & \text{for $j=0$}.
\end{cases}
\end{gather*}

{\em Case 3}. If $i=1$, then
\begin{gather*}\label{eq:generators for e'-3}
(\hat{e}_j, \hat{f}_j) =
\begin{cases}
\left([e_{0},e_{1}]_{\bq_{0,1}},[f_{1},f_{0}]_{\bq_{0,1}^{-1}}\right) & \text{for $j=0$},\\
(e_{j+1},f_{j+1}) & \text{for $j\neq 0$}.
\end{cases}
\end{gather*}}
Here $\bq_{a,b}=\bq(\alpha_a,\alpha_b)$ for $a,b\in I$, and $[X,Y]_t = XY-tYX$ for $X,\,Y\in\mathcal{U}(\epsilon)$, $t\in\Bbbk^\times$.

\begin{thm}[{\cite[Theorem 4.3]{KY}}]\label{thm:folding homomorphism}
Under the above hypothesis, there exists a homomorphism of $\Bbbk$-algebras 
\begin{equation*}
\xymatrixcolsep{2pc}\xymatrixrowsep{3pc}
\xymatrix{
\phi^\e_{\e'} : \mathcal{U}(\e') \ \ar@{->}[r] & \ \mathcal{U}(\e),
}
\end{equation*}
such that for $l\in \mathbb{I}'$ and $j\in I'$
\begin{equation*}
\phi^\e_{\e'}(k'_{\delta'_l})=\hat{k}_{\delta'_l},\quad \phi^\e_{\e'}(e'_j)=\hat{e}_j,\quad \phi^\e_{\e'}(f'_j)=\hat{f}_j.
\end{equation*}
\end{thm}

More generally, suppose that the sequence $\e'=(\e'_1,\dots,\e'_{n-r})$ with $1\le r\le n-3$ is obtained from $\e$ by removing $\e_{i_1},\dots,\e_{i_r}$ for some $i_1<\dots<i_r$ .
For $0\leq k\leq r$, let $\e^{(k)}$ be a sequence such that
\begin{itemize}
\item[(1)] $\e^{(0)}=\e$, $\e^{(r)}=\e'$, 

\item[(2)] $\e^{(k)}$ is obtained from $\e^{(k-1)}$ by removing $\e_{i_k}$ for $1\leq k\leq r$.
\end{itemize}
We define a homomorphism of $\Bbbk$-algebras
\begin{equation*}
\xymatrixcolsep{2pc}\xymatrixrowsep{3pc}
\xymatrix{
\phi^\e_{\e'} : \mathcal{U}(\e') \ \ar@{->}[r] & \ \mathcal{U}(\e),
}
\end{equation*}
by $\phi^{\e}_{\e'}=\phi^{\e^{(0)}}_{\e^{(1)}}\circ \phi^{\e^{(1)}}_{\e^{(2)}}\circ \dots \circ\phi^{\e^{(r-1)}}_{\e^{(r)}}$.

For a $\mathcal{U}(\e)$-module $V$ in $\mc{C}(\e)$, let
\begin{equation}\label{eq:truncation-1}
\mf{tr}^\e_{\e'}(V) = \bigoplus_{\substack{\mu\in {\rm wt}(V) \\ (\mu|\de_{i_1})=\dots=(\mu|\de_{i_r})=0}}V_\mu.
\end{equation}
We denote by $\pi^\e_{\e'} : V \longrightarrow \mf{tr}^\e_{\e'}(V)$ the natural projection. For any $\mathcal{U}(\e)$-modules $V, W$ in $\mc{C}(\e)$ and $f \in {\rm Hom}_{\mathcal{U}(\e)}(V,W)$, let
\begin{equation*}
\xymatrixcolsep{2pc}\xymatrixrowsep{3pc}
\xymatrix{
\mf{tr}^\e_{\e'}(f) : \mf{tr}^\e_{\e'}(V) \ \ar@{->}[r] &  \mf{tr}^\e_{\e'}(W)
}
\end{equation*}
be the $\Bbbk$-linear map given by $\mf{tr}^\e_{\e'}(f)(v)=f(v)$ for $v\in \mf{tr}^\e_{\e'}(V)$. Then we have the following commutative diagram of $\Bbbk$-vector spaces:
\begin{equation*}\label{eq:truncation commuting diagram}
\xymatrixcolsep{4pc}\xymatrixrowsep{3pc}\xymatrix{
V \ar@{->}[d]_{\pi^\e_{\e'}}\ar@{->}[r]^{f} &  W \ar@{->}[d]^{\pi^\e_{\e'}} \\
\mf{tr}^\e_{\e'}(V) \ar@{->}[r]^{\mf{tr}^\e_{\e'}(f)} & \mf{tr}^\e_{\e'}(W)
}
\end{equation*}
\begin{prop}[{\cite[Propositions 4.4]{KY}}]\label{prop:truncation}
Under the above hypothesis,
\begin{itemize}
\item[(1)] the space $\mf{tr}^\e_{\e'}(V)$ is invariant under the action of $\mathcal{U}(\e')$ via $\phi^{\e}_{\e'}$, and hence a $\mathcal{U}(\e')$-module in $\mc{C}(\e')$,

\item[(2)] the map $\mf{tr}^\e_{\e'}(f) : \mf{tr}^\e_{\e'}(V) \longrightarrow \mf{tr}^\e_{\e'}(W)$ is $\mathcal{U}(\e')$-linear,

\item[(3)] the space $\mf{tr}^\e_{\e'}(V\ot W)$ is naturally isomorphic to the tensor product of $\mf{tr}^\e_{\e'}(V)$ and $\mf{tr}^\e_{\e'}(W)$ as a $\mathcal{U}(\e')$-module.

\end{itemize}
\end{prop}

Hence we have a functor 
\begin{equation*}
\xymatrixcolsep{2pc}\xymatrixrowsep{3pc}
\xymatrix{
\mf{tr}^\e_{\e'} : \mc{C}(\e)  \ar@{->}[r] &  \mc{C}(\e'),}
\end{equation*}
which we call {\em truncation}. Note that $\mf{tr}^\e_{\e'}$ is exact by its definition and monoidal by Proposition \ref{prop:truncation} in the sense of \cite[Appendix A.1]{KKK}.
We may also define $\mf{tr}^\e_{\e'} : \ring{\mc{C}}(\e) \longrightarrow \ring{\mc{C}}(\e')$ in the same way.

\begin{prop}[{\cite[Propositions 4.5 and 4.6]{KY}}]\label{prop:truncation of poly and fund}
Let $M'$ and $N'$ be the numbers of $j$'s with $\e'_j=0$ and $\e'_j=1$ in $\e'$, respectively.
\begin{itemize}
\item[(1)] For $\la\in \cP_{M|N}$, $\mf{tr}^\e_{\e'}(V_\e(\la))$ is non-zero if and only if $\la\in \cP_{M'|N'}$.
In this case, we have 
$\mf{tr}^\e_{\e'}(V_\e(\la)) \cong V_{\e'}(\la)$,
as a $\ring\mathcal{U}(\e')$-module.  

\item[(2)] For $l\in\mathbb{Z}_+$ and $x\in \Bbbk^\times$, $\mf{tr}^\e_{\e'}(\mathcal{W}_{l,\e}(x))$ is non-zero if and only if $(l)\in \cP_{M'|N'}$. In this case, we have 
$\mf{tr}^\e_{\e'}(\mathcal{W}_{l,\e}(x)) \cong \mathcal{W}_{l,\e'}(x)$ as a $\mathcal{U}(\e')$-module.

\end{itemize}
\end{prop}
\begin{rem}\label{rem:isomorphism on truncated space}
{\rm
We may assume that the isomorphism
\[
\mathfrak{tr}^{\epsilon}_{\epsilon^\prime}\left(\mathcal{W}_{\ell_1,\epsilon}(x_1)\otimes\cdots\otimes\mathcal{W}_{\ell_t,\epsilon}(x_t)\right) \cong \mathcal{W}_{\ell_1,\epsilon^\prime}(x_1)\otimes\cdots\otimes\mathcal{W}_{\ell_t,\epsilon^\prime}(x_t)
\]
induced from Proposition \ref{prop:truncation of poly and fund}(2), when restricted as a $\ring{\mathcal{U}}(\e')$-linear map, gives the isomorphism in Proposition \ref{prop:truncation of poly and fund}(1) on each $\mf{tr}^\e_{\e'}(V_\e(\la))$.
}
\end{rem}

\subsection{Fundamental representations and $R$ matrices}\label{subsec:R matrix for fund}
Let $l,m\in\mathbb{Z}_+$ and $x,y\in \Bbbk^{\times}$ be given. 
As a $\ring\mathcal{U}(\e)$-module,
\begin{equation}\label{eq:decomp of W_l and W_m}
\mathcal{W}_{l,\e}(x)\otimes \mathcal{W}_{m,\e}(y) \cong 
\bigoplus_{t\in H(l,m)} V_\e((l+m-t,t)),
\end{equation}
where $H(l,m)=\{\,t\,|\,0\leq t\leq \min\{l,m\}, (l+m-t,t)\in \cP_{M|N}\}$ (cf. \cite[Remark 3.5]{KY}).
Let $\e''=(\e''_1,\dots,\e''_{n''})$ be a sequence of $0,1$'s with $n''\gg n$ such that
\begin{itemize}
\item[(1)] $\e$ is a subsequence of $\e''$, 

\item[(2)] as a $\ring\mathcal{U}(\e'')$-module
\begin{equation*}\label{eq:decomp of W_l and W_m-e''}
\mathcal{W}_{l,\e''}(x)\otimes \mathcal{W}_{m,\e''}(y) \cong 
\bigoplus_{0\leq t\leq \min\{l,m\}} V_{\e''}((l+m-t,t)),
\end{equation*}

\item[(3)] if $\e'=\e_{M''|0}$ with $M''=|\{\,i\,|\,\e''_i=0\,\}|$, then as a $\ring\mathcal{U}(\e')$-module
\begin{equation}\label{eq:decomp of W_l and W_m-e'}
\mathcal{W}_{l,\e'}(x)\otimes \mathcal{W}_{m,\e'}(y) \cong 
\bigoplus_{0\leq t\leq \min\{l,m\}} V_{\e'}((l+m-t,t)).
\end{equation}
\end{itemize}
Note that we have the following:
\begin{equation*} 
\xymatrixcolsep{0pc}\xymatrixrowsep{2pc}\xymatrix{
& \mathcal{W}_{l,\e^{''}}(x)\otimes \mathcal{W}_{m,\e''}(y)   \ar@{->}_{\pi^{\e''}_{\e}}[dl]\ar@{->}^{\pi^{\e''}_{\e'}}[dr] &  \\
\mathcal{W}_{l,\e}(x)\otimes \mathcal{W}_{m,\e}(y) & &  \mathcal{W}_{l,\e'}(x)\otimes \mathcal{W}_{m,\e'}(y)
}
\end{equation*}

For $0\leq t\leq \min\{l,m\}$, let $v'(l,m,t)$ be the highest weight vector of $V_{\e'}((l+m-t,t))$ in $\mathcal{W}_{l,\e'}(x)\otimes \mathcal{W}_{m,\e'}(y)$ such that 
\begin{equation*}\label{eq:extremal wt vec}
\begin{split}
& v'(l,m,t)\in \mc{L}_{l,\e'}\ot \mc{L}_{m,\e'},\\
& v'(l,m,t)\equiv | l\be_1 \rangle\ot | (m-t)\be_1 + t\be_2 \rangle \pmod{q\mc{L}_{l,\e'}\ot \mc{L}_{m,\e'}},
\end{split}
\end{equation*}
where $\mc{L}_{s,\e'}$ is the lower crystal lattice of $\mathcal{W}_{s,\e'}$ spanned by $|{\bf m}\rangle$.
We also define $v'(m,l,t)$ in the same manner. 
By Proposition \ref{prop:truncation of poly and fund}(1), we may regard 
\begin{equation*}
\begin{split}
&V_{\e'}((l+m-t,t))\subset V_{\e''}((l+m-t,t))\quad (0\leq t\leq \min\{l,m\}), \\
&V_{\e}((l+m-t,t))\subset V_{\e''}((l+m-t,t)) \quad (t\in H(l,m)),
\end{split}
\end{equation*}
as a $\Bbbk$-space.

For $0\leq t\leq \min\{l,m\}$, 
let 
$\mc P^{l,m}_{t} : \mathcal{W}_{l,\e''}(x)\otimes \mathcal{W}_{m,\e''}(y) \longrightarrow \mathcal{W}_{m,\e''}(y)\otimes \mathcal{W}_{l,\e''}(x)$
be a $\ring\mathcal{U}(\e'')$-linear map given by $\mc P^{l,m}_{t}(v'(l,m,t'))=\delta_{tt'}v'(m,l,t')$. 
For $t\in H(l,m)$, let
\begin{equation}\label{eq:projection of extremal component}
\xymatrixcolsep{2pc}\xymatrix{
\mc P^{l,m}_{t} : \mathcal{W}_{l,\e}(x)\otimes \mathcal{W}_{m,\e}(y) \ar@{->}[r] & \mathcal{W}_{m,\e}(y)\otimes \mathcal{W}_{l,\e}(x)}
\end{equation}
be its restriction onto $\mathcal{W}_{l,\e}(x)\otimes \mathcal{W}_{m,\e}(y)$.
Note that $\mc{P}^{l,m}_t$ on $\mathcal{W}_{l,\e}(x)\otimes \mathcal{W}_{m,\e}(y)$ is independent of the choice of $\e''$ and $x,y$.

\begin{rem}{\rm
When we define $\mc{P}^{l,m}_t$ \eqref{eq:projection of extremal component},
we may take $\e'=\e_{0|N''}$ with $N''=|\{\,i\,|\,\e''_i=1\,\}|$ such that \eqref{eq:decomp of W_l and W_m-e'} holds.
Recall that $\mathcal{U}(\e')\cong U^{\prime}_{-q^{-1}}(A_{N''-1}^{(1)})$ (up to a central element).
In this case, we take $v'(l,m,t)$ to be the highest weight vector of $V_{\e'}((l+m-t,t))$ in $\mathcal{W}_{l,\e'}(x)\otimes \mathcal{W}_{m,\e'}(y)$ such that
\begin{equation*}\label{eq:extremal wt vec-2}
\begin{split}
& v'(l,m,t)\in \mc{L}_{l,\e'}\ot \mc{L}_{m,\e'},\\
& v'(l,m,t)\equiv | \be_{1,t} + \be_{m+1,l+m-t} \rangle\ot | \be_{1,m} \rangle \pmod{q^{-1}\mc{L}_{l,\e'}\ot \mc{L}_{m,\e'}}, 
\end{split}
\end{equation*}
where $\be_{a,b}=\sum_{a\leq i\leq b}\be_i$ and $\mc{L}_{l,\e'}$ is the lower crystal lattice over the subring of $f(q)\in \Bbbk$ regular at $q=\infty$.
}
\end{rem}

Let $z_1, z_2$ be indeterminates and
let 
\begin{equation*}
\begin{split}
(\mathcal{W}_{l,\e})_{\rm aff}=\Bbbk[z_1^{\pm 1}]\ot\mathcal{W}_{l,\e},\quad
(\mathcal{W}_{m,\e})_{\rm aff}=\Bbbk[z_2^{\pm 1}]\ot\mathcal{W}_{m,\e},\\
(\mathcal{W}_{l,\e})^{{}^{\wedge}}_{\rm aff}=\Bbbk(z_1)\ot \mathcal{W}_{l,\e},\quad
(\mathcal{W}_{m,\e})^{{}^{\wedge}}_{\rm aff}=\Bbbk(z_2)\ot \mathcal{W}_{m,\e}.
\end{split}
\end{equation*}
Note that
\begin{equation*}
\begin{split}
\Bbbk(z_1,z_2)\ot_{\Bbbk[z_1^{\pm 1},z_2^{\pm 1}]} ((\mathcal{W}_{l,\e})_{\rm aff}\ot (\mathcal{W}_{m,\e})_{\rm aff}) \cong (\mathcal{W}_{l,\e})^{{}^{\wedge}}_{\rm aff} \ot (\mathcal{W}_{m,\e})^{{}^{\wedge}}_{\rm aff}.
\end{split}
\end{equation*}

\begin{thm}\label{thm:irreducibility of fund tensor func}
For $l, m\in \mathbb{Z}_+$, $(\mathcal{W}_{l,\e})^{{}^{\wedge}}_{\rm aff} \ot (\mathcal{W}_{m,\e})^{{}^{\wedge}}_{\rm aff}$ is an irreducible representation of $\Bbbk(z_1,z_2)\ot \mathcal{U}(\e)$. 
\end{thm}
\pf By \cite[Theorem 4.7]{KY}, $\mathcal{W}_{l,\e}\ot \mathcal{W}_{m,\e}$ is irreducible. Hence the irreducibility of $\Bbbk(z_1,z_2)\ot_{\Bbbk[z_1^{\pm 1},z_2^{\pm 1}]} ((\mathcal{W}_{l,\e})_{\rm aff}\ot (\mathcal{W}_{m,\e})_{\rm aff})$ follows from \cite[Lemma 3.4.2]{KMN}.
\qed
\vskip 2mm

Suppose that $M\neq N$.
Let $\mc{R}^{\rm univ}_{l,m}$ denote the universal $R$ matrix acting on $(\mathcal{W}_{l,\e})_{\rm aff}\ot (\mathcal{W}_{m,\e})_{\rm aff}$ \eqref{eq:R univ}:
\begin{equation*}\label{eq:R matrix on W_l tensor W_m}
\xymatrixcolsep{3pc}\xymatrixrowsep{3pc}\xymatrix{
(\mathcal{W}_{l,\e})_{\rm aff} \ot (\mathcal{W}_{m,\e})_{\rm aff}\ \ar@{->}^{\hskip -2cm\mc{R}^{\rm univ}_{l,m}}[r] &  \
\Bbbk[[z_1/z_2]]\ot_{\Bbbk[z_1/z_2]}\left((\mathcal{W}_{m,\e})_{\rm aff}\,{\ot}\, (\mathcal{W}_{l,\e})_{\rm aff}\right). 
}
\end{equation*}
Let $s=\max H(l,m)=\min\{l,m,n-1\}$  and $V_s=V_{\e}((l+m-s,s))$. 
Then
\begin{equation*}
\mc{R}^{\rm univ}_{l,m}\big\vert_{V_s} = \varphi_{l,m}(z_1/z_2){\rm id}_{V_s},
\end{equation*}
for some $\varphi_{l,m}(z_1/z_2)\in \Bbbk[[z_1/z_2]]\setminus\{0\}$. 
Put $z=z_1/z_2$ and let 
\begin{equation*}
c_s(z)=\prod^{\min\{l,m\}}_{i=s+1}\frac{1-q^{l+m-2i+2}z}{z-q^{l+m-2i+2}},
\end{equation*}
where we understand $c_s(z)=1$ in case of $s=\min\{l,m\}$.
Note that we have $s=\min\{l,m\}$ for $M\geq 2$.

Now, we define the normalized $R$ matrix by 
\begin{equation}\label{eq:R_{l,m}}
\mc{R}^{\rm norm}_{l,m} =  \varphi_{l,m}(z)^{-1} c_s(z)\mc{R}^{\rm univ}_{l,m}.
\end{equation}
Here we understand $\varphi_{l,m}(z)^{-1}$ as an element in the quotient field $\Bbbk(\!(z_1/z_2)\!)$ of $\Bbbk[[z_1/z_2]]$.

Since $\mc{R}^{\rm norm}_{l,m}\big\vert_{V_s}= c_s(z) {\rm id}_{V_s}$, 
we have by Theorem \ref{thm:irreducibility of fund tensor func} a $\Bbbk[z_1^{\pm 1},z_2^{\pm 1}]\ot\mathcal{U}(\e)$-linear map
\begin{equation}\label{eq:R matrix on W_l tensor W_m}
\xymatrixcolsep{3pc}\xymatrixrowsep{3pc}\xymatrix{
(\mathcal{W}_{l,\e})_{\rm aff} \ot (\mathcal{W}_{m,\e})_{\rm aff}\ \ar@{->}^{\hskip -2cm\mc{R}^{\rm norm}_{l,m}}[r] &  \
\Bbbk(z_1,z_2)\ot_{\Bbbk[z_1^{\pm 1},z_2^{\pm 1}]}\left((\mathcal{W}_{m,\e})_{\rm aff}\,{\ot}\, (\mathcal{W}_{l,\e})_{\rm aff}\right),
}
\end{equation}
or a $\Bbbk(z_1,z_2)\ot \mathcal{U}(\e)$-linear map
\begin{equation}\label{eq:R matrix on W_l tensor W_m-2}
\xymatrixcolsep{3pc}\xymatrixrowsep{3pc}\xymatrix{
((\mathcal{W}_{l,\e})_{\rm aff} \ot (\mathcal{W}_{m,\e})_{\rm aff})^{{}^{\wedge}}\ \ar@{->}^{\!\!\mc{R}^{\rm norm}_{l,m}}[r] &  \
((\mathcal{W}_{m,\e})_{\rm aff} \ot (\mathcal{W}_{l,\e})_{\rm aff})^{{}^{\wedge}},
}
\end{equation}
where 
$
((\mathcal{W}_{l,\e})_{\rm aff} \ot (\mathcal{W}_{m,\e})_{\rm aff})^{{}^{\wedge}}
=
\Bbbk(z_1,z_2){\ot}_{\Bbbk [z_1^{\pm 1},z_2^{\pm 1}]}((\mathcal{W}_{l,\e})_{\rm aff} \ot (\mathcal{W}_{m,\e})_{\rm aff})
$.

\begin{cor}\label{cor:characterization of R^{norm}}
Under the above hypothesis, $\mc{R}^{\rm norm}_{l,m}$ is a unique $\Bbbk(z_1,z_2)\ot \mathcal{U}(\e)$-linear map \eqref{eq:R matrix on W_l tensor W_m} or \eqref{eq:R matrix on W_l tensor W_m-2}
such that $\mc{R}^{\rm norm}_{l,m}\big\vert_{V_s}= c_s(z){\rm id}_{V_s}$.
\end{cor}

\pf It follows immediately from Theorem \ref{thm:irreducibility of fund tensor func}.
\qed\vskip 2mm

We have the following spectral decomposition of $\mc{R}^{\rm norm}_{l,m}$.

\begin{thm}\label{thm:spectral decomposition}
For $l, m\in\mathbb{Z}_+$, we have
\begin{equation*}\label{eqn:spectral-decomp}
\mc{R}^{\rm norm}_{l,m} = \sum_{t\in H(l,m)}\rho_t(z) \mc P^{l,m}_t, \quad 
\rho_t(z)=\prod_{i=t+1}^{\min\{l,m\}}\dfrac{1-q^{l+m-2i+2}z}{z-q^{l+m-2i+2}},
\end{equation*}
where $z=z_1/z_2$ and $\rho_{\min\{l,m\}}(z)=1$.
\end{thm}
\pf It is shown in \cite{KOS} that there is a $\Bbbk(z_1 ,z_2)\otimes\mathcal{U}(\e)$-linear map 
$R_{l,m}: ((\mathcal{W}_{l,\e})_{\rm aff} \ot (\mathcal{W}_{m,\e})_{\rm aff})^{{}^{\wedge}}
\longrightarrow
((\mathcal{W}_{m,\e})_{\rm aff} \ot (\mathcal{W}_{l,\e})_{\rm aff})^{{}^{\wedge}}$ such that $R_{l,m}\big\vert_{V_s}=c_s(z){\rm id}_{V_s}$, where $s=\min\{l,m,n-1\}$. 
The spectral decomposition of $R_{l,m}$ is given in \cite[Theorem 5.2]{KY} as follows:
\begin{equation*}
{R}_{l,m} = \sum_{t\in H(l,m)}\rho_t(z) \mc P^{l,m}_t, \quad 
\rho_t(z)=\prod_{i=t+1}^{\min\{l,m\}}\dfrac{1-q^{l+m-2i+2}z}{z-q^{l+m-2i+2}}.
\end{equation*} 
By Corollary \ref{cor:characterization of R^{norm}}, we have $R_{l,m}=\mc{R}^{\rm norm}_{l,m}$, which implies the decomposition of $\mc{R}^{\rm norm}_{l,m}$.
\qed
\vskip 2mm
\begin{rem}\label{rem:spectral decomp is the same}
{\rm
If we replace $z$ and $q$ with $z^{-1}$ and $-q^{-1}$, respectively in Theorem \ref{thm:spectral decomposition}, then it recovers the formula \cite{DO} when $\e=\e_{0|n}$. It is more important to observe that the spectral decomposition of $\mc{R}^{\rm norm}_{l,m}$ is independent of the choice of $\e$ if $n$ is large enough. This will play a crucial role in the remaining of the paper.
}
\end{rem}

For $t\geq 2$, let $z_1,\dots ,z_t$ be indeterminates and let $l_1,\dots ,l_t\in \mathbb{Z}_+$ given.
Let $W_i=(\mathcal{W}_{l_i,\e})_{\rm aff}=\Bbbk[z_i^\pm]\ot \mathcal{W}_{l_i,\e}$ for $1\leq i\leq t$. Let $\mf{S}_t$ be the group of permutations on $t$ letters. 
Since $\mc{R}^{\rm norm}_{l.m}$ satisfies the Yang-Baxter equation, we can define a $\Bbbk[z_1^\pm,\dots, z_t^\pm]\ot\mathcal{U}(\e)$-linear map
\begin{equation*}
\xymatrixcolsep{2pc}\xymatrixrowsep{3pc}\xymatrix{
\mc{R}_{l_1,\dots, l_t, w}^{\rm norm} : W_1\ot\cdots \ot W_t \ar@{->}[r] & 
\Bbbk(z_1,\dots,z_t)\ot_{\Bbbk[z_1^\pm,\dots, z_t^\pm]}W_{w(1)}\ot\cdots \ot W_{w(t)}
}
\end{equation*}
or a $\Bbbk(z_1,\dots,z_t)\ot\mathcal{U}(\e)$-linear map
\begin{equation*}\label{eq:R matrix on W_l's}
\xymatrixcolsep{2pc}\xymatrixrowsep{3pc}\xymatrix{
\mc{R}_{l_1,\dots, l_t, w}^{\rm norm} : (W_1\ot\cdots \ot W_t)^{{}^\wedge} \ar@{->}[r] & 
(W_{w(1)}\ot\cdots \ot W_{w(t)})^{{}^\wedge},
}
\end{equation*}
for $w\in \mf{S}_t$ by a composition of $\mc{R}^{\rm norm}_{l,m}$'s with respect to a reduced expression of $w$, where 
$(W_{w(1)}\ot\cdots \ot W_{w(t)})^{{}^\wedge}=
\Bbbk(z_1,\dots,z_t)\ot_{\Bbbk[z_1^\pm,\dots, z_t^\pm]}W_{w(1)}\ot\cdots \ot W_{w(t)}$. 
We put 
\begin{equation}\label{eq:R matrix for longest permutation}
\mc{R}_{l_1,\dots ,l_t}^{\rm norm}=\mc{R}_{l_1,\dots, l_t,w_0}^{\rm norm},
\end{equation} 
where $w_0$ is the longest element in $\mf{S}_t$.

\section{Irreducible modules}\label{sec:irr char}
From now on, we assume that $M\neq N$ for $\e$.
\subsection{Renormalized $R$ matrix and fusion construction}
Let $V, W$ be $\mathcal{U}(\e)$-modules in $\mc{C}(\e)$.
Let $\mc{R}^{\rm univ}_{V,W}$ be the universal $R$ matrix on $V_{\rm aff}\ot W_{\rm aff}$ \eqref{eq:R univ}. Following \cite{KKKO15}, we say that $\mathcal{R}_{V,W}^{\mathrm{univ}}$ is {\em rationally renormalizable}
if there exists $a\in\Bbbk(\!(z_{1}/z_{2})\!)^\times$ such that 
\begin{equation*}
\xymatrixcolsep{2pc}\xymatrixrowsep{3pc}\xymatrix{
a\mathcal{R}_{V,W}^{\mathrm{univ}} : V_{\rm aff}\otimes W_{\rm aff}\ar@{->}[r] & 
 W_{\rm aff}\otimes V_{\rm aff}.
}
\end{equation*}
Suppose that $\Bbbk$ is the algebraic closure of $\mathbb{Q}(q)$ in $\bigcup_{m>0}\mathbb{C}(\!(q^{\frac{1}{m}})\!)$.
If such $a$ exists, then we can choose $a$ (unique up to multiplication by a power of $z_1 / z_2$) such that 
$a\mathcal{R}_{V,W}^{\mathrm{univ}}\big\vert_{z_1=c_1,z_2=c_2}$ is nonzero for any $c_{1},c_{2}\in\Bbbk^{\times}$.
Then we have a nonzero $\mathcal{U}(\epsilon)$-linear map
\begin{equation*}
\xymatrixcolsep{2pc}\xymatrixrowsep{3pc}\xymatrix{
\mathbf{r}_{V,W}=a\mathcal{R}_{V,W}^{\mathrm{univ}}|_{z_{1}=z_{2}=1} : V \otimes W \ar@{->}[r] &  W \otimes V,
}
\end{equation*}
which we call the {\em renormalized $R$ matrix}.

The renormalized $R$ matrices $\mathbf{r}_{V,W}$'s also satisfy the hexagon property and the Yang-Baxter equation up to a multiplication by (possibly zero) constant, since they are specializations of scalar multiples of $\mc{R}^{\rm univ}$ \cite{KKKO15}. 

Let $V$, $W_{1}$, and $W_{2}$ be  $\mathcal{U}(\epsilon)$-modules in $\mc{C}(\e)$, and let $f\in {\rm Hom}_{\mathcal{U}(\e)}(W_{1},W_{2})$. 
By the definition of $\mc{R}^{\rm univ}$, the following diagram 
\begin{equation*}
\xymatrixcolsep{4pc}\xymatrixrowsep{3pc}\xymatrix{
V\otimes W_1  \ar@{->}[d]_{{\rm id}_V\ot f} \ar@{->}^{{\bf r}_{V,W_1}}[r] & W_1\ot V  \ar@{->}[d]^{f\ot {\rm id}_V} \\
V\otimes W_2  \ar@{->}^{{\bf r}_{V,W_2}}[r] & W_2\ot V
}
\end{equation*}
is commutative up to a constant multiple, provided ${\bf r}_{V,W_1}$ and ${\bf r}_{V,W_2}$ exist.

\begin{rem}{\rm
When $\mathcal{U}(\e)$ is the usual quantum affine algebra, that is, either $\e=\e_{n|0}$ or $\e_{0|n}$, $\mathcal{R}_{V,W}^{\mathrm{univ}}$ is always rationally renormalizable for any irreducible $V$ and $W$ \cite[Section 2.2]{KKKO15}. But we do not know yet whether it is also true in case of $\mc{C}(\e)$ for arbitrary $\e$.
}
\end{rem}

\begin{lem}[cf.~{\cite[Propositions 2.11 and 2.12]{KKOP}}]\label{lem:renormalizable cases}
Under the above hypothesis, $\mathcal{R}_{V,W}^{\mathrm{univ}}$ is rationally
renormalizable in one of the following cases:
\begin{enumerate}

\item $V$ (resp. $W$) is a submodule or a quotient of $V_0$ (resp. $W_0$) and $\mathcal{R}_{V_0,W}^{\mathrm{univ}}$ (resp. $\mathcal{R}_{V,W_0}^{\mathrm{univ}}$) is rationally renormalizable,

\item $V=V_{1}\ot V_{2}$ (resp. $W=W_{1}\ot W_{2}$) and both $\mathcal{R}_{V_{1},W}^{\mathrm{univ}}$
and $\mathcal{R}_{V_{2},W}^{\mathrm{univ}}$ (resp. $\mathcal{R}_{V,W_{1}}^{\mathrm{univ}}$
and $\mathcal{R}_{V,W_{2}}^{\mathrm{univ}}$) are rationally renormalizable. 

\end{enumerate}
\end{lem}
\pf
(1) It is clear. (2) It follows from the hexagon property of $\mc{R}^{\rm univ}$ (cf.~ \cite[32.2]{Lu93} and Remark \ref{rem:properties of Theta}).
\qed\vskip 2mm

\begin{thm}[cf.~{\cite[Theorems 3.2 and 3.12]{KKKO15}}]\label{thm:KKKO15-main}
Let $V$, $W$ be irreducible $\mathcal{U}(\epsilon)$-modules in $\mc{C}(\e)$. 
Suppose that $\mc{R}^{\mathrm{univ}}_{V,V}$, $\mc{R}^{\mathrm{univ}}_{W,W}$ and $\mc{R}^{\mathrm{univ}}_{V,W}$ are rationally renormalizable and 
\[
\mathbf{r}_{V,V}\in \Bbbk\,\mathrm{id}_{V^{\ot 2}}
\quad \text{or} \quad
\mathbf{r}_{W,W}\in \Bbbk\,\mathrm{id}_{W^{\ot 2}}.
\]
Then $\mathrm{Im}(\mathbf{r}_{V,W})$ is irreducible, and isomorphic to the head of $V\otimes W$ and the socle of $W\otimes V$.
\end{thm}
\pf 
We follow the argument in \cite{KKKO15}. But we give a self-contained proof here for the reader's convenience. Let us assume $\mathbf{r}_{W,W}\in \Bbbk\mathrm{id}_{W^{\ot 2}}$ since a similar argument works for the other case.

Let $S\subset W\ot V$ be a non-zero submodule. 
By assumption, there exist $a,b\in \Bbbk(\!(z_1/z_2)\!)^\times$ such that we have $\mathbf{r}_{V,W}=a\mc{R}^{\rm{univ}}_{V,W}|_{z_1=z_2=1}$ and $\mathbf{r}_{W,W}=b\mc{R}^{\rm{univ}}_{W,W}|_{z_1=z_2=1}$.
By Lemma \ref{lem:renormalizable cases}, $\mc{R}^{\rm{univ}}_{S,W}$ is also rationally renormalizable by $ab$, that is,
\[
ab\mc{R}^{\rm{univ}}_{S,W}\,:\,S_\mathrm{aff} \otimes W_\mathrm{aff} \rightarrow W_\mathrm{aff} \otimes S_\mathrm{aff}.
\]
Multiplying the rows of the diagram by $ab$ and specializing at 1, we obtain a commutative diagram
\begin{equation*}
\xymatrixcolsep{4pc}\xymatrixrowsep{3pc}\xymatrix{
S \ot W \ar@{^{(}->}[d] \ar@{->}^{c\,\mathbf{r}_{S,W}}[rr]& & W\ot S \ar@{^{(}->}[d] \\
W\ot V\ot W \ar@{->}^{{\rm id}_{W}\ot \mathbf{r}_{V,W}}[r] & W\ot W \ot V
\ar@{->}^{\mathrm{id}_{W^{\ot 2}}\ot{\rm id}_{V}}[r] & W\ot W \ot V
}
\end{equation*}
for some $c$, where we use the assumption $\mathbf{r}_{W,W} = \mathrm{id}_{W\otimes W}$ up to a multiplication by a non-zero scalar. 
Therefore we have $S\ot W \subset W \ot {\mathbf{r}_{V,W}^{-1}}(S)$.

By \cite[Lemma 3.10]{KKKO15}, there exists a submodule $K$ of $V$ such that $S\subset W\ot K$ and $K\ot W\subset {\mathbf{r}_{V,W}^{-1}}(S)$. Since $S\neq 0$, we have $K\neq 0$ and hence $K=V$, which implies that $V\ot W\subset {\mathbf{r}_{V,W}^{-1}}(S)$ and ${\rm Im}(\mathbf{r}_{V,W})\subset S$. It follows that ${\rm Im}(\mathbf{r}_{V,W})$ is a unique simple submodule of $W\otimes V$ and so equals to the socle of $W\otimes V$.
\qed

\begin{thm}\label{thm:main-1-renorm}
Let $V_1,\,\dots,\,V_t$ be irreducible $\mathcal{U}(\epsilon)$-modules in $\mc{C}(\e)$. Suppose that $\mathcal{R}^{\mathrm{univ}}_{V_i,V_j}$ are rationally renormalizable and $\mathbf{r}_{V_i,V_i}\in \Bbbk\,\mathrm{id}_{V_i^{\ot 2}}$ for any $1\leq i,j\leq t$. Let
\[
\xymatrixcolsep{2pc}\xymatrixrowsep{3pc}\xymatrix{
\mathbf{r} : V_1 \otimes \cdots \otimes V_t \ \ar@{->}[r] & \ V_t \otimes \cdots \otimes V_1
}
\]
be the composition of ${\bf r}_{V_i,V_j}$ associated to a reduced expression of the longest element in $\mf{S}_t$.
Then ${\rm Im}\,{\bf r}$ is irreducible if it is not zero.
\end{thm}

\pf
Use induction on $t$. It is true for $t=2$ by Theorem \ref{thm:KKKO15-main}.

Suppose that $t\geq 3$. Let ${\bf r}=\mathbf{r}_t$ be the map in the statement and $\mathbf{r}_{t-1}$ denote the map corresponding to the first $t-1$ factors.
Consider the following commutative diagram:
\begin{equation*}
\xymatrixcolsep{5pc}\xymatrixrowsep{2pc}
\xymatrix{
V_1\ot\cdots\ot V_t \ar@{->}[r]^{\mathbf{r}_t} \ar@{->}[d]^{\mathbf{r}_{t-1}\ot {\rm id}_{V_t}} & \ V_t\ot\cdots\ot V_1 \\
{\rm Im}(\mathbf{r}_{t-1})\ot V_t \ar@{^{(}->}[d] & \\
V_{t-1}\ot\cdots\ot V_1\ot V_t  \ar@{->}[ruu]_{\quad \hat{\mathbf{r}}_{t-1}\circ\dots\circ\hat{\mathbf{r}}_1} & \\
}
\end{equation*}
where 
$\hat{\mathbf{r}}_i= {\rm id}^{\otimes t-i-1}\ot \mathbf{r}_{V_i,V_t} \ot {\rm id}^{i-1}$ for $1\leq i\leq t-1$. 

Assume $\mathbf{r}_t \neq 0$. Then $\mathbf{r}_{t-1}$ has a nonzero image $L$, which is irreducible by the induction hypothesis. 
Applying the hexagon property repeatedly, we have the following diagram:
\begin{equation*}
\xymatrixcolsep{6pc}\xymatrixrowsep{3pc}\xymatrix{
L \ot V_t \ar@{^{(}->}[d] \ar@{->}^{\mathbf{r}_{L,V_t}}[r] &  V_t\ot L \ar@{^{(}->}[d] \\
V_{t-1}\ot\cdots\ot V_1\ot V_t \ar@{->}^{\hat{\mathbf{r}}_{t-1}\circ\cdots\circ\hat{\mathbf{r}}_1}[r] & V_{t}\ot V_{t-1}\ot\cdots\ot V_1
}
\end{equation*}
Note that it commutes up to multiplication by a non-zero scalar since both horizontal maps are nonzero.
Thus the image of $\mathbf{r}_t$ is equal to that of $\mathbf{r}_{L,V_t}$, which is simple as seen in Theorem \ref{thm:KKKO15-main}.
\qed\vskip 2mm

Let $l,m\in \mathbb{Z}_+$ be given. By Theorem \ref{thm:spectral decomposition}, we have
\begin{equation}\label{eq:renorm R for fund}
\xymatrixcolsep{2pc}\xymatrixrowsep{3pc}\xymatrix{
{\bf r}_{l,m}:=d_{l,m}(z)\mc{R}^{\rm norm}_{l,m} : (\mathcal{W}_{l,\e})_{\rm aff} \ot (\mathcal{W}_{m,\e})_{\rm aff}\ \ar@{->}[r] &  \
(\mathcal{W}_{m,\e})_{\rm aff} \ot (\mathcal{W}_{l,\e})_{\rm aff},
}
\end{equation}
where $z=z_1/z_2$ and 
\begin{equation}\label{eq:denominator}
d_{l,m}(z) = \prod_{i=1}^{\min\{l,m\}}(z-q^{l+m-2i+2}).
\end{equation}
Note that ${\bf r}_{l,m}$ is not zero for any specialization at $(z_1,z_2)=(x,y)\in (\Bbbk^\times)^2$. Hence $\mc{R}^{\rm univ}_{V,W}$ is rationally renormalizable when $V=\mathcal{W}_{l,\e}$ and $W=\mathcal{W}_{m,\e}$.

Suppose that $V_i=\mathcal{W}_{l_i,\e}(c_i)$ for some $l_i\in \mathbb{Z}_+$ and $c_i\in \Bbbk^\times$ for $1\leq i\leq t$.
We have ${\bf r}_{V_i,V_j}={\bf r}_{l_i,l_j}$ for $1\leq i,j\leq t$ given in \eqref{eq:renorm R for fund}. 
Note that ${\bf r}_{l_i,l_i} \in \Bbbk\, \mathrm{id}_{V_i^{\ot 2}}$ by putting $z=c_{i}/c_{i}=1$ in the spectral decomposition of Theorem \ref{thm:spectral decomposition}. 
The following follows immediately from Theorem \ref{thm:main-1-renorm}.

\begin{cor}\label{thm:main-1}
Suppose that $c_1,\dots,c_t\in \Bbbk^\times$ are given such that $c_i/c_j$ is not a zero of $d_{l_i,l_j}(z_i/z_j)$ for any $1\leq i<j\leq t$. Let 
\begin{equation}\label{eq:R matrix on standard module}
\xymatrixcolsep{2pc}\xymatrixrowsep{3pc}\xymatrix{
R : 
\mathcal{W}_{l_1,\e}(c_1)\ot \dots \otimes \mathcal{W}_{l_s,\e}(c_t) \ar@{->}[r] \ &
\mathcal{W}_{l_t,\e}(c_t)\ot \dots \otimes \mathcal{W}_{l_1,\e}(c_1)}
\end{equation}
be the specialization of $\mc{R}_{l_1,\dots l_t}^{\rm norm}$ in \eqref{eq:R matrix for longest permutation} at $(z_1,\dots,z_t)=(c_1,\dots,c_t)$. 
Then ${\rm Im}R$ is irreducible if it is not zero. 
\end{cor}

As another application, we have the following, which plays an important role in later sections.

\begin{prop}[cf.~{\cite[Proposition 3.2.9]{KKKO18}}]\label{prop:Hom between fundamental tensor}
Let $V$, $W$ be as in Theorem \ref{thm:KKKO15-main}.
Then for $f\in {\rm Hom}_{\mathcal{U}(\e)}(V\otimes W, W\otimes V)$, the
image of $f$ lies in the socle of $W\otimes V$.
\end{prop}
\pf
The same argument in the proof of \cite[Proposition 3.2.9]{KKKO18} for modules over quiver Hecke algebras also applies to our case.
\qed

\begin{cor}\label{cor:hom-dim}
For $l,m\geq1$ and $x,y\in\Bbbk^{\times}$, we have
\[
\mathrm{Hom}_{\mathcal{U}(\epsilon)}\left(\mathcal{W}_{l,\e}(x)\otimes\mathcal{W}_{m,\e}(y),\mathcal{W}_{m,\e}(y)\otimes\mathcal{W}_{l,\e}(x)\right)=\Bbbk\,\mathbf{r}_{l,m}.
\]
\end{cor}
\pf Again by Theorem \ref{thm:spectral decomposition}, we have ${\bf r}_{l,l}\in \Bbbk\,{\rm id}_{\mathcal{W}_{l,\e}(x)^{\ot 2}}$. 
Let $f$ be a non-zero $\mathcal{U}(\e)$-linear map from $\mathcal{W}_{l,\e}(x)\otimes\mathcal{W}_{m,\e}(y)$ to $\mathcal{W}_{m,\e}(y)\otimes\mathcal{W}_{l,\e}(x)$. 
Then ${\rm Im}f$ is isomorphic to the socle of $\mathcal{W}_{m,\e}(y)\otimes\mathcal{W}_{l,\e}(x)$ by Proposition \ref{prop:Hom between fundamental tensor}, which is isomorphic to the head of $\mathcal{W}_{l,\e}(x)\otimes\mathcal{W}_{m,\e}(y)$ and irreducible by Theorem \ref{thm:KKKO15-main}. Since the head of $\mathcal{W}_{l,\e}(x)\otimes\mathcal{W}_{m,\e}(y)$ is irreducible, ${\rm Ker}\,f$ is equal to ${\rm Ker}\,{\bf r}_{l,m}$. Hence by Schur's lemma, $f$ is equal to ${\bf r}_{l,m}$ up to a multiplication by a non-zero scalar. 
\qed

\subsection{Irreducible modules in $\mc{C}(\e)$ and truncation}
Let $\mc{P}^+$ be the set of $({\mb l},{\mb c})$ such that
\begin{itemize}
\item[(1)] ${\mb l}=(l_1,\dots,l_t)\in \mathbb{Z}_+^t$, and $\mb{c}=(c_1,\dots,c_t)\in (\Bbbk^\times)^t$ for some $t\geq 1$,

\item[(2)] $c_i/c_j$ is not a zero of $d_{l_i,l_j}(z_i/z_j)$ for any $1\leq i<j\leq t$, when $t\geq 2$.
\end{itemize}
For $({\mb l},{\mb c})\in \mc{P}^+$, we put 
\begin{itemize}
\item[$\bullet$] $\mathcal{W}_{\mb{l},\e}(\mb{c})=\mathcal{W}_{l_1,\e}(c_1)\ot\dots\ot \mathcal{W}_{l_t,\e}(c_t)$,

\item[$\bullet$] $\mathcal{W}_{\mb{l}^{\rm rev},\e}(\mb{c}^{\rm rev})=\mathcal{W}_{l_t,\e}(c_t)\ot\dots\ot \mathcal{W}_{l_1,\e}(c_1)$,

\item[$\bullet$] $R_{\mb{l},\e}(\mb c) : \mathcal{W}_{\mb{l},\e}(\mb{c}) \longrightarrow \mathcal{W}_{\mb{l}^{\rm rev},\e}(\mb{c}^{\rm rev})$ : the $\mathcal{U}(\e)$-linear map in  \eqref{eq:R matrix on standard module},

\item[$\bullet$] $\mathcal{W}_{\e}({\mb l},{\mb c})$ : the image of $R_{\mb{l},\e}(\mb c)$,

\item[$\bullet$] $\mc{P}^+(\e)=\{\,({\mb l},{\mb c})\,|\,\mathcal{W}_{\e}({\mb l},{\mb c})\neq 0\,\}$.
\end{itemize}
We assume that $\mathcal{W}_{\e}({\mb l},{\mb c})=\mathcal{W}_{l_1,\e}(c_1)$ when $t=1$.
Recall that $\mathcal{W}_{\e}({\mb l},{\mb c})$ is irreducible if it is not zero by Corollary \ref{thm:main-1}.
Put 
\begin{equation*}
\mc{E}=\bigsqcup_{r\geq 4}\mc{E}_r,
\end{equation*}
where $\mc{E}_r$ is the set of $\e=(\e_1,\dots,\e_r)$ such that
(1) $\e_i\in\{0,1\}$ for $1\leq i\leq r$,
(2) $|\{\,i\,|\,\e_i=0\,\}|\neq |\{\,i\,|\,\e_i=1\,\}|$.
Let $<$ be a partial order on $\mc{E}$ such that for $\e\in \mc{E}_r$ and $\e'\in \mc{E}_{r'}$
\begin{equation}\label{eq:partial order on e}
\e' < \e \ \Longleftrightarrow \ \text{$r'<r$ and $\e'$ is a subsequence of $\e$}.
\end{equation}

\begin{lem}\label{lem:truncation of Verma and R matrix}
Let $\e, \e'\in \mc{E}$ be given such that $\e'<\e$.  
For $(\mb{l},\mb{c})\in\mc{P}^+(\e)$, we have
$$\mf{tr}^\e_{\e'}(R_{\mb{l},\e}(\mb c))=R_{\mb{l},\e'}(\mb c).$$
\end{lem}
\pf If we define $\mf{tr}^\e_{\e'}$ on the weight spaces in $P_{\rm af}$ in the same way as in \eqref{eq:truncation-1}, then we also have $\mf{tr}^\e_{\e'}((\mathcal{W}_{l,\e})_{\rm aff})\cong(\mathcal{W}_{l,\e'})_{\rm aff}$. By \eqref{eq:R matrix on W_l tensor W_m-2}, we have the following commuting diagram:
\begin{equation}\label{eq:affine truncation diagram}
\xymatrixcolsep{5pc}\xymatrixrowsep{3pc}\xymatrix{
((\mathcal{W}_{l_1,\e})_{\rm aff}\ot\cdots \ot (\mathcal{W}_{l_t,\e})_{\rm aff})^{{}^\wedge} \ar@{->}[r]^{\hskip -0.5cm\mc{R}^{\rm norm}_{\mb{l},\e}}\ar@{->}[d]_{\pi^\e_{\e'}} & 
((\mathcal{W}_{l_t,\e})_{\rm aff}\ot\cdots \ot (\mathcal{W}_{l_1,\e})_{\rm aff})^{{}^\wedge}\ar@{->}[d]^{\pi^\e_{\e'}}
\\
((\mathcal{W}_{l_1,\e'})_{\rm aff}\ot \cdots \ot (\mathcal{W}_{l_t,\e'})_{\rm aff})^{{}^\wedge} \ar@{->}[r]^{\hskip -.3cm\mf{tr}^\e_{\e'}(\mc{R}^{\rm norm}_{\mb{l},\e})} & 
((\mathcal{W}_{l_t,\e'})_{\rm aff}\ot\cdots \ot (\mathcal{W}_{l_1,\e'})_{\rm aff})^{{}^\wedge}},
\end{equation}
where $\mc{R}^{\rm norm}_{\mb{l},\e}=\mc{R}^{\rm norm}_{l_1,\dots,l_t}$ \eqref{eq:R matrix for longest permutation} with respect to $\e$.
By Corollary \ref{cor:characterization of R^{norm}}, we have 
\begin{equation*}
\mf{tr}^\e_{\e'}\left(\mc{R}^{\rm norm}_{(l,m),\e}\right)=\mc{R}^{\rm norm}_{(l,m),\e'},
\end{equation*}
where $\mc{R}^{\rm norm}_{(l,m),\e}$ denotes $\mc{R}^{\rm norm}_{l,m}$ \eqref{eq:R matrix on W_l tensor W_m} with respect to $\e$.
Since $\mc{R}^{\rm norm}_{l_1,\dots,l_t}$ is a composition of $\mc{R}^{\rm norm}_{l_i,l_j}$'s and $\mf{tr}^\e_{\e'}$ preserves composition, we have 
\begin{equation}\label{eq:truncation of R matrix}
\mf{tr}^\e_{\e'}(\mc{R}^{\rm norm}_{\mb{l},\e})=\mc{R}^{\rm norm}_{\mb{l},\e'}.
\end{equation}
Specializing \eqref{eq:affine truncation diagram} at $\mb{c}=(c_1,\dots,c_t)$, we have the following commuting diagram:
\begin{equation}\label{eq:finite truncation diagram}
\xymatrixcolsep{4pc}\xymatrixrowsep{3pc}\xymatrix{
\mathcal{W}_{\mb{l},\e}(\mb{c}) \ar@{->}[d]_{\pi^\e_{\e'}}\ar@{->}[r]^{\hskip -5mm R_{\mb{l},\e}(\mb c)} &  \mathcal{W}_{\mb{l}^{\rm rev},\e}(\mb{c}^{\rm rev}) \ar@{->}[d]^{\pi^\e_{\e'}} \\
\mathcal{W}_{\mb{l},\e'}(\mb{c}) \ar@{->}[r]^{\hskip -5mm R_{\mb{l},\e'}(\mb c)} & \mathcal{W}_{\mb{l}^{\rm rev},\e'}(\mb{c}^{\rm rev})
}
\end{equation}
by \eqref{eq:truncation of R matrix}. This implies $\mf{tr}^\e_{\e'}(R_{\mb{l},\e}(\mb c))=R_{\mb{l},\e'}(\mb c)$.
\qed

\begin{thm}\label{thm:truncation sends simple to simple}
Let $\e, \e'\in\mc{E}$ be given such that $\e'<\e$.  
For $(\mb{l},\mb{c})\in\mc{P}^+(\e)$, we have $$\mf{tr}^\e_{\e'}(\mathcal{W}_\e(\mb{l},\mb{c}))\cong\mathcal{W}_{\e'}(\mb{l},\mb{c}).$$
\end{thm}
\pf It follows immediately from Lemma \ref{lem:truncation of Verma and R matrix}.
\qed
\vskip 2mm

Let $(\mb{l},\mb{c})\in \mc{P}^+(\e)$ be given. As a $\ring{\mathcal{U}}(\e)$-module, 
$\mathcal{W}_{\e}(\mb{l},\mb{c})$ is completely reducible and it is a direct sum of $V_\e(\la)$ for $\la\in\cP_{M|N}$. 
Let   
\begin{equation}\label{eq:branching multiplicity}
m_\la^{(\mb{l},\mb{c})}(\e)=\dim_{\Bbbk}{\rm Hom}_{\ring{\mathcal{U}}(\e)}(\mathcal{W}_{\e}(\mb{l},\mb{c}),V_{\e}(\la))
\end{equation}
be the multiplicity of $V_\e(\la)$ in $\mathcal{W}_\e(\mb{l},\mb{c})$ for $\la\in \cP_{M|N}$.

\begin{lem}\label{lem:branching multiplicity}
Let $(\mb{l},\mb{c})\in \mc{P}^+$ and $\la\in \cP$ be given.
If $m_\la^{(\mb{l},\mb{c})}(\e)\neq 0$ for some $\e\in \mc{E}$, then we have 
\begin{equation*}
m_\la^{(\mb{l},\mb{c})}(\e)=m_\la^{(\mb{l},\mb{c})}(\e'),
\end{equation*}
for any $\e'\in \mc{E}$ such that $\e\leq \e'$.
\end{lem}
\pf Since $m_\la^{(\mb{l},\mb{c})}(\e)\neq 0$, we have $\la\in \cP_{M|N}$. 
This implies $\la\in \cP_{M'|N'}$, where $M'$ and $N'$ are the numbers of $0$ and $1$ in $\e'$, respectively. By \eqref{eq:finite truncation diagram} (with $\e$ and $\e'$ exchanged), Proposition \ref{prop:truncation of poly and fund} and Theorem \ref{thm:truncation sends simple to simple}, we have
\begin{equation*}
\dim_{\Bbbk}{\rm Hom}_{\ring{\mathcal{U}}(\e)}(\mathcal{W}_{\e}(\mb{l},\mb{c}),V_{\e}(\la))
=\dim_{\Bbbk}{\rm Hom}_{\ring{\mathcal{U}}(\e')}(\mathcal{W}_{\e'}(\mb{l},\mb{c}),V_{\e'}(\la)).
\end{equation*}
\qed

\begin{thm}\label{thm:branching multiplicity}
For $(\mb{l},\mb{c})\in \mc{P}^+$ and $\la\in \cP$, 
there exists $m^{(\mb{l},\mb{c})}_\la\in\mathbb{Z}_+$ such that 
\begin{itemize}
\item[(1)] if $m^{(\mb{l},\mb{c})}_\la=0$, then $m^{(\mb{l},\mb{c})}_\la(\e)=0$ for all $\e\in \mc{E}$,

\item[(2)] if $m^{(\mb{l},\mb{c})}_\la\neq 0$, then $m^{(\mb{l},\mb{c})}_\la(\e)=m^{(\mb{l},\mb{c})}_\la$ for all $\e\in \mc{E}$ with $m^{(\mb{l},\mb{c})}_\la(\e)\neq 0$.
\end{itemize}
\end{thm}
\pf Define $m^{(\mb{l},\mb{c})}_\la=0$ if $m^{(\mb{l},\mb{c})}_\la(\e)=0$ for all $\e\in \mc{E}$. 
Suppose that both $m^{(\mb{l},\mb{c})}_\la(\e')$ and $m^{(\mb{l},\mb{c})}_\la(\e'')$ are non-zero for some $\e', \e''\in \mc{E}$. Let us take $\e'''\in \mc{E}$ such that $\e'<\e''$ and $\e''<\e'''$. 
By Lemma \ref{lem:branching multiplicity}, we have 
$m^{(\mb{l},\mb{c})}_\la(\e')=m^{(\mb{l},\mb{c})}_\la(\e'')=m^{(\mb{l},\mb{c})}_\la(\e''')$. Hence $m^{(\mb{l},\mb{c})}_\la(\cdot)$ is constant for any $\e\in \mc{E}$ with $m^{(\mb{l},\mb{c})}_\la(\e)\neq 0$, which we denote by $m^{(\mb{l},\mb{c})}_\la$.
\qed
\vskip 2mm

By Theorem \ref{thm:branching multiplicity}, the branching multiplicity \eqref{eq:branching multiplicity} is independent of the choice of $\e$ if it is non-zero. In particular, if $m^{(\mb{l},\mb{c})}_\la(\e)\neq 0$, then 
\begin{equation*}
m^{(\mb{l},\mb{c})}_\la(\e)=m^{(\mb{l},\mb{c})}_\la(\e_{M'|0})=m^{(\mb{l},\mb{c})}_\la(\e_{0|N'}),
\end{equation*}
for some $M'$ or $N'$. Therefore, finding the usual character of $\mathcal{W}_\e(\mb{l},\mb{c})$ is equivalent to finding the one of the corresponding irreducible module over the usual quantum affine algebra of type $A$.

\subsection{Inverse limit category}\label{subsec:inverse limit}
Suppose that 
\begin{equation}\label{eq:e infinite}
\e^\infty=(\e_i)_{i\geq 1}= (\e_1,\dots,\e_k,\dots)
\end{equation} 
is given, where $\e_i\in \{0,1\}$ for $i\geq 1$. We also assume that $\epsilon^\infty$ has infinitely many $1$'s for convenience.
Take an ascending chain $(\epsilon^{(k)})_{k\ge 1}$ of subsequences of $\epsilon^\infty$ in $\mc{E}$ with respect to the partial order \eqref{eq:partial order on e} such that 
\begin{equation}\label{eq:cond-epsilon-infty}
\epsilon^\infty = \lim_k  \epsilon^{(k)}.
\end{equation} 
For example, one may choose $\epsilon^{(k)}=(\epsilon_1, \epsilon_2, \dots, \epsilon_{2k+3})$. Let $M_k,\,N_k$ denote the number of $0$'s, $1$'s in $\epsilon^{(k)}$ respectively. Let us consider the inverse limit category associated to $\{\,\mc{C}(\e^{(k)})\,\}_{k\geq 1}$  (cf.~\cite{EA}, \cite[Section 5.1]{Sc}). We follow the presentation \cite{EA}.
\begin{df}\label{def:inverse limit category}
{\rm 
Let $$\mc{C}(\e^\infty)=\varprojlim \mc{C}(\e^{(k)})$$ be the category defined as follows:
\begin{itemize}
\item[(1)] an object is a pair $\mathbb{V}=((V_k)_{k\geq 1},(f_{k})_{k\geq 1})$ such that 
\begin{equation*}
V_k\in \mc{C}(\e^{(k)}),\quad 
f_k : \mf{tr}^{k+1}_k(V_{k+1})\stackrel{\cong}{\longrightarrow} V_k \quad (k\geq 1),
\end{equation*}
where $\mf{tr}^{k+1}_k=\mf{tr}^{\e^{(k+1)}}_{\e^{(k)}}$ and $f_k$ is an isomorphism of $\mathcal{U}(\epsilon^{(k)})$-modules,

\item[(2)] a morphism from 
$\mathbb{V}=((V_k)_{k\geq 1},(f_{k})_{k\geq 1})$ to 
$\mathbb{W}=((W_k)_{k\geq 1},(g_{k})_{k\geq 1})$ is a sequence $\phi=(\phi_k)_{k\geq 1}$ such that $\phi_k\in {\rm Hom}_{\mathcal{U}(\epsilon^{(k)})}(V_k,W_k)$ and the following diagram is commutative for all $k$:
\begin{equation*}
\xymatrixcolsep{4pc}\xymatrixrowsep{3pc}\xymatrix{
\mf{tr}^{k+1}_k(V_{k+1})\ar@{->}[d]_{\mf{tr}^{k+1}_k(\phi_{k+1})}\ar@{->}[r]^{f_k} & V_k  \ar@{->}[d]^{\phi_k} \\
\mf{tr}^{k+1}_k(W_{k+1}) \ar@{->}[r]^{g_k} & W_k
}
\end{equation*}
\end{itemize}
}
\end{df}\vskip 2mm

Then $\mc{C}(\e^\infty)$ is an abelian category.
Moreover, it is monoidal (or a tensor category) with respect to $\ot$, where we define
\begin{equation*}
\mathbb{V}\ot \mathbb{W} = ((V_k\ot W_k)_{k\geq 1},(f_k\ot g_k)_{k\geq 1}),
\end{equation*}
for $\mathbb{V}=((V_k)_{k\geq 1},(f_{k})_{k\geq 1})$ and 
$\mathbb{W}=((W_k)_{k\geq 1},(g_{k})_{k\geq 1})$. 
Given $k\geq 1$, let 
\begin{equation*}
\mf{tr}_k : \mc{C}(\e^\infty)\longrightarrow \mc{C}(\e^{(k)})
\end{equation*}
be a functor given by $\mf{tr}_k(\mathbb{V})=V_k$ and $\mf{tr}_k(\phi)=\phi_k$ for $\mathbb{V}=((V_k)_{k\geq 1},(f_{k})_{k\geq 1})$ and a morphism $\phi=(\phi_k)_{k\geq 1} : \mathbb{V} \longrightarrow \mathbb{W}$.

The following can be proved by standard arguments.
\begin{lem}\label{lem:universal mapping property of inverse limit}
Given a category $\mathcal{C}$ with a family of functors 
$\{\,F_k:\mathcal{C}\longrightarrow\mc{C}(\epsilon^{(k)})\,\}_{k\geq 1}$
such that $\mf{tr}_{k}^{k+1}\circ F_{k+1} \cong F_{k}$ for $k\geq 1$, there exists a functor 
\begin{equation*}
\xymatrixcolsep{2pc}\xymatrixrowsep{3pc}\xymatrix{
F:=\varprojlim F_k : \mc{C}  \ar@{->}[r] & \ \mc{C}(\epsilon^\infty)},
\end{equation*}
such that $\mf{tr}_{k}\circ F \cong F_{k}$ for all $k\geq 1$.
Moreover, if $F_k$ is exact, then so is $F$. If $\mc{C}$ is a monoidal category and $F_k$ is  monoidal, then so is $F$. 
\end{lem} 
\begin{rem}{\rm For any finite subsequence $\epsilon^\prime$ of $\epsilon^\infty$, there exists $k$ such that $\epsilon^\prime <\epsilon^{(k)}$ so that we may define a truncation functor $\mathfrak{tr}_{\epsilon^\prime}\coloneqq \mathfrak{tr}^{\epsilon^{(k)}}_{\epsilon^\prime} \circ \mathfrak{tr}_{k}$. Here the choice of $k$ is irrelevant because truncation functors are transitive by definition.

If we take another ascending chain $(\widetilde{\epsilon}^{(k)})_{k\ge 1}$ of subsequences of $\epsilon^\infty$ satisfying (\ref{eq:cond-epsilon-infty}) and construct the inverse limit category $\widetilde{\mathcal{C}}(\epsilon^\infty)=\varprojlim \mathcal{C}(\widetilde{\epsilon}^{(k)})$, then we have an equivalence
\[
\mathcal{C}(\epsilon^{\infty}) \cong \widetilde{\mathcal{C}}(\epsilon^\infty)
\] using the truncation functor above and Lemma~\ref{lem:universal mapping property of inverse limit}. Therefore the inverse limit category $\mathcal{C}(\epsilon^\infty)$ does not depend on the choice of an ascending chain, and the requirement  $\epsilon^{(k)}\in \mathcal{E}$ is not a strong one in this sense.
}
\end{rem}

\begin{lem}\label{lem:equiv-infty}
Suppose that there is a family of functors $\left\{F_k : \mathcal{C}\longrightarrow \mathcal{C}(\epsilon^{(k)})\right\}_{k\geq 1}$ satisfying
\begin{itemize}
\item[(1)] $\mathfrak{tr}^{k+1}_{k}\circ F_{k+1}\cong F_{k}$ for $k\geq 1$, 

\item[(2)] there exists $N$ such that $F_n$ is an equivalence of categories for $n\geq N$.
\end{itemize}
Then the resulting functor $F=\varprojlim F_k$ is an equivalence of categories.
In particular, if $\mathfrak{tr}^{n+1}_{n}$ is an equivalence for all $n\geq N$, then $\mathfrak{tr}_{n}$ is an equivalence for $n\geq N$.
\end{lem}
\pf
Take a quasi-inverse $G_k$ of $F_k$ so that $G_{k}\circ\mathfrak{tr}^{k+1}_{k} \cong G_{k+1}$ for $k\ge N$. 
We claim that $G  \coloneqq G_N\circ\mathfrak{tr}_{N}$ is a quasi-inverse of $F$. Clearly $G\circ F\cong \mathrm{id}_{\mathcal{C}}$. 

Given an object $\mathbb{V}$ in $\mc{C}(\e^\infty)$, we have
\[
\mf{tr}_k(F\circ G(\mathbb{V}))=
\mf{tr}_k(F(G_N(\mf{tr}_N(\mathbb{V}))))=F_k(G_N(\mf{tr}_N(\mathbb{V}))).
\]
For $k\geq N$, we have $G_k \cong G_N \circ \mathfrak{tr}^{k}_{N}$ and hence $\mathrm{id}_{\mathcal{C}(\epsilon^{(k)})} \cong F_k \circ G_N \circ \mathfrak{tr}^{k}_{N}$, that is, 
\[
\mf{tr}_k(\mathbb{V}) \cong F_k \circ G_N \circ \mathfrak{tr}^{k}_{N} (\mf{tr}_k(\mathbb{V})) \cong F_k (G_N (\mf{tr}_N(\mathbb{V}))).
\]
For $k<N$, as $F_N \circ G_N \cong \mathrm{id}$,
\[
F_k(G_N(\mf{tr}_N(\mathbb{V})))\cong \mathfrak{tr}^{N}_{k} \circ F_N \circ G_N (\mathfrak{tr}_N (\mathbb{V}))
\cong\mathfrak{tr}_{k}^{N}\circ \mathfrak{tr}_N (\mathbb{V})\cong \mathfrak{tr}_{k}(\mathbb{V}).
\]
From the definition of morphisms in $\mathcal{C}(\epsilon^\infty)$, these isomorphisms are natural. Hence $F\circ G \cong \mathrm{id}_{\mathcal{C}(\epsilon^\infty)}$.
\qed

\begin{lem}\label{lem:simplicity}
Let $\mathbb{V}$ be an object in $\mc{C}(\e^\infty)$.
If $\mathbb{V}$ is non-zero and $\mf{tr}_k(\mathbb{V})$ is irreducible or zero for all $k$,
then $\mathbb{V}$ is a simple object.
\end{lem}
\pf Let $\mathbb{U}$ be a non-zero subobject of $\mathbb{V}$. We may assume that $\mf{tr}_k(\mathbb{U})$ is a $\mathcal{U}(\e^{(k)})$-submodule of $\mf{tr}_k(\mathbb{V})$. Let $k_0$ (resp. $k_1$) be the smallest one such that $\mf{tr}_k(\mathbb{U})\neq 0$ (resp. $\mf{tr}_k(\mathbb{V})\neq 0$). 
Since $\mf{tr}_k(\mathbb{V})$ is irreducible for $k\geq k_1$ and $\mf{tr}_k(\mathbb{U})$ is a submodule of $\mf{tr}_k(\mathbb{V})$, we have $k_0\geq k_1$, and $\mf{tr}_k(\mathbb{V})=\mf{tr}_k(\mathbb{U})$ for $k\geq k_0$, which also implies that $\mf{tr}^{k_0}_k(\mf{tr}_{k_0}(\mathbb{U}))=\mf{tr}^{k_0}_k(\mf{tr}_{k_0}(\mathbb{V}))$ for $k\leq k_0$. Therefore, we have $k_0=k_1$ and $\mathbb{U}$ is isomorphic to $\mathbb{V}$.
\qed\vskip 2mm

\begin{rem}{\rm
We expect that the category $\mc{C}(\e^\infty)$ when $\e$ is homogeneous is closely related to the representations of the quantum affine algebra $U_q(\widehat{\mf{sl}}_\infty)$ introduced in \cite{Her11}.
}
\end{rem}

\subsection{$R$ matrix and simple objects in $\mc{C}(\e^\infty)$}
For $l\in\mathbb{Z}_+$ and $x\in \Bbbk^\times$, we let 
\begin{equation*}
\mathcal{W}_{l,\e^\infty}(x) = \left((\mathcal{W}_{l,\e^{(k)}}(x))_{k\geq 1},(f_k)_{k\geq 1}\right),
\end{equation*}
where $f_k$ is an isomorphism in Proposition \ref{prop:truncation of poly and fund}(2). Then it is an object in $\mc{C}(\e^\infty)$. Define its affinization by
\begin{equation*}
\begin{split}
(\mathcal{W}_{l,\e^\infty})_{\rm aff} &= 
\left(((\mathcal{W}_{l,\e^{(k)}})_{\rm aff})_{k\geq 1},(\mathrm{id}_{\Bbbk[z^{\pm 1}]}\ot f_k)_{k\geq 1}\right).\\
\end{split}
\end{equation*}
By \eqref{eq:affine truncation diagram} and \eqref{eq:truncation of R matrix}, there exists a well-defined morphism (in the sense of Definition~\ref{def:inverse limit category}(2))
\begin{equation}\label{R^{norm} at infty}
\xymatrixcolsep{4pc}\xymatrixrowsep{3pc}\xymatrix{
((\mathcal{W}_{l,\e^\infty})_{\rm aff} \ot (\mathcal{W}_{m,\e^\infty})_{\rm aff})^{{}^\wedge}\ \ar@{->}[r]^{\mc{R}^{\rm norm}_{(l,m),\e^\infty}} &  \
((\mathcal{W}_{m,\e^\infty})_{\rm aff} \ot (\mathcal{W}_{l,\e^\infty})_{\rm aff})^{{}^\wedge},
}
\end{equation}
for $l,m\in\mathbb{Z}_+$, given by 
$\mc{R}^{\rm norm}_{(l,m),\e^\infty}=\left(\mc{R}^{\rm norm}_{(l,m),\e^{(k)}}\right)_{k\geq 1}$,
where $\mc{R}^{\rm norm}_{(l,m),\e^{(k)}}$ is the map given in \eqref{eq:R matrix on W_l tensor W_m-2} with respect to $\e^{(k)}$.

For $({\mb l},{\mb c})\in \mc{P}^+$, 
let $\mathcal{W}_{\mb{l},\e^\infty}(\mb{c})=\mathcal{W}_{l_1,\e^\infty}(c_1)\ot\dots\ot \mathcal{W}_{l_t,\e^\infty}(c_t)$. 
By using \eqref{R^{norm} at infty}, we have a morphism 
\begin{equation*}
\xymatrixcolsep{2pc}\xymatrixrowsep{3pc}\xymatrix{
R_{\mb{l},\e^\infty}(\mb c) : \mathcal{W}_{\mb{l},\e^\infty}(\mb{c})\ \ar@{->}[r] & \ \mathcal{W}_{\mb{l}^{\rm rev},\e^\infty}(\mb{c}^{\rm rev})},
\end{equation*}
given by $R_{\mb{l},\e^\infty}(\mb c)=(R_{\mb{l},\e^{(k)}}(\mb c))_{k\geq 1}$. 
Let 
\begin{equation*}
\mathcal{W}_{\e^\infty}(\mb{l},\mb{c})={\rm Im}R_{\mb{l},\e^{\infty}}(\mb{c}),
\end{equation*}  
which is a well-defined object in $\mc{C}(\e^\infty)$ by Theorem \ref{thm:truncation sends simple to simple}.
Then we have the following analogue of Theorem \ref{thm:main-1}:
\begin{thm}\label{thm:suff-cond-nonzero}
For $(\mb{l},\mb{c})\in \mc{P}^+$, $\mathcal{W}_{\e^\infty}(\mb{l},\mb{c})$ is a nonzero simple object in $\mc{C}(\e^\infty)$. 
\end{thm}
\pf We have
$\mathcal{W}_{\e^\infty}(\mb{l},\mb{c})
=\left(\left(\mathcal{W}_{\e^{(k)}}(\mb{l},\mb{c})\right)_{k\geq 1},(h_k)_{k\geq 1}\right)$, where $(h_k)_{k\geq 1}$ is induced from the one associated to $\mathcal{W}_{\mb{l}^{\rm rev},\e^\infty}(\mb{c}^{\rm rev})$.
Suppose that $\mathcal{W}_{\e^\infty}(\mb{l},\mb{c})$ is not zero, and let $k_0$ be the smallest one such that $\mathcal{W}_{\e^{(k_0)}}(\mb{l},\mb{c})$ is not zero.
Since $\mathcal{W}_{\e^{(k)}}(\mb{l},\mb{c})$ is simple for $k\geq k_0$ by Theorem \ref{thm:main-1}, $\mathcal{W}_{\e^\infty}(\mb{l},\mb{c})$ is simple by Lemma \ref{lem:simplicity}. 

Next, we claim that $\mathcal{W}_{\e^\infty}(\mb{l},\mb{c})$ is non-zero for all $(\mb{l},\mb{c})\in\mc{P}^+$. 
Take a sufficiently large $k$ such that $l_i< N_k$ for $1\leq i\leq t$, and put $\e'=\e_{0|N_k}$. 
Since $\e'<\e^{(k)}$, we have by Theorem \ref{thm:truncation sends simple to simple} that 
\begin{equation*}
\mf{tr}^{\e^{(k)}}_{\e'}(\mathcal{W}_{\e^{(k)}}(\mb{l},\mb{c}))=\mathcal{W}_{\e'}(\mb{l},\mb{c}).
\end{equation*}  
Since $\mathcal{W}_{\e'}(\mb{l},\mb{c})\neq 0$ by our choice of $\e^{(k)}$ and \cite{AK,Kas02}, we have $\mathcal{W}_{\e^{(k')}}(\mb{l},\mb{c})\neq 0$ for all $k'\geq k$. Hence  $\mathcal{W}_{\e^\infty}(\mb{l},\mb{c})$ is non-zero.
\qed
\vskip 2mm

We define a category $\ring{\mc{C}}(\e^\infty)$ in a similar way, where we replace $\mathcal{U}(\e^{(k)})$ and $\mc{C}(\e^{(k)})$ with $\ring{\mathcal{U}}(\e^{(k)})$ and $\ring{\mc{C}}(\e^{(k)})$, respectively. 
For $\la\in \cP$, let
\begin{equation*}
V_{\e^\infty}(\la)=\left( (V_{\e^{(k)}}(\la) )_{k\geq 1},(\ring{f}_k)_{k\geq 1}\right),
\end{equation*}
where $\ring{f}_k$ is the isomorphism in Proposition \ref{prop:truncation of poly and fund}(1), and we assume that $V_{\e^{(k)}}(\la) =0$ if $\la\not\in \cP_{M_k|N_k}$.
By Lemma \ref{lem:simplicity}, $V_{\e^\infty}(\la)$ is a simple object in $\ring{\mc{C}}(\e^\infty)$. 

\begin{thm}
For $(\mb{l},\mb{c})\in \mc{P}^+$, $\mathcal{W}_{\e^\infty}(\mb{l},\mb{c})$ is a semisimple object in $\ring{\mc{C}}(\e^\infty)$, and  
\begin{equation*}
\mathcal{W}_{\e^\infty}(\mb{l},\mb{c}) = \bigoplus_{\la\in\cP}V_{\e^\infty}(\la)^{\oplus m^{(\mb{l},\mb{c})}_\la},
\end{equation*}
where $\la$ runs over the partitions of $l_1+\dots+ l_t$ for $\mb{l}=(l_1,\dots,l_t)$, and $m^{(\mb{l},\mb{c})}_\la$ is the multiplicity in Theorem \ref{thm:branching multiplicity}.
\end{thm}
\pf We have as a $\ring{\mathcal{U}}(\e^{(k)})$-module
\begin{equation*}
\mathcal{W}_{\e^{(k)}}(\mb{l},\mb{c}) = \bigoplus_{\la\in\cP_{M_k|N_k}}V_{\e^{(k)}}(\la)^{\oplus m^{(\mb{l},\mb{c})}_\la(\e^{(k)})} \quad (k\geq 1).
\end{equation*}
Hence the decomposition of $\mathcal{W}_{\e^\infty}(\mb{l},\mb{c})$ follows from Theorem \ref{thm:branching multiplicity} (cf.~Remark \ref{rem:isomorphism on truncated space}).
\qed

\section{Generalized quantum affine Schur-Weyl duality functor}\label{sec:QSWD}
\subsection{Quiver Hecke algebras}
Let us briefly review necessary background for (symmetric) quiver Hecke algebras following the convention in \cite{KKK,KKKO15} with a slight modification of notations.

Let ${\bf k}$ be a field.
Let $J$ be a set, $\mathbb{Z}[J]$ the free abelian group generated by $J$, and  
 $\mathbb{N}[J]$  the subset of $\mathbb{Z}[J]$ consisting of $\sum_{i\in J}c_{i}\cdot i$ with $c_i\in\mathbb{Z}_+$. 
For $\ell\ge 2$, we assume that the symmetric group $\mf{S}_\ell$, which is generated by $s_{i}=(i,i+1)$ for $i=1,\dots,\ell-1$, acts on $J^\ell$ by place permutation.

For $\beta=\sum_{i\in J}c_{i}\cdot i\in \mathbb{N}[J]$ with ${\rm ht}(\beta):=\sum_{i\in J}c_i=\ell$, set
\[
J^{\beta}=\left\{\, \nu=(\nu_{1},\dots,\nu_{\ell})\in J^{\ell}\,|\,\nu_{1}+\cdots+\nu_{
\ell}=\beta\,\right\}.
\]
Let $u, v$ be formal variables. Let $(Q_{ij})_{i,j\in J}$ be the matrix with entries 
$Q_{ij}=Q_{ij}(u,v)\in\mathbf{k}[u,v]$ such that $Q_{ij}(u,v)=Q_{ji}(v,u)$ for $i\neq j$, and $Q_{ii}(u,v)=0$.

Suppose that $A=(a_{ij})_{i,j\in J}$ is a symmetrizable generalized Cartan matrix. Then there exist positive integers $s_i$ ($i\in J$) such that $s_ia_{ij}=s_ja_{ji}=:{\bf a}_{ij}$ for $i,j\in J$.
We further assume that the coefficient of $u^{p}v^{q}$ ($p,q\in\mathbb{Z}_+$) in $Q_{ij}(u,v)$ vanishes unless ${\bf a}_{ii}p+{\bf a}_{jj}q=-2 {\bf a}_{ij}$.

\begin{df}
{\rm
The {\em quiver Hecke algebra} or {\em Khovanov-Lauda-Rouquier algebra} $R(\beta)$ at $\beta$ associated with $(Q_{ij})_{i,j\in J}$ is the associative $\bf{k}$-algebra with $1$ generated
by $e(\nu)$ $(\nu\in J^\beta)$, $x_k$ $(k=1,\dots,\ell)$, and $\tau_m$ $(m=1,\dots,\ell-1)$ subject to the following relations:
{\allowdisplaybreaks
\begin{align*}
& e(\nu)e(\nu^{\prime})=\delta_{\nu\nu^{\prime}}e(\nu),\quad 
\sum_{\nu\in J^{\beta}}e(\nu)=1,\\
& x_{k}x_{k^{\prime}}=x_{k^{\prime}}x_{k},\quad 
x_{k}e(\nu)=e(\nu)x_{k},\quad \tau_{m}e(\nu)=e(s_{m}(\nu))\tau_{m},\\
& \left(\tau_{m}x_{k}-x_{s_{m}(k)}\tau_{m}\right)e(\nu)=
\begin{cases}
-e(\nu) & \text{if $k=m$ and $\nu_{m}=\nu_{m+1}$},\\
e(\nu) & \text{if $k=m+1$ and $\nu_{m}=\nu_{m+1}$},\\
0 & \text{otherwise,}
\end{cases}\label{eq:x and tau commutation} \\
& \tau_{m}^{2}e(\nu)=Q_{\nu_{m},\nu_{m+1}}(x_{m},x_{m+1})e(\nu),\quad \tau_{m}\tau_{m^{\prime}}=\tau_{m^{\prime}}\tau_{m} \quad  \text{ if }\left|m-m^{\prime}\right|>1,\\
&\left(\tau_{m+1}\tau_{m}\tau_{m+1}-\tau_{m}\tau_{m+1}\tau_{m}\right)e(\nu)=\\ & \quad\quad\quad\quad
\begin{cases}
\dfrac{Q_{\nu_{m},\nu_{m+1}}(x_{m},x_{m+1})-Q_{\nu_{m},\nu_{m+1}}(x_{m+2},x_{m+1})}{x_{m}-x_{m+2}}e(\nu) & \text{if }\nu_{m}=\nu_{m+2},\\
0 & \text{otherwise}.
\end{cases}
\end{align*}}}
\end{df}

The algebra $R(\beta)$ carries a $\mathbb{Z}$-grading defined by $\deg e(\nu)=0$, $\deg x_{k}e(\nu)=2$, and $\deg\tau_{m}e(\nu)=-a_{\nu_m\nu_{m+1}}$.
We set 
\begin{equation}\label{eq:R}
\begin{split}
R(\ell)&=\bigoplus_{\substack{\beta\in J^\ell, \mathrm{ht}(\beta)=\ell}}R(\beta) \quad (\ell\ge 1),\ \ R=\bigoplus_{\ell\ge 0}R(\ell),
\end{split}
\end{equation}
where $R(0)={\bf k}$.
For $\beta\in \mathbb{N}[J]$ with $\mathrm{ht}(\beta)=\ell$, we have the following decomposition:
\[
R(\beta)=\bigoplus_{\nu\in J^{\beta},w\in \mf{S}_{\ell}}\mathbf{k}\left[x_{1},\dots,x_{\ell}\right]e(\nu)\tau_{w},
\]
where $\tau_{w}=\tau_{i_{1}}\cdots\tau_{i_{l}}$
is defined after fixing a reduced expression $w=s_{i_{1}}\cdots s_{i_{l}}$
for each $w\in \mf{S}_{\ell}$.

For $\beta_1, \beta_2\in \mathbb{N}[J]$ with ${\rm ht}(\beta_1)=\ell_1$ and ${\rm ht}(\beta_2)=\ell_2$, let 
\begin{equation*}
e(\beta_1,\beta_2)=\sum_{\nu}e(\nu)\in R(\beta_1+\beta_2),
\end{equation*}
where the sum is over $\nu=(\nu_1,\dots,\nu_{\ell_1+\ell_2})\in J^{\beta_1+\beta_2}$ such that $(\nu_1,\dots,\nu_{\ell_1})\in J^{\beta_1}$ and $(\nu_{\ell_1+1},\dots,\nu_{\ell_1+\ell_2})\in J^{\beta_2}$. Then we have a homomorphism of ${\bf k}$-algebras 
\begin{equation*}
\xymatrixcolsep{2pc}\xymatrixrowsep{3pc}\xymatrix{
R(\beta_1)\ot R(\beta_2) \ \ar@{->}[r] & \ e(\beta_1,\beta_2)R(\beta_1+\beta_2) e(\beta_1,\beta_2),}
\end{equation*}
given by sending
\begin{gather*}
e(\nu)\ot e(\eta) \longmapsto e(\nu \ast \eta),\\
x_{k_1}\ot 1  \longmapsto  x_{k_1} e(\beta_1,\beta_2),\quad 
1\ot x_{k_2}   \longmapsto x_{\ell_1+k_2}e(\beta_1,\beta_2), \\
\tau_{m_1}\ot 1 \longmapsto  \tau_{m_1}e(\beta_1,\beta_2),\quad 
1\ot\tau_{m_2} \longmapsto \tau_{\ell_1+m_2} e(\beta_1,\beta_2),
\end{gather*}
for $1\le k_i\le \ell_i$ and $1\le m_i\le \ell_i-1$ ($i=1,2$),
where $\nu\ast \eta=(\nu_1,\dots,\nu_{\ell_1},\eta_1,\dots,\eta_{\ell_2})$ is the concatenation of $\nu=(\nu_1,\dots,\nu_{\ell_1})$ and $\eta=(\eta_1,\dots,\eta_{\ell_2})$. 
Let $M_i$ be an $R(\beta_i)$-module for $i=1,2$. The convolution product of $M_1$ and $M_2$ is defined by
\begin{equation*}
M_1\circ M_2 = R(\beta_1+\beta_2)\stackrel[R(\beta_1)\ot R(\beta_2)]{}{\otimes} (M_1\ot M_2),
\end{equation*}  
which is an $R(\beta_1+\beta_2)$-module.

Let $R(\beta)$-gmod denote the category of finite-dimensional graded $R(\beta)$-modules for $\beta\in \mathbb{N}[J]$, and let $R\text{-gmod}=\bigoplus_{\beta\in\mathbb{N}[J]}R(\beta)\text{-gmod}$.
Then $R\text{-gmod}$ is monoidal with respect to the convolution product with the unit $R(0)={\bf k}$. Moreover, there is
a degree shifting functor $q$ on $R\text{-gmod}$ defined by 
$(qM)_{k}=M_{k-1}$ for $M=\bigoplus_{k\in\mathbb{Z}}M_k$. 

Next, let us recall the notion of (renormalized) $R$ matrix for quiver Hecke algebras.
Let $\beta\in \mathbb{N}[J]$ be given with ${\rm ht}(\beta)=\ell\ge 2$.
For $1\leq m\leq \ell-1$, let 
\begin{equation*}
\varphi_me(\nu)=
\begin{cases}
(\tau_mx_m-x_m\tau_m)e(\nu) & \text{if $\nu_m=\nu_{m+1}$},\\
\tau_me(\nu) & \text{if $\nu_m\neq \nu_{m+1}$},
\end{cases}
\end{equation*}
for $\nu=(\nu_1,\dots,\nu_\ell)\in J^\beta$.
For $\ell_1, \ell_2\ge 1$ with $\ell_1+\ell_2=\ell$, put
\begin{equation*}
\varphi_{[\ell_1,\ell_2]}=
\varphi_{k_1}\dots\varphi_{k_{\ell_1\ell_2}},
\end{equation*}
where $w=s_{k_1}\dots s_{k_{\ell_1\ell_2}}\in \mf{S}_{\ell}$ is a reduced expression of the permutation $w$ such that $w(a)=a+\ell_2$ for $1\le a\le \ell_1$ and $w(b)=b-\ell_1$ for $\ell_1+1\le b\le \ell$. 

Suppose that $\beta_1, \beta_2\in \mathbb{N}[J]$ are given such that $\beta=\beta_1+\beta_2$ and ${\rm ht}(\beta_i)=\ell_i$ ($i=1,2$).
Let $M_i$ be an $R(\beta_i)$-module for $i=1,2$. 
Then we have a homomorphism of $R(\beta)$-modules:
\begin{equation*}
\xymatrixcolsep{2pc}\xymatrixrowsep{0pc}\xymatrix{
R_{M_1,M_2} : M_1\circ M_2 \ \ar@{->}[r] & \ M_2\circ M_1 \\
\quad\quad\quad\  a\ot(u_1\ot u_2) \ \ar@{|->}[r] & \ a\varphi_{[\ell_2,\ell_1]}\ot(u_2\ot u_1)}
\end{equation*}
for $a\in R(\beta)$ and $u_i\in M_i$ $(i=1,2)$.
The renormalized $R$ matrix on $M_1\circ M_2$ can be defined by using an affinization of $M_i$'s. More precisely, for an indeterminate $z$ and $R(\gamma)$-module $N$, let $N_z:={\bf k}[z]\ot N$ be the ${\bf k}[z]\ot R(\gamma)$-module defined by
\begin{equation*}
\begin{split}
& e(\nu)(a\ot u)= a\ot (e(\nu)u),\quad x_k(a\ot u)= (za)\ot u + a \ot (x_ku),\quad \tau_m(a\ot u) = a\ot (\tau_m u),
\end{split}
\end{equation*} 
for $a\ot u \in {\bf k}[z]\ot N$. 
Then we define ${\bf r}_{M_1,M_2} : M_1\circ M_2 \longrightarrow M_2\circ M_1$ by
\begin{equation*}
\xymatrixcolsep{2pc}\xymatrixrowsep{0pc}\xymatrix{
{\bf r}_{M_1,M_2}=z^{-s}R_{(M_1)_z,M_2}\vert_{z=0}},
\end{equation*}
where $s$ is the largest non-negative integer such that 
${\rm Im}R_{(M_1)_z,M_2}\subset z^s R_{M_2,(M_1)_z}$. We call ${\bf r}_{M_1,M_2}$ the renormalized $R$ matrix on $M_1\circ M_2$. 

Note that the renormalized $R$ matrices satisfy the Yang-Baxter equation. In particular, this implies that given $\beta_i\in \mathbb{N}[J]$ and $R(\beta_i)$-module $M_i$ ($i=1,\dots,t$), 
we have a well-defined homomorphism of $R(\beta_1+\dots+\beta_t)$-modules
\begin{equation}\label{eq:renorm R matrix for QHA}
\xymatrixcolsep{2pc}\xymatrixrowsep{0pc}\xymatrix{
{\bf r}_{w} : M_1\circ\dots\circ M_t \ \ar@{->}[r] & \ M_{w(1)}\circ \dots\circ M_{w(t)}},
\end{equation}
for each $w\in \mf{S}_t$, which is given by a composition of ${\bf r}_{M_i,M_j}$'s with respect to a reduced expression of $w$.

\subsection{Quiver Hecke algebras of type $A$}
Let $\Gamma$ be the quiver given by

\begin{equation*}\label{eq:quiver Gamma}
\xymatrixcolsep{2.5pc}\xymatrixrowsep{0pc}\xymatrix{
\cdots   \ar@{->}[r] &  \circ \ar@{->}[r] & \circ \ar@{->}[r]  & \circ \ar@{->}[r] & \circ \ar@{->}[r] &   \circ \ar@{->}[r] & \cdots \\
& {}_{-2}  & {}_{-1} & {}_0 & {}_{1} & {}_2 & & }.
\end{equation*}
Suppose that $J=\mathbb{Z}$. Let $P_{ij}(u,v)=(u-v)^{d_{ij}}$, where $d_{ij}$ is the number of arrows from $i$ to $j$ in $\Gamma$, and take $(Q_{ij})_{i,j\in J}$ as 
\begin{equation}\label{eq:Q polynomials for Gamma}
Q_{ij}=Q_{ij}(u,v)=P_{ij}(u,v)P_{ji}(v,u) \quad (i\neq j).
\end{equation}
The resulting quiver Hecke algebra $R^J(\beta)$ at $\beta\in \mathbb{N}[J]$ is associated to $(Q_{ij})_{i,j\in J}$ and the generalized Cartan matrix $(a_{ij})_{i,j\in\mathbb{Z}}$ with 
\[
a_{ij}=\begin{cases}
2 & \text{if } i=j, \\
-1 & \text{if } i=j\pm 1, \\
0 & \text{otherwise.}
\end{cases}
\]

One may regard a pair of integers $(a,b)$ with $a\leq b$ as a positive root of type $A_\infty$, say $\beta_{(a,b)}:=1\cdot a + 1\cdot(a+1) + \dots + 1\cdot b\in \mathbb{N}[J]$.
Let us call such $(a,b)$ a segment and
$\ell=b-a+1$ its length. Define the lexicographic order
on the set of segments, that is,
\[
(a,b)\geq (a^\prime,b^\prime) \Longleftrightarrow a>a^\prime \text{ or } (a=a^\prime, b\geq b^\prime).
\]
We call a finite sequence of segments $\left( (a_1,b_1),\dots,(a_t,b_t)\right)$ a multisegment, and say that it is ordered if $(a_k,b_k)\geq (a_{k+1},b_{k+1})$ for all $1\le k\le t-1$.

For a segment $(a,b)$ of length $\ell$, let $L(a,b)=\mathbf{k} u(a,b)$
be a one-dimensional $R^J\left(\beta_{(a,b)}\right)$-module given
by
\begin{align*}
x_{k}u(a,b) & =\tau_{m}u(a,b)=0,\quad
e(\nu)u(a,b)  =\begin{cases}
u(a,b) & \text{if }\nu=(a,a+1,\dots,b),\\
0 & \text{otherwise},
\end{cases}
\end{align*}
for $1\le k\le \ell$, $1\le m\le \ell-1$, and $\nu\in J^{\beta_{(a,b)}}$
We also denote $L(a,a)$ by $L(a)$.

Fix $\ell\ge 2$. 
Suppose that $\left( (a_1,b_1),\dots,(a_t,b_t)\right)$ is an ordered multisegment with $\ell_k$ the length of $(a_k,b_k)$ and $\sum_{k}\ell_k=\ell$.
By \eqref{eq:renorm R matrix for QHA}, 
we have a homomorphism of $R^J(\beta)$-modules
\begin{equation}\label{eq:maximal intertwiner for QHA}
\xymatrixcolsep{2pc}\xymatrixrowsep{0pc}\xymatrix{
{\bf r}_{w_0} : L(a_1,b_1)\circ\dots\circ L(a_t,b_t)  \ar@{->}[r] & \ L(a_t,b_t)\circ\dots\circ L(a_1,b_1)},
\end{equation}
where $\beta=\sum_{i=1}^t\beta_{(a_i,b_i)}$ and $w_0$ is the longest element in $\mf{S}_t$.

\begin{prop}[{\cite{KP},\cite[Proposition 4.7]{KKK}}]
\label{prop:classification of simples for QHA}
Under the above hypothesis, we have
\begin{itemize}
\item[(1)] the correspondence from $\left( (a_1,b_1),\dots,(a_t,b_t)\right)$ to $\mathrm{hd}\left(L(a_{1},b_{1})\circ\cdots\circ L(a_{t},b_{t})\right)$, the head of $L(a_{1},b_{1})\circ\cdots\circ L(a_{t},b_{t})$, is a bijection from the set of ordered multisegments to the set of isomorphism classes (up to grading shifts) of finite-dimensional irreducible graded $R^J(\ell)$-modules,

\item[(2)] the image of ${\bf r}_{w_0}$ in \eqref{eq:maximal intertwiner for QHA} is irreducible and isomorphic to $\mathrm{hd}\left(L(a_{1},b_{1})\circ\cdots\circ L(a_{t},b_{t})\right)$.
\end{itemize}
\end{prop}

Finally we recall some particular exact sequences involving renormalized $R$-matrices.

\begin{prop}[{\cite[Proposition 4.3]{KKK}}]\label{prop:QHAexact}
 Let $(a,b), (a^{\prime},b^{\prime})$ be two segments with lengths $\ell$ and $\ell'$ respectively such that $(a,b)\geq(a^{\prime},b^{\prime})$. 
\begin{enumerate}
\item If one of the following holds: $a^{\prime}<a\leq b\leq b^{\prime}$, $a>b^{\prime}+1$, $a=a^{\prime}\leq b^{\prime}\leq b$, then
$L(a,b)\circ L(a^{\prime},b^{\prime})$ is irreducible and
\[
\xymatrixcolsep{2pc}\xymatrixrowsep{0pc}\xymatrix{
L(a,b)\circ L(a^{\prime},b^{\prime}) \ar@{->}[r]^{\mathbf{r}} & L(a^{\prime},b^{\prime})\circ L(a,b)
}
\]
is an isomorphism.
\item If $a^{\prime}<a\leq b^{\prime}<b$, then we have the following exact
sequence 
$$
\begin{tikzcd}[row sep=small]
0 \arrow[r] & L(a^{\prime},b)\circ L(a,b^{\prime}) \arrow[r] \arrow[rd, phantom, ""{coordinate, name=Z}] & L(a,b)\circ L(a^{\prime},b^{\prime})
\arrow[dl, rounded corners,  "\hskip -3.5cm\mathbf{r}",
to path={ -- ([xshift=2ex]\tikztostart.east)
|- (Z) [near end]\tikztonodes
-| ([xshift=-2ex]\tikztotarget.west)
-- (\tikztotarget)}] \\
& L(a^{\prime},b^{\prime}) \circ L(a,b) \arrow[r] & L(a^\prime,b)\circ L(a,b^\prime) \arrow[r] & 0
\end{tikzcd}.
$$

\item If $a=b^{\prime}+1$, then we have an exact sequence 
$$
\begin{tikzcd}[column sep=small]
0\arrow[r] 
& L(a^\prime,b) \arrow[r]
& L(a,b)\circ L(a^\prime, b^\prime) \arrow[r,"\mathbf{r}"]
& L(a^\prime, b^\prime)\circ L(a,b) \arrow[r]
& L(a^\prime, b) \arrow[r]
& 0
\end{tikzcd}.
$$
\end{enumerate}
Here we ignore the grading and denote by $\mathbf{r}$ the renormalized $R$-matrix on $L(a,b)\circ L(a',b')$.
\end{prop}

\subsection{Generalized quantum affine Schur-Weyl duality functor}\label{subsec:QASWD}
Let $R^J(\beta)$ be the quiver Hecke algebra at $\beta\in \mathbb{N}[J]$ associated to $(Q_{ij})_{i,j\in J}$ in \eqref{eq:Q polynomials for Gamma} and let $R^J(\ell)$ and $R^J$ denote the corresponding ones in \eqref{eq:R}.
The goal of this subsection is to show that there is a generalized quantum affine Schur-Weyl duality functor $\mc{F}_\e$ from $R^J$-gmod to $\mc{C}(\e)$, which can be constructed by the same method as in \cite[Section 3]{KKK} when $\e=\e_{0|N}$.

Suppose that ${\bf k}=\Bbbk$ is the algebraic closure of $\mathbb{Q}(q)$ in $\bigcup_{m>0}\mathbb{C}(\!(q^{\frac{1}{m}})\!)$.
Let $X: J\rightarrow\mathbf{k}^{\times}$ be given by $X(i)=q^{-2i}$ for $i\in J$. 
Note that $X(i)/X(j)=q^{-2(i-j)}$ is a zero of the denominator of $d_{1,1}(z)=z-q^{2}$ if and only if $j=i+1$ for $i\in J$.

Fix $\ell\ge 2$. Let $X_1,\dots,X_\ell$ be indeterminates.
For $\nu=(\nu_1,\dots,\nu_\ell)\in J^\ell$,
put $X(\nu)=(X(\nu_1),\dots,X(\nu_\ell))\in ({\bf k}^\times)^\ell$ and
\begin{equation*}
{\mathbb{O}}_{\nu}=\mathbf{k}\left\llbracket X_{1}-X(\nu_{1}),\dots,X_{\ell}-X(\nu_{\ell})\right\rrbracket.
\end{equation*}

For $\beta\in \mathbb{N}[J]$ with ${\rm ht}(\beta)=\ell$, we define $\mathbb{O}_{\beta}$
to be a ${\bf k}$-algebra with a ${\bf k}$-basis $\{\,fe(\nu)\,|\, f\in \mathbb{O}_{\nu}\, (\nu\in J^\beta)\,\}$, where $e(\nu)$ commutes with $\mathbb{O}_{\nu}$ and  $e(\nu)e(\nu')=\delta_{\nu\nu'}e(\nu)$ for $\nu, \nu'\in J^\beta$. Similarly, we define $\mathbb{K}_{\beta}$, where $\mathbb{O}_{\nu}$ is replaced by its field of quotients, say $\mathbb{K}_{\nu}$. We have
\begin{equation}\label{eq:O and K}
\begin{split}
\mathbb{O}_{\beta}=\bigoplus_{\nu\in J^\beta}\mathbb{O}_{\nu}e(\nu),\ \
\mathbb{K}_{\beta}=\bigoplus_{\nu\in J^\beta}\mathbb{K}_{\nu}e(\nu).
\end{split}
\end{equation}

Let ${\bf k}[\mf{S}_\ell]$ denote the group algebra of $\mf{S}_\ell$ over ${\bf k}$ with a basis $\left\{r_w\,|\,w\in \mathfrak{S}_{\ell}\right\}$.
We assume that $\mf{S}_\ell$ acts on $\mathbb{K}_\beta$ by $w(X_i)=X_{w^{-1}(i)}$ for $w\in \mf{S}_\ell$ and $1\le i\le \ell$, and hence $\mathbb{K}_\beta\ot_{\bf k}{\bf k}[\mf{S}_\ell]$ is a ${\bf k}$-algebra with $r_w f =w(f) r_w$ for $f\in \mathbb{K}_\beta$ and $w\in \mf{S}_\ell$.
Recall that
we have an embedding of ${\bf k}$-algebras
\begin{equation}\label{eq:embedding of QHA}
\xymatrixcolsep{3pc}\xymatrixrowsep{0pc}\xymatrix{
R^J(\beta)  \ \ar@{^{(}->}[r] & \ \mathbb{K}_\beta\ot_{\bf k}{\bf k}[\mf{S}_\ell]
},
\end{equation}
or we may identify $R^J(\beta)$ as a subalgebra of 
$\mathbb{K}_\beta\ot_{\bf k}{\bf k}[\mf{S}_\ell]$ by letting

\begin{equation}\label{eq:poly-repn-QHA}
\begin{split}
e(\nu) &= e(\nu),\\
e(\nu) x_k &=e(\nu)X(\nu_k)^{-1}(X_k-X(\nu_k)),\\
e(\nu)\tau_m &=
\begin{cases}
e(\nu)(r_{m}-1) \left( \dfrac{1}{x_m-x_{m+1}} \right)  & \text{if $\nu_m=\nu_{m+1}$},\\ 
e(\nu) r_m  P_{\nu_m,\nu_{m+1}}(x_{m+1},x_{m}) & \text{if $\nu_m\neq \nu_{m+1}$},
\end{cases}
\end{split}
\end{equation}
for $1\le k\le \ell$ and $1\leq m\leq \ell-1$, where $r_m=r_{s_m}$ \cite[Proposition 3.12]{R}.

Now, suppose that $\e$ is given. 
Let $V=(\mathcal{W}_{1,\epsilon})_{\mathrm{aff}}$ and regard $V^{\ot \ell}$ as a $\mathbf{k}[X_{1}^{\pm1},\dots,X_{\ell}^{\pm1}]\ot\mathcal{U}(\e)$-module,
where $X_{i}$ acts as $z$ on the $i$th component \eqref{eq:automorphism z}.
Put
\begin{equation*}
\begin{split} 
{V}^{\ot \beta}_{\mathbb{O}}= \mathbb{O}_\beta\ot_K V^{\ot \ell},\quad
{V}^{\ot \beta}_{\mathbb{K}}= {\mathbb{K}_\beta}\ot_{\mathbb{O}_\beta} {V}^{\ot \beta}_{\mathbb{O}},
\end{split}
\end{equation*}
where $K={\bf k}[X_{1}^{\pm1},\dots,X_{\ell}^{\pm1}]$.
We have
\begin{equation*}
\begin{split}
{V}^{\ot \beta}_{\mathbb{O}}=\bigoplus_{\nu\in J^{\beta}}{V}^{\ot \nu}_{\mathbb{O}},
\quad\text{where}\quad {V}^{\ot \nu}_{\mathbb{O}}
={\mathbb{O}}_{\nu}e(\nu)\ot_{K}V^{\ot \ell}.
\end{split}
\end{equation*}
Then $V^{\ot \beta}_{\mathbb{K}}$ is also a right $\mathbb{K}_\beta\ot_{\bf k}{\bf k}[\mf{S}_\ell]$-module, where $1\ot r_m$ ($1\le m\le \ell-1$) acts on $V^{\ot \beta}_{\mathbb{K}}$ by $\mc{R}^{\rm norm}_{1,1}$ \eqref{eq:R matrix on W_l tensor W_m} on 
the $(m,m+1)$-component in $V^{\ot \ell}$, that is,
\begin{equation*}
\xymatrixcolsep{3pc}\xymatrixrowsep{0pc}\xymatrix{
\mathbb{K}_\nu e(\nu)\ot_K V^{\ot \ell}  \ \ar@{->}[r] & \ \mathbb{K}_{s_m(\nu)} e(s_m(\nu))\ot_K V^{\ot \ell} \\
f e(\nu)\ot (v_1\ot \dots\ot v_\ell) \ar@{|->}[r] & s_m(f) e(s_m(\nu))\ot (\dots\ot \mc{R}^{\rm norm}_{1,1}(v_m\ot v_{m+1})\ot \dots)
},
\end{equation*}
for $\nu\in J^\beta$, $f\in \mathbb{K}_\nu$ and $v_1\ot\dots\ot v_\ell\in V^{\ot \ell}$. Hence $V^{\ot \beta}_{\mathbb{K}}$ is a right $R^J(\beta)$-module by \eqref{eq:embedding of QHA}. 
Therefore, $V^{\ot \beta}_{\mathbb{K}}$ is a $(\mathcal{U}(\e),R^J(\beta))$-bimodule since $V^{\ot \beta}_{\mathbb{K}}$ is a left $\mathcal{U}(\e)$-module and the action of $\mathcal{U}(\e)$ commutes with that of $R^J(\beta)$.

\begin{prop}\label{prop:(U,QH)-bimodule}
For $\beta\in\mathbb{N}[J]$, $V^{\ot \beta}_{\mathbb{O}}$ is invariant under the action of $R^J(\beta)$. Hence it is a $(\mathcal{U}(\e),R^J(\beta))$-bimodule.
\end{prop}
\pf By Theorem \ref{thm:spectral decomposition}, we have
\begin{equation*}
\mc{R}^{\rm norm}_{1,1}(z)= \mc{P}^{1,1}_1 + \frac{1-q^2z}{z-q^2}\mc{P}^{1,1}_0,
\end{equation*}
which is the same as in the case of $\e=\e_{0|n}$.
This enables us to apply literally the same argument as in \cite[Theorem 3.3]{KKK} by replacing $q$ in \cite{KKK} with $-q^{-1}$ (cf.~Remark~\ref{rem:spectral decomp is the same}).
\qed\vskip 2mm

For a graded $R^J(\beta)$-module $M$ (not necessarily finite-dimensional), we define
\begin{equation*}
\mc{F}_{\e,\beta}(M)= V^{\otimes\beta}_{\mathbb{O}}\ot_{R^J(\beta)}M.
\end{equation*}
Then $\mc{F}_{\e,\beta}$ is a functor from the category of graded $R^J(\beta)$-modules to that of $\mathcal{U}(\e)$-modules. We also let 
\begin{equation}\label{eq:functor F}
\mc{F}_{\e,\ell}=\bigoplus_{\beta\in \mathbb{N}[J], {\rm ht}(\beta)=\ell}\mc{F}_{\e,\beta}, \quad 
\mc{F}_\e=\bigoplus_{\ell\ge 0}\mc{F}_{\e,\ell},
\end{equation}
which is a functor from the category of graded modules over $R^J(\ell)$ (resp. $R^J$) to that of $\mathcal{U}(\e)$-modules.

\begin{thm}[cf. {\cite[Theorem 3.4]{KKK}}]\label{thm:SWD functor F}
The functors in \eqref{eq:functor F} induce exact functors
\begin{equation*}
\xymatrixcolsep{2pc}\xymatrixrowsep{3pc}\xymatrix{
\mc{F}_{\e,\ell} : R^{J}(\ell)\text{\rm -gmod}  \ar@{->}[r] & \ \mc{C}^\ell(\e)},\quad
\xymatrixcolsep{2pc}\xymatrixrowsep{3pc}\xymatrix{
\mc{F}_{\e} : R^{J}\text{\rm -gmod}  \ar@{->}[r] & \ \mc{C}(\e)}.
\end{equation*}
respectively, and $\mathcal{F}_{\epsilon}$ is monoidal.
\end{thm}
\pf Thanks to Theorem \ref{thm:spectral decomposition}, we may apply the same arguments in the proof of \cite[Theorem 3.4]{KKK} as in Proposition \ref{prop:(U,QH)-bimodule}.
\qed\vskip 2mm

Next, let us describe the image of irreducible $R^J$-modules under $\mc{F}_\e$.
To compute $\mathcal{F}_\e(L(a,b))$ first, we need the following lemma
generalizing \cite[Lemma B.1]{AK} of quantum affine algebras. 

\begin{lem}\label{lem:AKlemma} 
For $\ell \geq 4$, we have the following exact sequence of $\mathcal{U}(\e)$-modules:
$$
\begin{tikzcd}[sep=10]
0 \arrow[r] & \mathcal{W}_{\ell,\e}(1) \arrow[r] 
& \mathcal{W}_{1,\e}(q^{1-\ell})\otimes \mathcal{W}_{\ell-1,\e}(q) \arrow[rrr, "{}^{\mathcal{R}^{\rm norm}_{1,\ell-1}(q^{-\ell})}"]
& &
& \mathcal{W}_{\ell-1,\e}(q)\otimes \mathcal{W}_{1,\e}(q^{1-\ell}) \arrow[r]
& \mathcal{W}_{\ell,\e}(1) \arrow[r] & 0.
\end{tikzcd}
$$
\end{lem}
\pf See Appendix \ref{appendix A}.
\qed

\begin{prop}[cf. {\cite[Proposition 4.9]{KKK}}]\label{prop:F of L(a,b)}
For a segment $(a,b)$ of length $\ell$, we have
\[
\mathcal{F}_\e(L(a,b))\cong\mathcal{W}_{\ell,\epsilon}(q^{-a-b}).
\]
\end{prop}
\pf
With the help of Proposition \ref{prop:QHAexact}, Lemma \ref{lem:AKlemma} and Corollary \ref{cor:hom-dim}, the proof of \cite[Proposition 4.9]{KKK} applies in the same manner except that the induction argument proceeds in $\ell$.
\qed

\begin{lem}\label{lem:F(r matrix)}
Let $(a,b)$ and $(a^{\prime},b^{\prime})$ be segments of lengths $\ell$ and $\ell'$, respectively, such that
$(a,b)\geq(a^{\prime},b^{\prime})$. Set
$c=q^{-a-b}$ and $c^{\prime}=q^{-a^{\prime}-b^{\prime}}$. Then $c/c^{\prime}$
is not a zero of the denominator $d_{\ell,\ell^{\prime}}(z)$ of $\mathcal{R}_{\ell,\ell^{\prime}}^{\mathrm{norm}}(z)$,
and the map
\[
\xymatrixcolsep{7pc}\xymatrixrowsep{0pc}\xymatrix{
\mathcal{W}_{\ell,\e}(c)\otimes\mathcal{W}_{\ell^{\prime},\e}(c^{\prime})\ar@{->}[r]^{\mathcal{F}_\e\left({\bf{r}}_{L(a,b),L(a^{\prime},b^{\prime})}\right)} & \mathcal{W}_{\ell^{\prime},\e}(c^{\prime})\otimes\mathcal{W}_{\ell,\e}(c)
}
\]
is equal to a nonzero constant multiple of $\mathcal{R}_{\ell,\ell^{\prime}}^{\mathrm{norm}}(c/c^{\prime})$
except the following case:
\begin{equation*}\label{eq:Fr_exception}
a^{\prime}<a\leq b^{\prime}<b,\quad 
\,M=1,\quad N\leq b^{\prime}-a+1,
\end{equation*}
in which case $\mathcal{F}_{\epsilon}({\bf{r}}_{L(a,b),L(a^{\prime},b^{\prime})})=0$.
\end{lem}
\pf Let ${\bf r}={\bf{r}}_{L(a,b),L(a^{\prime},b^{\prime})}$.
The assumption $(a,b)\geq(a^{\prime},b^{\prime})$ implies $a^{\prime}+b^{\prime}-a-b\leq\ell^{\prime}-\ell$.
Hence $c/c^{\prime}=q^{a^{\prime}+b^{\prime}-a-b}$ is not a zero of $d_{\ell,\ell^{\prime}}(z)$ since
\[
d_{\ell,\ell^{\prime}}(z)=\left(z-q^{\ell+\ell^{\prime}}\right)\left(z-q^{\ell+\ell^{\prime}-2}\right)\cdots\left(z-q^{\left|\ell-\ell^{\prime}\right|+2}\right).
\]
The second assertion follows once we determine when $\mathcal{F}_{\epsilon}\left(\mathrm{\mathbf{r}}_{L(a,b),L(a^{\prime},b^{\prime})}\right)$
is nonzero since  
\[
\mathrm{Hom}_{\mathcal{U}(\epsilon)}\left(\mathcal{W}_{\ell,\e}(c)\otimes\mathcal{W}_{\ell^{\prime},\e}(c^{\prime}),\mathcal{W}_{\ell^{\prime},\e}(c^{\prime})\otimes\mathcal{W}_{\ell,\e}(c)\right)=\mathbf{k}\mathcal{R}_{\ell,\ell^{\prime}}^{\mathrm{norm}}(c/c^{\prime}),
\]
by Corollary \ref{cor:hom-dim}.  

According to Proposition \ref{prop:QHAexact},
$\mathbf{r}$ is not an isomorphism if and only if $a^{\prime}<a\leq b^{\prime}<b$
or $a=b^{\prime}+1$. In such cases, we apply $\mathcal{F}_\e$ to the exact sequences in Proposition \ref{prop:QHAexact}(2), (3) to get an exact
sequence

\begin{center}
\begin{tikzcd}[row sep=small]
0 \arrow[r] & \mathcal{W}_{\ell_1,\e}(q^{-a'-b})\ot \mathcal{W}_{\ell_2,\e}(q^{-a-b'}) \arrow[r] \arrow[rd, phantom, ""{coordinate, name=Z}] & \mathcal{W}_{\ell,\e}(c)\ot \mathcal{W}_{\ell',\e}(c')
\arrow[dl, rounded corners,  "\hskip -5cm\mc{F}_\e(\mathbf{r})",
to path={ -- ([xshift=2ex]\tikztostart.east)
|- (Z) [near end]\tikztonodes
-| ([xshift=-2ex]\tikztotarget.west)
-- (\tikztotarget)}] \\
& \mathcal{W}_{\ell',\e}(c')\ot \mathcal{W}_{\ell,\e}(c) \arrow[r] & 
\mathcal{W}_{\ell_1,\e}(q^{-a' -b})\ot \mathcal{W}_{\ell_2,\e}(q^{-a-b'}) \arrow[r] & 0
\end{tikzcd}
\end{center}
where $\ell_{1}=b-a^{\prime}+1$
and $\ell_{2}=b^{\prime}-a+1$.
Now $\mathcal{F}_\e(\mathbf{r})=0$ if and only if 
\[
\dim
\left(\mathcal{W}_{\ell_{1},\e}(q^{-a^{\prime}-b})\otimes\mathcal{W}_{\ell_{2},\e}(q^{-a-b^{\prime}})\right)
=
\dim
\left(\mathcal{W}_{\ell,\e}(c)\otimes\mathcal{W}_{\ell^{\prime},\e}(c^{\prime})\right),
\]
and this equality holds if and only if the decompositions  of the modules on both sides in the above equation as a $\mathring{\mathcal{U}}(\e)$-module coincide (cf.~\eqref{eq:decomp of W_l and W_m}). This happens exactly when
$a^\prime < a\leq b^\prime <b,\, M=1 \text{ and } N\leq b^\prime -a+1$.
\qed

Finally we obtain our version of \cite[Theorem 4.11]{KKK}.
\begin{thm}\label{thm:F of simple for QHA}
Let $\left( (a_1,b_1),\dots,(a_t,b_t)\right)$ be an ordered multisegment with $\ell_k$ the length of $(a_k,b_k)$ and $\sum_{k}\ell_k=\ell$, and let $L$ be the corresponding irreducible $R^{J}(\ell)$-module in $R^J(\ell)$-${\rm gmod}$.
If $N=|\{\,i\,|\,\e_i=1\,\}|$ is greater than $\max\{\ell_1,\dots,\ell_t\}$, then 
$\mathcal{W}_\e({\mb l},{\mb c})$ is non-zero and
\begin{equation*}
\mc{F}_\e(L) \cong \mathcal{W}_\e({\mb l},{\mb c}),
\end{equation*}
for $({\mb l},{\mb c})\in \mc{P}^+(\e)$ with ${\mb l}=(\ell_1,\dots,\ell_t)$ and ${\mb c}=(q^{-a_1-b_1},\dots,q^{-a_t-b_t})$.
\end{thm}
\pf Let ${\bf r}$ be the map in \eqref{eq:maximal intertwiner for QHA}.
Recall that $L\cong\mathrm{Im}(\mathbf{r})$ by Proposition \ref{prop:classification of simples for QHA}. Since $\mc{F}$ is exact, we have 
$$\mc{F}(L)\cong\mc{F}({\rm Im}(\bf{r}))\cong\rm{Im}\mc{F}(\bf{r}).$$
By Lemma \ref{lem:F(r matrix)}, $\mathcal{F}(\mathbf{r})$ is equal to $R_{\mb{l},\e}(\mb c)$ up to multiplication by a nonzero scalar. 
Note that the condition $N>\max\{\ell_1,\dots,\ell_t\}$ excludes the exceptional case in Lemma~\ref{lem:F(r matrix)}. Indeed, whenever $a_j <a_i \leq b_j <b_i$ for some  $i<j$, we have $b_j-a_i+1<b_j - a_j +1=\ell_j<N$. This implies that $\mc{F}({\bf r})$ is non-zero and hence ${\rm Im}\mc{F}({\bf r})=\mathcal{W}_\e({\mb l},{\mb c})$ is non-zero.
\qed\vskip 2mm

\begin{thm}\label{thm:existence of SWD functor}
Let $\e, \e'\in \mc{E}$ be given such that $\e'<\e$.
Then for an $R^J$-module $M$ in $R^J${\rm -gmod}, we have an isomorphism of $\mathcal{U}(\e')$-modules
\begin{equation*}
\mf{tr}^{\e}_{\e'}\circ\mc{F}_\e(M)\cong \mc{F}_{\e'}(M),
\end{equation*}
which is natural in $M$.
\end{thm}
\pf Let $\beta\in \mathbb{N}[J]$ be given with ${\rm ht}(\beta)=\ell$.
Let $M$ be an $R^{J}(\beta)$-module.
Since the weight space decomposition of $\mc{F}_{\e,\beta}(M)$ is given by
\[
\mc{F}_{\e,\beta}(M)=\bigoplus_{\la\in P}\left(\mc{F}_{\e,\beta}(M)\right)_{\la},\quad\left(\mc{F}_{\e,\beta}(M)\right)_{\la}=\left(V^{\ot\beta}_{\mathbb{O}}\right)_{\la}\ot_{R^{J}(\beta)}M,
\]
we have
\begin{equation}\label{eq:tr of F(M)}
\mf{tr}_{\e'}^{\e}\left(\mc{F}_{\e,\beta}(M)\right)=
\left(\mf{tr}_{\e'}^{\e}\left({V}^{\ot\beta}_{\mathbb{O}}\right)\right)\ot_{R^{J}(\beta)}M.
\end{equation}
Since 
${V}^{\ot\beta}_{\mathbb{O}}
=\bigoplus_{\nu\in J^{\beta}}\mathbb{O}_{\nu}e(\nu)\ot_K V^{\ot \ell}$, 
we have as a $\mathcal{U}(\epsilon^\prime)$-module
\begin{equation}\label{eq:tr of SWD module}
\begin{aligned}
\mf{tr}_{\e'}^{\e}\left({V}^{\ot \beta}_{\mathbb{O}}\right) & 
=\bigoplus_{\nu\in J^{\beta}}
\mathbb{O}_{\nu}e(\nu)\ot_K \mf{tr}_{\e'}^{\e}\left( V^{\ot \ell}\right)\\
& \cong
\bigoplus_{\nu\in J^{\beta}}
\mathbb{O}_{\nu}e(\nu)\ot_K \left(\mf{tr}_{\e'}^{\e}\left(\mathcal{W}_{1,\e}\right)_{\rm aff}\ot \cdots\ot \mf{tr}_{\e'}^{\e}\left(\mathcal{W}_{1,\e}\right)_{\rm aff}\right)\\
& \cong\bigoplus_{\nu\in J^{\beta}}
\mathbb{O}_{\nu}e(\nu)\ot_K \left(\left(\mathcal{W}_{1,\e'}\right)_{\rm aff}\ot \cdots\ot\left(\mathcal{W}_{1,\e'}\right)_{\rm aff}\right)\\
& = \bigoplus_{\nu\in J^{\beta}}\mathbb{O}_{\nu}e(\nu)\ot_K {V'}^{\ot \ell}=(V')^{\ot\beta}_{\mathbb{O}},
\end{aligned}
\end{equation}
where $V'=\left(\mathcal{W}_{1,\e'}\right)_{\rm aff}$ and $K={\bf k}[X_{1}^{\pm1},\dots,X_{\ell}^{\pm1}]$.
Therefore it follows from \eqref{eq:tr of F(M)} and \eqref{eq:tr of SWD module} that 
$\mf{tr}^{\e}_{\e'}\circ\mc{F}_{\e,\beta}(M)\cong \mc{F}_{\e',\beta}(M)$, which is also natural in $M$.
\qed\vskip 2mm 

Hence we have the following diagram: 
\begin{equation*}\label{eq:diagram of SWDF}
\begin{tikzcd}
& \mc{C}(\epsilon) \arrow[dd, "\mathfrak{tr}^{\epsilon}_{\epsilon^\prime}"] \\
R^{J}\text{-gmod} \arrow[ru,"\mc{F}_\e"] \arrow[rd,"\mc{F}_{\e'}"'] & \\
& \mc{C}(\epsilon^{\prime})
\end{tikzcd}
\end{equation*}
such that $\mf{tr}^{\e}_{\e'}\circ\mc{F}_\e \cong \mc{F}_{\e'}$.
\vskip 2mm

\begin{cor}
Let $\e^\infty$ and $(\epsilon^{(k)})_{k\ge 1}$ be as in \eqref{eq:e infinite} and \eqref{eq:cond-epsilon-infty}, respectively.
There exists an exact monoidal functor 
\begin{equation*}\label{eq:F infinite rank}
\xymatrixcolsep{2pc}\xymatrixrowsep{3pc}\xymatrix{
\mc{F}_{\e^\infty}=\varprojlim \mc{F}_{\e^{(k)}} : R^{J}\text{\rm -gmod}  \ar@{->}[r] & \ \mc{C}(\epsilon^\infty)},
\end{equation*}
such that $\mf{tr}_{k}\circ \mc{F}_{\e^\infty} \cong \mc{F}_{\e^{(k)}}$ for all $k\geq 1$.
\end{cor}
\pf It follows from Lemma \ref{lem:universal mapping property of inverse limit} and Theorem \ref{thm:existence of SWD functor}.
\qed 

\begin{rem}{\rm
Note that a generalized quantum affine Schur-Weyl duality functor for usual quantum affine algebras is given in a more general setting \cite[Section 3.1]{KKK}.
We remark that this also applies to the case of $\mathcal{U}(\e)$-modules for arbitrary $\e$.

In other words, as in \cite[Section 3.1]{KKK} we can define functors to subcategories of $\mathcal{C}(\epsilon)$ associated to various choices of the triple $(J,X:J\rightarrow \mathbf{k}^{\times}, S:J\rightarrow \mathcal{S})$, where $\mc{S}$ is a set of simple objects in $\mc{C}(\e)$.
For instance,  given any orientation ${\mc Q}$ of the Dynkin diagram of $A_{m-1}$ and a height function $\xi$ on $\mc Q$, one may use the following data:
\begin{equation*}
\begin{split}
& J\subset \{1,\dots,m-1\}\times \mathbb{Z},\\
& X(i,p)=q^{-p-m},\quad S(i,p)=\mathcal{W}_{i,\e}(1) \quad \text{ for } (i,p)\in J,
\end{split}
\end{equation*}
as given in \cite[Section 3]{KKK2} to obtain a functor from $R^J$-gmod to $\mathcal{C}_{\mc Q}(\e)$. Here $\mathcal{C}_{\mc Q}(\e)$ is a category of $\mathcal{U}(\e)$-modules, which can be viewed as an analog of the category $\mathcal{C}_{\mc Q}$ for the quantum affine algebra of type $A_{m-1}^{(1)}$ introduced by Hernandez-Leclerc \cite{HL10}. We should remark that $m$ is not necessarily equal to $n$, the length of $\e$, whenever $\mc{S}$ is well-defined. The functor is defined in a natural way for $\epsilon$ and compatible with the truncation due to Theorem \ref{thm:spectral decomposition} and Lemma \ref{lem:truncation of Verma and R matrix}. The type of the corresponding quiver Hecke algebra $R^J$ is of type $\mc{Q}^{\rm rev}$, the quiver obtained by reversing all the arrows of $\mc{Q}$, as in \cite[Theorem 4.3.1]{KKK2} since the denominators of normalized $R$ matrices \eqref{eq:denominator} associated to simple objects in $\mc{S}$ coincide with those as in the case of $U_q'\left(A_{m-1}^{(1)}\right)$ (cf.~Remark \ref{rem:spectral decomp is the same}).

We also have a functor from the category of not necessarily graded finite-dimensional modules over quiver Hecke algebras. This is directly related to the more classical quantum affine Schur-Weyl duality \cite{CP}, which is dealt with in the next section.
}
\end{rem}

\section{Duality}\label{sec:duality}
\subsection{Subcategory $\mc{C}_J(\e)$}
Let us introduce a subcategory of $\mc{C}(\e)$ to describe the image of the generalized quantum affine Schur-Weyl duality functor $\mc{F}_\e$.
\begin{df}{\rm \mbox{}
\begin{itemize}
\item[(1)] Let $\mc{C}_{J}(\epsilon)$ be the full subcategory
of $\mc{C}(\epsilon)$ whose objects have composition factors
appearing as a composition factor of a tensor product of modules of
the form $\mathcal{W}_{1,\e}(q^{-2j})$ for $j\in J$. 
For $\ell\in\mathbb{Z}_+$, let $\mc{C}^\ell_{J}(\e)=\mc{C}_{J}(\e)\cap\mc{C}^\ell(\e)$.

\item[(2)] Let $\mc{P}^+_{J}(\e)$ be a subset of $\mc{P}^+(\epsilon)$
consisting of elements $(\mb{l},\mb{c})$ with ${\mb l}=(l_1,\dots,l_t)\in \mathbb{Z}_+^t$, and $\mb{c}=(c_1,\dots,c_t)\in (\Bbbk^\times)^t$ for some $t\geq 1$ such that 
\begin{equation}\label{eq:cond for P_J^+}
c_{i}\in q^{l_{i}-1+2\mathbb{Z}}\quad (i=1,\dots ,t).
\end{equation}

\end{itemize}
}
\end{df}

\begin{lem}\label{lem:SWD functo on J}
For an $R^J$-module $M$ in $R^J$-{\rm gmod}, we have $\mc{F}_\e(M)\in \mc{C}_J(\e)$.
Hence we have a functor
\begin{equation*}
\xymatrixcolsep{2pc}\xymatrixrowsep{3pc}\xymatrix{
\mc{F}_{\e} : R^{J}\text{\rm -gmod}  \ar@{->}[r] & \ \mc{C}_J(\e)}.
\end{equation*}
\end{lem}
\pf Recall that the convolution product $M_1\circ M_2$ is exact both in $M_1$ and $M_2$ by \cite[Proposition 2.16]{KL}, and  $L(a,b)$ is a composition factor of a convolution of product of $L(i)$'s for $i\in J$ by Proposition \ref{prop:QHAexact}(3).
By Proposition \ref{prop:F of L(a,b)}, we have $\mc{F}_\e(L(a,b))=\mathcal{W}_{l,\e}(q^{-a-b})\in \mc{C}_J(\e)$, and hence $\mathcal{W}_{\e}(l,c)\in \mc{C}_J(\e)$ for $(l,c)\in \mc{P}_J^+(\e)$.

Let $0\subset M_1\subset\dots\subset M_r=M$ be a composition series of $M$. Since $\mc{F}_\e$ is exact, $0\subseteq \mc{F}_\e(M_1)\subseteq\dots\subseteq \mc{F}_\e(M_r)=\mc{F}_\e(M)$ is a filtration of $\mc{F}_\e(M)$ whose subquotients are either zero or of the form $\mathcal{W}_{\epsilon}(\mb{l},\mb{c})$ for $({\mb l},{\mb c})\in \mc{P}^+_J(\e)$ by Proposition~\ref{prop:classification of simples for QHA} and Theorem~\ref{thm:F of simple for QHA}. Now note that $\mathcal{W}_{\epsilon}(\mb{l},\mb{c})$ is a quotient of $\mathcal{W}_{{\mb l},\e}({\mb c})$, and hence a composition factor of a tensor product of $\mathcal{W}_{1,\e}(q^{-2j})$'s.
\qed\vskip 2mm

By Theorem \ref{thm:SWD functor F} and Lemma \ref{lem:SWD functo on J}, we have $\mc{F}_{\e}=\bigoplus_{\ell\in\mathbb{Z}_+}\mc{F}_{\e,\ell}$, where
\begin{equation}\label{eq:SWD functo on J ell}
\xymatrixcolsep{2pc}\xymatrixrowsep{3pc}\xymatrix{
\mc{F}_{\e,\ell} : R^{J}(\ell)\text{\rm -gmod}  \ar@{->}[r] & \ \mc{C}^\ell_J(\e)}.
\end{equation}

\begin{prop}\label{prop:classification-irrep-C_Z}
Let $V$ be a $\mathcal{U}(\e)$-module in $\mc{C}_{J}(\e)$.
Then $V$ is irreducible if and only if  
such that $V\cong\mathcal{W}_{\e}(\mb{l},\mb{c})$ for some $(\mb{l},\mb{c})\in\mc{P}^+_{J}(\e)$.
\end{prop}
\pf Suppose that $V$ is a composition factor of $\mathcal{W}_{1,\e}(q^{-2j_1})\ot\dots\ot\mathcal{W}_{1,\e}(q^{-2j_t})$ for some $j_1,\dots,j_t\in J$. 
Then $V$ is isomorphic to $\mc{F}_\e(M)$ for some composition factor $M$ of $L(j_1)\circ\dots\circ L(j_t)$. By Proposition \ref{prop:classification of simples for QHA} and Theorem \ref{thm:F of simple for QHA}, $V\cong \mathcal{W}_{\e}(\mb{l},\mb{c})$ for some $(\mb{l},{\mb c})\in \mc{P}^+_J(\e)$. Conversely, it is clear by definition that $\mathcal{W}_\e(\mb{l},\mb{c})$ belongs to $\mc{C}_J(\e)$.
\qed

\vskip 2mm 

The category $\mc{C}_J(\e)$ is compatible with the truncation in the following sense.

\begin{lem}\label{lem:C_J under truncation}
Let $\e, \e'\in \mc{E}$ be given such that $\e'<\e$.   
For a $\mathcal{U}(\e)$-module $V$ in $\mc{C}_J(\e)$, let $V'=\mf{tr}^{\e}_{\e'}(V)$. 
Then we have
\begin{itemize}
\item[(1)] $V'$ is a $\mathcal{U}(\e')$-module in $\mc{C}_J(\e')$,

\item[(2)] $\ell(V')\le \ell(V)$, where $\ell(V)$ denotes the length of the composition series of $V$.
\end{itemize}
\end{lem}
\pf Since $\mf{tr}^{\e}_{\e'}$ is exact, it follows from Theorem~\ref{thm:truncation sends simple to simple} and Proposition~\ref{prop:classification-irrep-C_Z}.
\qed

\subsection{Affine Hecke algebras and quantum affine Schur-Weyl duality functor}

There is a quantum affine Schur-Weyl duality between the affine Hecke algebra and quantum affine algebra of type $A$ \cite{CP}. In this subsection, we prove an analogue for $\mathcal{U}(\e)$ and explain how to identify it with our functor $\mc{F}_{\e,\ell}$.

First, we briefly recall that a quiver Hecke algebra of type $A$ is isomorphic to an affine Hecke algebra (of type $A$) after suitable completions on both sides \cite{BK, R}. Our exposition is based on \cite{Kim}.

\begin{defn}
Let $\ell\geq 2$. The {\em affine Hecke algebra} $\aha$ (of the symmetric group $\mathfrak{S}_{\ell}$)
is the associative $\mathbf{k}$-algebra with 1 generated by $X_{k}^{\pm1}$ ($1\leq k\leq\ell$) and $h_{m}$ ($1\leq m\leq\ell-1$) subject to the
following relations:
\begin{gather*}
 h_{m}h_{m+1}h_{m}=h_{m+1}h_{m}h_{m+1},\quad  h_{m}h_{m'}=h_{m'}h_{m}\quad (\left|m-m'\right|>1),\\
(h_{m}-q^{2})(h_{m}+1)=0,\\
 X_{k}X_{k'}=X_{k'}X_{k},\quad X_{k}X_{k}^{-1}=X_{k}^{-1}X_{k}=1,\\
 h_{m}X_{m}h_{m}=q^{2}X_{m+1}\quad  h_{m}X_{k}=X_{k}h_{m}\quad (k\neq m,m+1).
\end{gather*}
For convenience, we assume $H^{\mathrm{aff}}_{0}(q^2)=\mathbf{k}$ and $H^{\mathrm{aff}}_{1}(q^2)=\mathbf{k}[X^{\pm1}]$.

The subalgebra $\fha$ generated by $h_{m}$ ($1\leq m\leq\ell-1$)
is called the {\em finite Hecke algebra}. 
\end{defn}

Note that given a reduced expression $w=s_{i_1}\cdots s_{i_n}\in \mf{S}_{\ell}$, the element $h_w\coloneqq h_{i_1}\cdots h_{i_n}$ does not depend on the reduced expression, and the set $\left\{ h_w \right\}_{w\in \mathfrak{S}_\ell}$ is a basis for $\fha$.

Put
$$
\mathbb{O}_{\ell}=\bigoplus_{\beta\in J^\ell} \mathbb{O}_{\beta},\quad 
\mathbb{K}_{\ell}=\bigoplus_{\beta\in J^\ell} \mathbb{K}_{\beta} 
$$
(cf.~\eqref{eq:O and K}), and
consider the following completions:
\begin{align*}
\mathbb{O}\aha & :=\mathbb{O}_{\ell}\otimes_{\mathbf{k}\left[X_{1}^{\pm1},\dots,X_{\ell}^{\pm1}\right]}\aha\cong\mathbb{O}_{\ell}\otimes_{\mathbf{k}}\fha,\\
\mathbb{K}\aha & :=\mathbb{K}_{\ell}\otimes_{\mathbb{O}_{\ell}}\mathbb{O}\aha\cong\mathbb{K}_{\ell}\otimes_{\mathbf{k}}\fha,
\end{align*}
as ${\bf k}$-vector spaces. We simply write $fe(\nu)h=fe(\nu)\ot h$ for $fe(\nu)\in \mathbb{K}_\ell$ and $h\in \fha$.
We regard $\mathbb{O}H_{\ell}^{\mathrm{aff}}(q^{2})$ and $\mathbb{K}H_{\ell}^{\mathrm{aff}}(q^{2})$ as associative ${\bf k}$-algebras, where

\begin{equation*}\label{eq:mult on aha}
h_{m}fe(\nu)  =
s_{m}(f)e(s_{m}(\nu))h_{m}+\dfrac{(q^{2}-1)X_{m+1}\left(fe(\nu)-s_{m}(f)e(s_{m}(\nu))\right)}{X_{m+1}-X_{m}}
\end{equation*}
for $1\le m\le \ell-1$ and $fe(\nu)\in\mathbb{O}_{\ell}\text{ or }\mathbb{K}_{\ell}$. 
Note that we have a sequence of subalgebras
\begin{equation*}
\aha\subset\mathbb{O}\aha\subset\mathbb{K}\aha.
\end{equation*}

For $1\le m\le \ell-1$, define the intertwiner $\Phi_m\in \mathbb{K}H_{\ell}^{\mathrm{aff}}(q^{2})$ by
\begin{equation*}
\Phi_{m} =h_{m}-\dfrac{(q^{2}-1)X_{m+1}}{X_{m+1}-X_{m}}.
\end{equation*}
The following properties of $\Phi_{m}$ can be checked in a straightforward way:
\begin{equation}\label{eq:property of Phi}
\begin{gathered}
\Phi_{m}f(X_{m},X_{m+1})e(\nu)=f(X_{m+1},X_{m})e(s_{m}(\nu))\Phi_{m},\\
\Phi_{m}\Phi_{m'}=\Phi_{m'}\Phi_{m} \quad (\left|m-m'\right|>1),\\
\Phi_{m}\Phi_{m+1}\Phi_{m}=\Phi_{m+1}\Phi_{m}\Phi_{m+1}
\end{gathered}
\end{equation}
for $f(X_{m},X_{m+1})\in\mathbf{k}(X_{m},X_{m+1})$ and $\nu\in J^\ell$.
Moreover, we have
\[
\Phi_{m}^{2}e(\nu)=
\frac{X_{m+1}-q^{2}X_{m}}{X_{m}-X_{m+1}}\cdot\frac{X_{m}-q^{2}X_{m+1}}{X_{m+1}-X_{m}}e(\nu).
\]
Hence if we normalize $\Phi_m$ by
\[
\widetilde{\Phi}_{m}=
\frac{X_{m}-X_{m+1}}{X_{m+1}-q^{2}X_{m}}\Phi_{m},
\]
then we have $\widetilde{\Phi}_{m}^{2}=1$.

\begin{thm}[\cite{BK, R}]\label{thm:Psi}
Let $\widehat{R}^{J}(\ell)  =\mathbf{k}\left\llbracket x_{1},\dots,x_{\ell}\right\rrbracket \otimes_{\mathbf{k}\left[x_{1},\dots,x_{\ell}\right]}R^{J}(\ell)$ be the completion of $R^J(\ell)$ with multiplication naturally extended from $R^{J}(\ell)$.
Then there is an isomorphism of $\bf k$-algebras
\begin{equation*}\label{eq:Psi}
\xymatrixcolsep{2pc}\xymatrixrowsep{3pc}\xymatrix{
\Psi:\widehat{R}^{J}(\ell)  \ar@{->}[r] & \ \mathbb{O}\aha,}
\end{equation*}
given by 
\begin{align*}
\Psi(e(\nu)) & =e(\nu),\\
\Psi(e(\nu)x_{k}) & =e(\nu)X(\nu_{k})^{-1}(X_{k}-X(\nu_{k})),\\
\Psi(e(\nu)\tau_{m}) & =
\begin{cases}
e(\nu)(\widetilde{\Phi}_{m}-1) \left(\frac{X_{m}}{X(\nu_{m})}-\frac{X_{m+1}}{X(\nu_{m+1})}\right)^{-1} & \text{if }\nu_{m+1}=\nu_{m},\\
e(\nu)\widetilde{\Phi}_{m} \left(\frac{X_{m}}{X(\nu_{m})}-\frac{X_{m+1}}{X(\nu_{m+1})}\right) & \text{if }\nu_{m+1}=\nu_{m}+1,\\
e(\nu)\widetilde{\Phi}_{m} & \text{otherwise},
\end{cases}
\end{align*}
for $\nu\in J^\ell$, $1\le k\le \ell$ and $1\le m\le \ell-1$.
\end{thm}

\pf Let us give a proof briefly for the reader's convenience. 
First we see that all $\Psi(e(\nu)\tau_{m})$ indeed belong to $\mathbb{O}\aha$.
By \eqref{eq:property of Phi}, the embedding \eqref{eq:poly-repn-QHA}
implies the existence of $\Psi$ since we have in \eqref{eq:poly-repn-QHA}
\[
e(\nu)P_{\nu_m ,\nu_{m+1}}(x_m , x_{m+1})=
\begin{cases}
e(\nu)(x_m - x_{m+1}) & \text{if } \nu_{m+1}=\nu_m +1,\\
e(\nu) & \text{otherwise},
\end{cases}
\] and $e(\nu)x_k = e(\nu)(X_k / X(\nu_k) -1)$.

Recall that 
$\left\{\tau_w \right\}_{w\in \mathfrak{S}_\ell }$ is an $\mathbb{O}_\ell$-basis for $\widehat{R}^{J}(\ell)$. 
Since $\Psi$ is $\mathbb{O}_{\ell}$-linear by definition, 
it is enough to show that $\left\{ \Psi(\tau_{w})\right\} _{w\in \mathfrak{S}_{\ell}}$ is an $\mathbb{O}_\ell$-basis for $\mathbb{O}\aha$, which implies that $\Psi$ is bijective.

Note that $\left\{ h_w \right\}_{w\in\mathfrak{S}_\ell}$ is also an $\mathbb{O}_\ell$-basis for $\mathbb{O}\aha$. We claim that each $h_w$ is given by an $\mathbb{O}_\ell$-linear combination of $\Psi(\tau_w)$'s.
Moreover, since $\Psi$ is an algebra homomorphism, it is enough to consider $e(\nu)h_m$.

Suppose $\nu_{m} = \nu_{m+1}$. Then 
\begin{align*}
\Psi(e(\nu)\tau_{m}) & =e(\nu)(\widetilde{\Phi}_m -1)(X_m - X_{m+1})^{-1}X(\nu_m) \\
& =X(\nu_m) e(\nu) \left( (X_{m+1}- X_m )^{-1}\widetilde{\Phi}_m - (X_m -X_{m+1})^{-1} \right) \\
& =X(\nu_m) e(\nu) (X_{m+1}-X_m )^{-1} \left( \frac{X_m-X_{m+1}}{X_{m+1}-q^{2}X_m}h_m +\frac{(q^{2}-1)X_{m+1}}{X_{m+1}-q^{2}X_m} +1 \right) \\
& =X(\nu_m) e(\nu) (X_{m+1}-X_m ) \left( \frac{X_m-X_{m+1}}{X_{m+1}-q^{2}X_m}h_m + \frac{q^{2}(X_{m+1}-X_m )}{X_{m+1}-q^{2}X_m} \right) \\
& =q^{-2\nu_m} e(\nu) \left( -\frac{1}{X_{m+1}-q^{2}X_m}h_m + \frac{q^2}{X_{m+1}-q^{2}X_m} \right)
\end{align*}
so that 
\[
e(\nu) h_m = q^{2\nu_m}(q^2 X_m - X_{m+1})\Psi(e(\nu)\tau_m ) + q^2 e(\nu)
\]
as desired. The other cases can be checked similarly. This proves the claim.
\qed\vskip 2mm

Now, let us prove an analogue of quantum affine Schur-Weyl duality \cite{CP} for $\mathcal{U}(\e)$ and explain how to identify it with our functor $\mc{F}_{\e,\ell}$.

Recall that $V=(\mathcal{W}_{1,\e})_{\rm aff}$. Put 
\begin{equation*}\label{eq:cali V}
\mc{V}=\mathcal{W}_{1,\e}\subset V.
\end{equation*}
Assume that $V^{\ot \ell}=(\mc{V}_{\rm aff})^{\ot \ell}$ is a ${\bf k}[X_1^{\pm 1},\dots,X_\ell^{\pm 1}]$-module where $X_k$ acts on the $k$-th component as the automorphism $z$.
For $1\le m\le \ell-1$, let
\[
R_{m}: V^{\otimes\ell}\rightarrow\mathbf{k}(X_{m},X_{m+1})\otimes_{\mathbf{k}}V^{\otimes\ell}
\]
denote the map given by applying $\mathcal{R}_{1,1}^{\mathrm{norm}}(X_{m}/X_{m+1})$ \eqref{eq:R matrix on W_l tensor W_m}
on the $m$-th and $(m+1)$-st component of $\mathcal{V}^{\ot \ell}$. 
Then we can check the following.
\begin{lem}\label{lem:afa action on tensor power}
We have the following.
\begin{itemize}
\item[(1)]
$V^{\otimes\ell}$ has a right $\aha$-module structure
given by
\begin{align*}
\left(f\otimes v\right)X_{k} & =\left(X_{k}f\right)\otimes v,\\
\left(f\otimes v\right)h_{m} & =
\left(f\otimes v\right)\left(\frac{X_{m+1}-q^{2}X_{m}}{X_{m}-X_{m+1}}R_{m}+\frac{(q^{2}-1)X_{m+1}}{X_{m+1}-X_{m}}\right)
\end{align*}
for $f\otimes v\in V^{\otimes\ell}$, $1\le k\le \ell$, and $1\le m\le \ell-1$. 

\item[(2)] $\mathcal{V}^{\otimes\ell}\subset V^{\otimes\ell}$ is invariant under the action of $\fha\subset \aha$, and hence $V^{\otimes\ell}\cong\mathcal{V}^{\otimes\ell}\otimes_{\fha}\aha$ as an $\aha$-module.
\end{itemize}
\end{lem}

So we may regard $V^{\otimes\ell}_{\mathbb{O}}\coloneqq \mathbb{O}_{\ell}\otimes_{\mathbf{k}[X^{\pm1}_1,\dots,X^{\pm1}_{\ell}]}V^{\otimes\ell}$ as a right $\mathbb{O}\aha$-module via the following isomorphisms of ${\bf k}$-vector space,
\begin{equation}\label{eq:QHA AHA module structure}
\begin{split}
V^{\otimes\ell}_{\mathbb{O}}
&\cong \mathcal{V}^{\otimes\ell}\otimes_{\mathbf{k}}\mathbb{O}_{\ell}\cong\left(\mathcal{V}^{\otimes\ell}\otimes_{\fha}\fha\right)\otimes_{\mathbf{k}}\mathbb{O}_{\ell}\\
&\cong\mathcal{V}^{\otimes\ell}\otimes_{\fha}\left(\fha\otimes_{\mathbf{k}}\mathbb{O}_{\ell}\right)\cong\mathcal{V}^{\otimes\ell}\otimes_{\fha}\left(\mathbb{O}_{\ell}\otimes_{\mathbf{k}}\fha\right)\\
&\cong\mathcal{V}^{\otimes\ell}\otimes_{\fha}\mathbb{O}\aha.
\end{split}
\end{equation}
Hence $V^{\otimes\ell}_{\mathbb{O}}$ becomes a $\left(\mathcal{U}(\e),\mathbb{O}\aha\right)$-bimodule. 
We remark that the action of $\widetilde{\Phi}_{m}$ on $V^{\ot\ell}_{\mathbb{O}}$ coincides with that of $R_{m}$.
Moreover we can check the following.
\begin{lem}
The $\left(\mathcal{U}(\e),\widehat{R}^{J}(\ell)\right)$-bimodule structure on $V^{\otimes\ell}_{\mathbb{O}}$ induced from $\Psi$ in Theorem \ref{thm:Psi} coincides with the one defined in Section \ref{subsec:QASWD}.
\end{lem}

Let $\aha\text{-mod}_{J}$ be the category of finite-dimensional
$\aha$-modules such that eigenvalues of $X_{k}$ ($1\le k\leq\ell$)
are in the set $\left\{\,X(j)=q^{-2j}\,|\,j\in J\,\right\}$. 
Then 
\begin{equation*}
\aha\text{-mod}_{J} \cong \mathbb{O}\aha\text{-mod},
\end{equation*}
where $\mathbb{O}\aha\text{-mod}$ is the category of finite-dimensional
$\mathbb{O}\aha$-modules $(e(\nu)$ acts as the projection to the generalized
eigenspace), and $\cong$ means an equivalence of categories. 

On the other hand, let $\widehat{R}^{J}(\ell)\text{-mod}$ be the category of (not necessarily graded) finite-dimensional $\widehat{R}^{J}(\ell)$-modules. Note that 
\begin{equation*}
\widehat{R}^{J}(\ell)\text{-mod} \cong R^{J}(\ell)\text{-mod}_0,
\end{equation*}
where $R^{J}(\ell)\text{-mod}_0$ is the category of finite-dimensional $R^J(\ell)$-modules such that $x_{k}$ acts nilpotently.

Since we have an equivalence
\begin{equation*}
\mathbb{O}\aha\text{{\rm -mod}}\cong\widehat{R}^{J}(\ell)\text{{\rm -mod}},
\end{equation*}
which is induced from the isomorphism $\Psi$ in Theorem \ref{thm:Psi}, we have an equivalence
\begin{equation*}
\Psi^{*}: \aha\text{{\rm -mod}}_{J} \overset{\cong}{\longrightarrow} R^{J}(\ell)\text{{\rm -mod}}_0.
\end{equation*}
Summarizing the above arguments, we have the following.
\begin{prop}
We have a functor
\begin{equation*}
\xymatrixcolsep{2pc}\xymatrixrowsep{0pc}\xymatrix{
\mc{F}_{\e,\ell}^{\ast} : \aha\text{{\rm -mod}}_{J} \ar@{->}[r] & \ \mc{C}_J^\ell(\e) \\
M \ar@{|->}[r] & \ V_{\mathbb{O}}^{\otimes\ell}\ot_{\mathbb{O}\aha}M.}
\end{equation*}
We also have 
$$\mathcal{F}_{\e,\ell}^{\ast}(M)\cong\mathcal{F}_{\e,\ell}\circ\Psi^*(M),$$ 
where $\mc{F}_{\e,\ell}$ is given in \eqref{eq:functor F} and the isomorphism is natural in $M$.
\end{prop}

We see from \eqref{eq:QHA AHA module structure}
\begin{equation*}
\mc{F}_{\e,\ell}^{\ast}(M)=V_{\mathbb{O}}^{\otimes\ell}\ot_{\mathbb{O}\aha}M\cong\mathcal{V}^{\otimes\ell}\otimes_{\fha}M
\end{equation*}
for a left $\aha$-module $M$, which is related to the Schur-Weyl duality of finite type as follows (note that \cite{KY} deals with $H_{\ell}(q^{-2})$ where $h^{-1}_i$ is used as a generator instead of $h_i$).

\begin{prop}[{\cite[Theorem 3.1]{KY}}]\label{thm:finite-SWD} 
Let $\mathcal{R}:\mathcal{V}^{\otimes2}\longrightarrow\mathcal{V}^{\otimes2}$
be a $\mathring{\mathcal{U}}(\epsilon)$-linear map given by
\[
\mathcal{R}\left(\ket{\mathbf{e}_{i}}\otimes\ket{\mathbf{e}_{j}}\right)=\begin{cases}
qq_{i}\ket{\mathbf{e}_{i}}\otimes\ket{\mathbf{e}_{i}} & \text{if }i=j,\\
q\ket{\mathbf{e}_{j}}\otimes\ket{\mathbf{e}_{i}} & \text{if }i>j,\\
(q^{2}-1)\ket{\mathbf{e}_{i}}\otimes\ket{\mathbf{e}_{j}}+q\ket{\mathbf{e}_{j}}\otimes\ket{\mathbf{e}_{i}} & \text{if }i<j.
\end{cases}
\]
Then $\mathcal{V}^{\otimes\ell}$ has a $\left(\mathring{\mathcal{U}}(\epsilon),\fha\right)$-bimodule
structure given by $v.h_{m}=\mathcal{R}_{m}(v)$ for $1\le m\le \ell-1$, where $\mathcal{R}_{m}$
denotes the map given by applying $\mathcal{R}$ on the $m$-th and $(m+1)$-st components of $\mc{V}^{\ot \ell}$ and the identity elsewhere.
Hence we have a functor 
\[
\xymatrixcolsep{2pc}\xymatrixrowsep{0pc}\xymatrix{
\mathcal{J}_{\ell} : \fha\text{{\rm -mod}} \ar@{->}[r] & \ \ring{\mc{C}}^\ell(\e)  \\
M \ar@{|->}[r] & \ \mathcal{V}^{\otimes\ell}\otimes_{\fha}M},
\]
which is an equivalence for $\ell \leq n$.
\end{prop}

\begin{rem}{\rm
One can directly check that the right $\fha$-action on $\mathcal{V}^{\otimes\ell}$
in Proposition \ref{thm:finite-SWD} coincides with the one in Lemma \ref{lem:afa action on tensor power}(2).
Hence, as $\mathring{\mathcal{U}}(\epsilon)$-modules, 
\begin{equation*}\label{eq:two functors}
\mathcal{F}_{\e,\ell}^{\ast}(M)\cong\mathcal{J}_{\ell}(M).
\end{equation*}
}
\end{rem}

The following is an analogue of the result \cite[Section 4]{CP} for $\mathcal{U}(\e)$. 
\begin{thm}\label{thm:equiv-QASWD}
For $\ell < n$, the functor $\mc{F}_{\e,\ell}^{\ast}$ is an equivalence of categories.
\end{thm}
\pf
The proof is given in Appendix \ref{appendix B}.
\qed

\begin{cor}\label{cor:equiv for large rank}
For $\ell < n$, we have an equivalence of categories 
\begin{equation*}
\xymatrixcolsep{2pc}\xymatrixrowsep{3pc}\xymatrix{
\mc{F}_{\e,\ell} : R^{J}(\ell)\text{\rm -mod}_0  \ar@{->}[r] & \ \mc{C}^\ell_J(\e)}.
\end{equation*}
\end{cor}

\begin{cor}\label{cor:truncation is an equivalence}
Suppose that $\e=(\e_1,\dots,\e_n), \e'=(\e'_1,\dots,\e'_{n'})\in \mc{E}$ are given such that $\e'<\e$ and $\ell <  n'<n$.
Then $\mf{tr}^\e_{\e'}$ induces an equivalence of categories
\begin{equation*}
\xymatrixcolsep{2pc}\xymatrixrowsep{3pc}\xymatrix{
\mf{tr}^\e_{\e'} : \mc{C}_J^\ell(\e)  \ar@{->}[r] & \ \mc{C}_J^\ell(\e')}.
\end{equation*}
\end{cor}
\pf By Theorem \ref{thm:existence of SWD functor}, we have 
$\mf{tr}^{\e}_{\e'}\circ\mc{F}_{\e,\ell} \cong \mc{F}_{\e',\ell}$. Since $\mc{F}_{\e,\ell}$ and $\mc{F}_{\e',\ell}$ are equivalences by Corollary \ref{cor:equiv for large rank}, $\mf{tr}^{\e}_{\e'}$ is an equivalence.
\qed

\subsection{Equivalence between $R^{J}\text{\rm -mod}_0$ and $\mc{C}_J(\e^\infty)$}
In this subsection, we reformulate the generalized quantum affine Schur-Weyl duality functor as an equivalence between $R^{J}\text{\rm -mod}_0$ and a certain subcategory of $\mc{C}(\e^\infty)$ for any given $\e^\infty$.

\begin{df}{\rm 
Let $\e^\infty$ and $(\epsilon^{(k)})_{k\ge 1}$ be as in \eqref{eq:e infinite} and \eqref{eq:cond-epsilon-infty}, respectively.

\begin{itemize}
\item[(1)] 
Let $\mc{C}_J(\e^\infty)$ be the full subcategory
of $\mc{C}(\e^\infty)$ with objects $\mathbb{V}$ such that 

\noindent (i)  $\mf{tr}_k(\mathbb{V})\in \mc{C}_J(\e^{(k)})$ for $k\geq 1$, 

\noindent (ii) $\ell(\mf{tr}_k(\mathbb{V}))=\ell(\mf{tr}_{k+1}(\mathbb{V}))$ for all sufficiently large $k$.

\item[(2)]
For $\ell\in \mathbb{Z}_+$, let $\mc{C}^\ell_J(\e^\infty)$ be the full subcategory of $\mc{C}_J(\e^\infty)$ with objects $\mathbb{V}$ such that $\mf{tr}_k(\mathbb{V})\in \mc{C}^\ell_J(\e^{(k)})$ for $k\geq 1$.

\item[(3)]  
Let $\mc{P}^+_J$ be the set of $({\mb l},{\mb c})\in \mc{P}^+$ satisfying \eqref{eq:cond for P_J^+}. 

\end{itemize}
}
\end{df}
Note that $\mc{C}_J(\e^\infty)$ is monoidal.
We remark that $\mc{C}_J(\e^\infty)$ is not the inverse limit category associated to $\{\,\mc{C}_J(\e^{(k)})\,\}_{k \ge 1}$. 
Instead, we have the following.
\begin{lem}\label{lem:subcategory C J ell}
\mbox{}
\begin{itemize}
\item[(1)]
We have the following decomposition
\begin{equation*}\label{eq:decomp of C_J}
\mc{C}_J(\e^\infty)=\bigoplus_{\ell\in\mathbb{Z}_+}\mc{C}^\ell_J(\e^\infty).
\end{equation*}

\item[(2)]
For $\ell\in\mathbb{Z}_+$, we have  
$$\mc{C}^\ell_J(\e^\infty)=\varprojlim \mc{C}^\ell_J(\e^{(k)}),$$
which is the inverse limit category associated to $\{\,\mc{C}^\ell_J(\e^{(k)})\,\}_{k\geq 1}$ defined as in Definition \ref{def:inverse limit category}.
\end{itemize}
\end{lem}
\pf (1) Suppose that $\mathbb{V}$ is an object in $\mc{C}_J(\e^\infty)$. 
If we put $V_k=\mf{tr}_k(\mathbb{V})$ for $k\geq 1$, then we have $V_k=\bigoplus_{\ell}V_k^\ell$ where $V_k^\ell\in \mc{C}_J^\ell(\e^{(k)})$. Let $\mathbb{V}^\ell$ be an object in $\mc{C}_J^\ell(\e^\infty)$ such that $\mf{tr}_k(\mathbb{V}^\ell)=V_k^\ell$ for $k\geq 1$.

For $k\geq 1$, let $\ell_k=\max\{\,\ell\,|\,V_k^\ell\neq 0\,\}$. If the sequence $\{\ell_k\}_{k\geq 1}$ is unbounded, then so is $\{\ell(V_k)\}_{k\geq 1}$, which is a contradiction. 
Hence $\mathbb{V}=\bigoplus_\ell\mathbb{V}^\ell$, where $\mathbb{V}^\ell$ is zero except for a finitely many $\ell$'s.

(2) Let $\mc{A}_1$ and $\mc{A}_2$ denote the categories on the lefthand side and the righthand side of the statement respectively. It is clear that $\mc{A}_1$ is a subcategory of $\mc{A}_2$. 

Let $\mathbb{V}$ be an object in $\mc{A}_2$, and let $V_k=\mf{tr}_k(\mathbb{V})$ for $k\geq 1$. Given $k$, suppose that $\mathcal{W}_{\e^{(k)}}({\mb l},{\mb c})$ is a composition factor of $V_k$. Note that $\mathcal{W}_{\e^{(k)}}({\mb l},{\mb c})$ is a direct sum of $V_{\e^{(k)}}(\la)$'s for partitions $\la$ of $\ell$. If $k$ is sufficiently large, then we see from Proposition \ref{prop:truncation of poly and fund} that $\mf{tr}^{k}_{k-1}(V_{\e^{(k)}}(\la))$ is non-zero for all partitions $\la$ of $\ell$, and then $\mf{tr}^{k}_{k-1}(\mathcal{W}_{\e^{(k)}}({\mb l},{\mb c}))$ is non-zero. This implies that $\ell(V_k)=\ell(V_{k'})$ for $k'\ge k$, and hence $\mathbb{V}$ is an object in $\mc{A}_1$. Therefore $\mc{A}_1=\mc{A}_2$.
\qed\vskip 2mm

\begin{lem}\label{lem:simple object in restricted inverse limit}
Let $\mathbb{V}$ be a nonzero object in $\mc{C}_{J}(\e^\infty)$. Then $\mathbb{V}$ is simple if and only if  $\mf{tr}_k(\mathbb{V})$ is irreducible or zero for all $k\geq 1$.
\end{lem}
\pf By Lemma \ref{lem:simplicity}, the converse is true. Suppose that $\mathbb{V}$ is simple.
Then by \cite[Lemma 4.1.5]{EA}, any non-zero $\mf{tr}_k(\mathbb{V})$ is irreducible.
\qed

\begin{lem}\label{prop:classification-irrep-C_Z-2}
If $\mathcal{W}_{\e^\infty}(\mb{l},\mb{c})\cong \mathcal{W}_{\e^\infty}({\mb{l}^\ast},{\mb{c}}^\ast)$ for some $(\mb{l},\mb{c})$ and $({\mb{l}^\ast},{\mb{c}}^\ast)\in \mc{P}_J^+$, then $(\mb{l},\mb{c})$ and $ ({\mb{l}^\ast},{\mb{c}}^\ast)$ are equal up to permutation.
\end{lem}
\pf Take a sufficiently large $k$ such that the number of occurrence of $1$ in $\e:=\e^{(k)}$, say $N$, is large enough, and put $\e'=\e_{0|N}<\e$. Then $\mathcal{W}_{\e}(\mb{l},\mb{c})$ and $\mathcal{W}_{\e'}(\mb{l},\mb{c})$ are non-zero, and $\mf{tr}^{\e}_{\e'}(\mathcal{W}_{\e}(\mb{l},\mb{c}))=\mathcal{W}_{\e'}(\mb{l},\mb{c})$.
Similarly, we have $\mf{tr}^{\e}_{\e'}(\mathcal{W}_{\e}(\mb{l}^*,\mb{c}^*))=\mathcal{W}_{\e'}(\mb{l}^*,\mb{c}^*)$. Since $\mathcal{W}_{\e}(\mb{l},\mb{c})\cong \mathcal{W}_{\e}(\mb{l}^*,\mb{c}^*)$, we have $\mathcal{W}_{\e'}(\mb{l},\mb{c})\cong \mathcal{W}_{\e'}(\mb{l}^*,\mb{c}^*)$.  

Then by the classification of finite-dimensional irreducible modules over $\mathcal{U}(\e')$, the quantum affine algebra of type $A_{N-1}^{(1)}$, we have $t=t^*$ where $\boldsymbol{l}\in\mathbb{Z}^{t}_{+},\,\boldsymbol{l}^*\in \mathbb{Z}^{t^*}_{+}$ and $({\mb l}^*,{\mb c}^*)$ is a permutation of $({\mb l},{\mb c})$.
\qed\vskip 2mm

By \eqref{eq:SWD functo on J ell} and Lemmas \ref{lem:universal mapping property of inverse limit} (with $\mc{C}(\e^{(k)})$ replaced by $\mc{C}_J^\ell(\e^{(k)})$) and \ref{lem:subcategory C J ell}(2), we have functors
\begin{equation*}\label{eq:functors F infty}
\begin{split}
\xymatrixcolsep{2pc}\xymatrixrowsep{3pc}\xymatrix{
\mc{F}_{\e^\infty,\ell} : R^{J}(\ell)\text{\rm -mod}_0  \ar@{->}[r] & \ \mc{C}^\ell_J(\e^\infty)},\\
\xymatrixcolsep{2pc}\xymatrixrowsep{3pc}\xymatrix{
\mc{F}_{\e^\infty}=\bigoplus_{\ell}\mc{F}_{\e^\infty,\ell} : R^{J}\text{\rm -mod}_0  \ar@{->}[r] & \ \mc{C}_J(\e^\infty)},
\end{split}
\end{equation*}
which are exact and monoidal.

\begin{prop}
The functor $\mc{F}_{\e^\infty,\ell}$ is an equivalence of categories, and hence so is $\mc{F}_{\e^\infty}$.
\end{prop}
\pf
It follows directly from Corollaries \ref{cor:equiv for large rank} and \ref{lem:subcategory C J ell}(2), and Lemma~\ref{lem:equiv-infty} with $\mc{C}\left(\e^{(k)}\right)$ replaced by $\mc{C}_J^\ell\left(\e^{(k)}\right)$.
\qed

\begin{cor}
The functor $\mc{F}_{\e^\infty}$ induces a one-to-one correspondence between the sets of isomorphism classes of irreducible objects in $R^{J}\text{\rm -mod}_0$ and $\mc{C}_J(\e^\infty)$ by sending
\begin{equation*}
\xymatrixcolsep{2pc}\xymatrixrowsep{3pc}\xymatrix{
\mathrm{hd}\left(L(a_{1},b_{1})\circ\cdots\circ L(a_{t},b_{t})\right)   \ar@{->}[r] & \ \mathcal{W}_{\e^\infty}({\mb l},{\mb c})},
\end{equation*}
where $((a_1,b_1),\dots,(a_t,b_t))$ is an ordered multisegment and $({\mb l},{\mb c})\in \mc{P}_J^+$ is given by ${\mb l}=(b_1-a_1+1,\dots,b_t-a_t+1)$ and ${\mb c}=(q^{-a_1-b_1},\dots,q^{-a_t-b_t})$.
\end{cor}

\subsection{Super duality}\label{subsec:super duality}
Let $\e^\infty$ and $(\epsilon^{(k)})_{k\ge 1}$ be as in \eqref{eq:e infinite} and \eqref{eq:cond-epsilon-infty}, respectively.
Let us describe more explicitly a connection between $\mc{C}_J(\e^\infty)$ and the category of finite-dimensional representations of the usual quantum affine algebras $\mathcal{U}(\e)$ of type $A$, that is, $\e=\e_{n|0}$ and $\e=\e_{0|n}$.

Let
\begin{equation*}
\un{\e}^\infty=(\un{\e}_i)_{i\in \mathbb{Z}_{>0}}= (0,0,0,\dots),  \quad \
\ov{\e}^\infty=(\ov{\e}_i)_{i\in \mathbb{Z}_{>0}}= (1,1,1,\dots).
\end{equation*}
Recall that $M_k$ and $N_k$ are the number of $0$'s and $1$'s in $\e^{(k)}$, respectively.
Let $\{r_k\}_{k\geq 1}$ and $\{s_k\}_{k\geq 1}$ be the strictly increasing sequences satisfying the following:
\begin{equation*}
\e_{k|0}=\un{\e}^{(k)}< \e^{(r_k)}, \quad \e_{0|k}=\ov{\e}^{(k)}< \e^{(s_k)}.
\end{equation*}

Given $\ell\in \mathbb{Z}_+$, let us define functors ${\mc{S}}_{k|0}$ and ${\mc{S}}_{0|k}$
\begin{equation}\label{eq:triangle-1} 
\xymatrixcolsep{3pc}\xymatrixrowsep{1pc}\xymatrix{
& \mc{C}^\ell_J(\e^\infty)   \ar@{->}_{{\mc{S}}_{k|0}}[dl]\ar@{->}^{{\mc{S}}_{0|k}}[dr] &  \\
\mc{C}^\ell_J(\e_{k|0}) & &  \mc{C}^\ell_J(\e_{0|k})
}
\end{equation}
for $k\geq 1$ by
\begin{equation*}
{\mc{S}}_{k|0}=\mf{tr}^{\e^{(r_k)}}_{\e_{k|0}}\circ \mf{tr}_{r_k},\quad
{\mc{S}}_{0|k}=\mf{tr}^{\e^{(s_k)}}_{\e_{0|k}}\circ \mf{tr}_{s_k}.  
\end{equation*}
%
Applying Lemma \ref{lem:universal mapping property of inverse limit} to \eqref{eq:triangle-1} with $\mc{C}\left(\e^{(k)}\right)$ replaced by $\mc{C}^\ell_J(\e_{k|0})$ and $\mc{C}^\ell_J(\e_{0|k})$, we obtain functors $\mc{S}_{\infty|0}$ and $\mc{S}_{0|\infty}$.
\begin{equation*} 
\xymatrixcolsep{3pc}\xymatrixrowsep{1pc}\xymatrix{
& \mc{C}^\ell_J(\e^\infty)   \ar@{->}_{\mc{S}_{\infty|0}}[dl]\ar@{->}^{\mc{S}_{0|\infty}}[dr] &  \\
\mc{C}^\ell_J(\un{\e}^\infty) & &  \mc{C}^\ell_J(\ov{\e}^\infty)
}
\end{equation*}
such that $\mf{tr}_k\circ \mc{S}_{\infty|0}\cong \mc{S}_{k|0}$ and $\mf{tr}_k\circ \mc{S}_{0|\infty}\cong \mc{S}_{0|k}$ for $k\geq 1$.
Again by Lemma \ref{lem:subcategory C J ell}(1), we have functors of monoidal categories
\begin{equation*} 
\xymatrixcolsep{3pc}\xymatrixrowsep{1pc}\xymatrix{
& \mc{C}_J(\e^\infty)   \ar@{->}_{\mc{S}_{\infty|0}}[dl]\ar@{->}^{\mc{S}_{0|\infty}}[dr] &  \\
\mc{C}_J(\un{\e}^\infty) & &  \mc{C}_J(\ov{\e}^\infty)
}.
\end{equation*}
Then we have the following, which can be viewed as a quantum affine analogue of {\em super duality} of type $A$ \cite{CL}. 

\begin{thm}\label{thm:super duality}
The functors $\mc{S}_{\infty|0}$ and $\mc{S}_{0|\infty}$ are equivalences of monoidal categories.
\end{thm}
\pf It is enough to show that $\mc{S}_{\infty|0}: \mc{C}^\ell_J(\e^\infty)\longrightarrow \mc{C}^\ell_J(\un{\e}^\infty)$ is an equivalences of categories for each $\ell\in\mathbb{Z}_+$. 
First, note that $\mc{C}^\ell_J(\e^\infty)\cong \mc{C}^\ell_J\left(\e^{(k)}\right)$ for $k>\ell$ by Lemma \ref{lem:equiv-infty} and Corollary \ref{cor:truncation is an equivalence}. Also, $\mc{C}^\ell_J\left(\e^{(k)}\right) \cong \mc{C}^\ell_J(\e_{k|0})$ for all sufficiently large $k$ by Corollary \ref{cor:truncation is an equivalence}.
Hence $\mc{S}_{k|0} : \mc{C}^\ell_J(\e^\infty) \longrightarrow \mc{C}^\ell_J(\e_{k|0})$ is an equivalence for all sufficiently large $k$.
This implies that  $\mc{S}_{\infty|0}=\varprojlim \mc{S}_{k|0} : \mc{C}^\ell_J(\e^\infty)\longrightarrow \mc{C}^\ell_J(\un{\e}^\infty)$ is an equivalence of categories by Lemma \ref{lem:equiv-infty}.
\qed\vskip 2mm

\begin{cor}\label{cor:super duality}
Given $\ell\in\mathbb{Z}_+$, the following are equivalences of categories for a sufficiently large $k$:
\begin{equation*} 
\xymatrixcolsep{3pc}\xymatrixrowsep{2pc}\xymatrix{
& \mc{C}^\ell_J(\e^\infty)   \ar@{->}_{{\mc{S}}_{k|0}}[dl]\ar@{->}_{\mf{tr}_k}[d]\ar@{->}^{{\mc{S}}_{0|k}}[dr] &  \\
\mc{C}^\ell_J(\e_{k|0}) & \mc{C}^\ell_J(\e^{(k)}) &  \mc{C}^\ell_J(\e_{0|k})
.}
\end{equation*}
\end{cor}

\subsection{Grothendieck rings}

Let $\e^\infty$ be the sequence given as in \eqref{eq:e infinite}. 
Let ${\rm K}(\mc{C}_J(\e^\infty))$ be the Grothendieck group of $\mc{C}_J(\e^\infty)$. Then 
\begin{equation*}
{\rm K}(\mc{C}_J(\e^\infty))=\bigoplus_{\ell\in\mathbb{Z}_+}{\rm K}(\mc{C}^\ell_J\left(\e^\infty)\right)
\end{equation*}
by Lemma \ref{lem:subcategory C J ell}, and it has a well-defined ring structure with multiplication $$[\,\mathbb{V}\,]\cdot[\,\mathbb{W}\,]=[\,\mathbb{V}\ot \mathbb{W}\,]\in {\rm K}(\mc{C}^{\ell+\ell'}_J\left(\e^\infty)\right),$$ 
for $[\,\mathbb{V}\,]\in {\rm K}(\mc{C}^\ell_J\left(\e^\infty)\right)$ and 
$[\,\mathbb{W}\,]\in {\rm K}(\mc{C}^{\ell'}_J\left(\e^\infty)\right)$.
By Theorem \ref{thm:super duality}, we have ring isomorphisms from ${\rm K}(\mc{C}_J(\e^\infty))$ to ${\rm K}(\mc{C}_J(\un{\e}^\infty))$ and ${\rm K}(\mc{C}_J(\ov{\e}^\infty))$.

Let $S=\{\,(l,a)\,|\,l\in \mathbb{N},\, a\in  l-1+2\mathbb{Z}\}$ and
\begin{equation*}
R=\mathbb{Z}\left[\,t_{l,a}\,|\,(l,a)\in S \, \right]
\end{equation*}
be the polynomial ring generated by $t_{l,a}\in S$. Put ${\rm deg}(t_{l,a})=l$ and let $R^\ell$ denote the subgroup of $R$ generated by monomials of degree $\ell\in \mathbb{Z}_+$.

\begin{prop}
There is an isomorphism of rings 
\begin{equation*}
\xymatrixcolsep{3pc}\xymatrixrowsep{0pc}\xymatrix{
{\rm K}(\mc{C}_J({\e}^\infty)) \ar@{->}[r] & \ R  \\
\left[\,\mathcal{W}_{l,\e^\infty}(q^a)\,\right] \ar@{|->}[r] & \ t_{l,a}},
\end{equation*}
where ${\rm K}(\mc{C}_J^\ell({\e}^\infty))$ maps onto $R^\ell$.
In particular, ${\rm K}\left(\mc{C}_J\left(\e^{(k)}\right)\right)$ is a homomorphic image of $R$ given by sending $t_{l,a}$ to $[\mathcal{W}_{l,\e^{(k)}}(q^{a})]$ for $k\geq 1$, where $(\e^{(k)})_{k\ge 1}$ is given as in \eqref{eq:cond-epsilon-infty}.
\end{prop}
\pf 
For $k\geq 1$,
it is well-known \cite[Corollary 2]{FR} that  ${\rm K}(\mc{C}_J(\e_{0|k+1}))$ is isomorphic to  
$$R_k:=\mathbb{Z}[\,t_{l,a}\,|\,l=1,\dots,k, a\in l-1+2\mathbb{Z}\,]\subset R,$$ 
as a ring, where $\left[\mathcal{W}_{l,\e_{0|k+1}}(q^{a})\right]$ is mapped to $t_{l,a}$, and ${\rm K}(\mc{C}^\ell_J(\e_{0|k+1}))\cong R_k^\ell:=R^\ell\cap R_k$ for $\ell\in\mathbb{Z}_+$. 
Since $\mf{tr}^{\epsilon_{0|k+1}}_{\epsilon_{0|k}}$ induces a map $R_{k}\longrightarrow R_{k-1}$ given by $t_{k,a}=0$ for $a\in l-1+2\mathbb{Z}$ and its restriction $R^\ell_{k}\longrightarrow R^\ell_{k-1}$ stabilizes for all sufficiently large $k$, we have
\begin{equation*}
{\rm K}(\mc{C}_J^\ell(\ov{\e}^\infty)) = \varprojlim {\rm K}(\mc{C}^\ell_J(\e_{0|k})) \cong R^\ell.
\end{equation*}
Hence we have an isomorphism of rings 
${\rm K}(\mc{C}_J(\ov{\e}^\infty)) \longrightarrow R$ which maps 
$\left[\mathcal{W}_{l,\ov{\e}^\infty}(q^a)\right]$ to $t_{l,a}$
for $(l,a)\in S$. 
Therefore, the first isomorphism from ${\rm K}(\mc{C}_J({\e}^\infty))$ to $R$ follows from
$\mc{S}_{0|\infty}(\mathcal{W}_{l,\e^\infty}(q^a))=\mathcal{W}_{l,\ov{\e}^\infty}(q^a)$. The second homomorphism from $R$ to ${\rm K}(\mc{C}_J\left(\e^{(k)}\right))$ is induced from $\mf{tr}_k$.
\qed  

\begin{cor}
Let $\epsilon\in\mathcal{E}$ and $(\mb{l},\mb{c})\in \mc{P}_J^+(\e)$ be given. If $\chi$ is a polynomial in $R$ such that 
$$[\mathcal{W}_{\e_{0|k}}(\mb{l},\mb{c})] =\chi (\{\,[\mathcal{W}_{l,\e_{0|k}}(q^a)]\,|\,(l,a)\in S\,\})$$ for a sufficiently large $k$, then the following holds in ${\rm K}(\mc{C}_J(\e))$:
$$[\mathcal{W}_{\e}(\mb{l},\mb{c})] =\chi \left(\{\,[\mathcal{W}_{l,\e}(q^a)]\,|\,(l,a)\in S\,\}\right).$$
\end{cor}

\begin{rem}{\rm
The isomorphism of $\mathbb{Q}$-algebras in \eqref{eq:iso tilde} induces an equivalence between $\mc{C}_J(\un{\e}^\infty)$ and $\mc{C}_J(\ov{\e}^\infty)$, by which we identify their Grothendieck groups. Under this identification, the equivalence $\mc{S}_{\infty|0}\circ \mc{S}_{0|\infty}^{-1}$ induces an involution on ${\rm K}(\mc{C}_J(\un{\e}^\infty))$, which can be viewed as an affine analogue of the involution on the ring of symmetric functions sending a Schur function to another one of conjugate shape.}
\end{rem}

\subsection{KR modules and $T$-system}
Let $\epsilon\in\mathcal{E}$ be given. 
For $r, s\in \mathbb{N}$ and $c\in \mathbf{k}^\times$, let
\begin{equation*}
\mathcal{W}_{\e}^{r,s}(c) :=  \mathcal{W}_{\e}(\mb{l},\mb{c})\quad \text{where}\quad 
\mb{l}=(r,\dots,r)\in\mathbb{Z}^{s}_{+} \text{ and } \mb{c}=c(q^{-2(s-1)},\dots, q^{-2},1).
\end{equation*}
If $\e=\e_{0|n}$ or $\e_{n|0}$, then $\mathcal{W}_{\e}^{r,s}(c)$ is the usual Kirillov-Reshetikhin module of type $A_{n-1}^{(1)}$ (see Remark \ref{rem:KR module} for more detail).

We note that $\mathcal{W}_{\e}^{r,s}(c)\cong V_\e((s^r))$ as a $\ring{\mathcal{U}}(\e)$-module if is not zero. This can be seen easily from the fact in case of $\e_{0|n}$, Theorem \ref{thm:branching multiplicity} and Corollary \ref{cor:super duality}. Moreover, $\mathcal{W}^{r,s}_{\epsilon}(c)\in\mathcal{C}_{J}(\epsilon)$ if and only if $c\in q^{r-1+2\mathbb{Z}}$.
 
\begin{prop}
There exists a short exact sequence in $\mc{C}(\e)$ 
\begin{equation*}
\xymatrixcolsep{1.5pc}\xymatrixrowsep{0pc}\xymatrix{
0 \ar@{->}[r] 
& \stackrel[r'=r\pm 1]{}{\bigotimes}\mathcal{W}_{\e}^{r',s}(cq^{-1}) \ar@{->}[r] 
&  \mathcal{W}_{\e}^{r,s}(c)\ot \mathcal{W}_{\e}^{r,s}(cq^{-2})  \ar@{->}[r] 
& \mathcal{W}_{\e}^{r,s+1}(c)\ot \mathcal{W}_{\e}^{r,s-1}(cq^{-2})  \ar@{->}[r] 
&  0 },
\end{equation*}
where $c\in \mathbf{k}^\times$. Hence the following holds in ${\rm K}(\mc{C}(\e))$:
\begin{equation*}
\left[\,\mathcal{W}_{\e}^{r,s}(c)\,\right][\,\mathcal{W}_{\e}^{r,s}(cq^{-2})\,] 
= [\,\mathcal{W}_{\e}^{r,s+1}(c)\,][\,\mathcal{W}_{\e}^{r,s-1}(cq^{-2})\,]
+ [\,\mathcal{W}_{\e}^{r-1,s}(cq^{-1})\,][\,\mathcal{W}_{\e}^{r+1,s}(cq^{-1})\,].
\end{equation*}
\end{prop}
\pf 
First we take $c=q^{r-1}$ so that all the modules in the sequence belong to $\mathcal{C}_{J}(\epsilon)$. In case of $\e=\e_{0|n}$, the existence of such short exact sequence is well-known, which is called the {\em $T$-system} \cite{Her06,KNS,Na03}.
Suppose that $\e$ is arbitrary. 
Take $\e^\infty$ as in \eqref{eq:e infinite} with an ascending chain of subsequence $(\e^{(k)})_{k\ge 1}$ such that $\e^{(k)}=\e$ for some $k$. 
We consider $\mathcal{W}_{\e^\infty}^{r,s}(c):=\mathcal{W}_{\e^\infty}(\mb{l},\mb{c})\in \mc{C}_J(\e^\infty)$.
Since we have such a short exact sequence for $\e_{0|m}$ for a sufficiently large $m$, and hence for $\e^\infty$ by Corollary \ref{cor:super duality}, the result follows by applying $\mf{tr}_k$.
Then by applying an automorphism $\tau_c$ of $\mathcal{U}(\e)$ ($c\in \bf{k}^\times$) given by $\tau_c(e_i)=c^{\delta_{i0}}e_i$, $\tau_c(f_i)=c^{-\delta_{i0}}f_i$ and $\tau_c(k_\mu)=k_\mu$ for $i\in I$ and $\mu\in P$ to the above exact sequence,
we may conclude the same result for arbitrary $c\in \mathbf{k}^\times$.
\qed\vskip 2mm

\begin{rem}\label{rem:KR module}
{\rm
(1)  
Let $W^{(r)}_{s,a}$ denote the Kirillov-Reshetikhin module over $U^{\prime}_q\left(A_{n-1}^{(1)}\right)$ associated to $r=1,\dots, n-1$, $s\geq 1$ and $a\in \mathbf{k}^\times$. 
The Drinfeld polynomial of $W^{(r)}_{s,a}$ is given by $(P_i(u))_{1\le i\le n-1}$, where
\begin{equation*} 
P_i(u)=
\begin{cases}
\prod_{k=1}^s(1-aq^{2k-2}u) & \text{if $i=r$},\\
1 & \text{if $i\neq r$}.
\end{cases}
\end{equation*}
Note that if $\e=\e_{0|n}$ and $r<n$, then $\mathcal{W}_{r,\e}(c)$ corresponds to $W^{(r)}_{1,a}$, where
\begin{equation*}
a= -o(r)(-q)^{-n} \td{c}
\end{equation*}
by Remark \ref{rem:homogeneous fundamental} and \cite[Remark 3.3]{Na04-2}.
Here $\td{\ \mbox{}\ }: {\bf k} \longrightarrow {\bf k}$ is an automorphism such that $\td{q}=-q^{-1}$ and $o: I\setminus\{0\}\rightarrow \{\pm 1\}$ is a map such that $o(i)=-o(j)$ if $a_{ij}=-1$.
Therefore, $\mathcal{W}_{\e}^{r,s}(c)$ corresponds to $W^{(r)}_{s,a}$.

(2) When $\e=\e_{n|0}$, we have a short exact sequence of tensor products of KR modules where the shapes corresponding to the $\ring{\mathcal{U}}(\e)$-highest weight are the conjugate of the ones appearing in case of $\e=\e_{0|n}$.

(3) In \cite{KO}, it is shown that the Kirillov-Reshetikhin module $\mathcal{W}_{\e}^{r,s}(c)$ has a crystal base with a suitable choice of $c$, when $\e$ is standard, that is, $\e=\e_{M|N}$. 
An explicit description of its crystal and the associated combinatorial $R$ matrix is also given in \cite{KO}.

}
\end{rem}

\appendix{}
\section{}

\subsection{Proof of Lemma \ref{lem:AKlemma}}\label{appendix A}

In this section we prove Lemma~\ref{lem:AKlemma}. The proof consists of direct calculations as indicated in \cite[Lemma B.1]{AK}, but we give details for the reader's convenience as it is little more involved.

We claim that there exists an exact sequence of the following form for each $\ell\geq 2$: 
\begin{equation}\label{eq:exact sequence}
\begin{tikzcd}[sep=1.4em] 
0 \arrow[r] 
& \mathcal{W}_{\ell,\e}(1) \arrow[r, "\psi_1"] 
& \mathcal{W}_{1,\e}(q^{1-\ell})\otimes \mathcal{W}_{\ell-1,\e}(q^{1}) \arrow[r, "\mathcal{R}"]  
& \mathcal{W}_{\ell-1,\e}(q^{1})\otimes \mathcal{W}_{1,\e}(q^{1-\ell}) \arrow[r,"\psi_2"] 
& \mathcal{W}_{\ell,\e}(1) \arrow[r] & 0, 
\end{tikzcd}
\end{equation} 
for some $\mathcal{U}(\e)$-linear maps $\psi_1$ and $\psi_2$, where $\mc{R}=\mc{R}^{\rm norm}_{1,\ell-1}(q^{-\ell})$.
Recall from Theorem \ref{thm:spectral decomposition}
that
\begin{equation}\label{eq:R matrix for appendix}
\mathcal{R}_{1,\ell-1}^{\mathrm{norm}}(z)
=\mathcal{P}_{1}+\frac{1-zq^{\ell}}{z-q^{\ell}}\mathcal{P}_{0},
\end{equation}
which is equal to $\mc{P}_1$ when $z=q^{-\ell}$.
We may assume that $\epsilon_{1}=0$. 
Indeed, the result for arbitrary $\epsilon$ follows once we choose $\epsilon^{\prime}>\epsilon$
with $\epsilon_{1}^{\prime}=0$ and apply the truncation functor $\mathfrak{tr}_{\epsilon}^{\epsilon^{\prime}}$
to the exact sequence for $\epsilon^{\prime}$, 
keeping Propositions~\ref{prop:truncation},~\ref{prop:truncation of poly and fund} and Lemma~\ref{lem:truncation of Verma and R matrix} in mind. 

Recall that when $\epsilon_{1}=0$, the $\ring{\mathcal{U}}(\e)$-highest weight vectors of $V((\ell))$ and
$V((\ell-1,1))$ in the decomposition
$\mathcal{W}_{1,\epsilon}(x)\otimes\mathcal{W}_{\ell-1,\epsilon}(y)
= V_\epsilon((\ell))\oplus V_\epsilon((\ell-1,1))$ are given by
\begin{equation}\label{eq:maximal vectors}
\ket{\be_{1}}\otimes\ket{(\ell-1)\be_{1}},\quad 
\ket{\be_{1}}\otimes\ket{(\ell-2)\be_{1}+\be_{2}}-q^{\ell-1}\ket{\be_{2}}\otimes\ket{(\ell-1)\be_{1}}
\end{equation}
respectively. On the other hand, when $\epsilon_{1}=1$ the highest weight vectors
become more complicated.

Let us define $\psi_1$ and $\psi_2$ by
\begin{align*}
& \psi_1(\left|\mathbf{m}\right\rangle) 
=\sum_{1\le k\le n}\left|\be_{k}\right\rangle \otimes\left|\mathbf{m}-\be_{k}\right\rangle \left(\left[m_{k}\right]\prod_{k<j\leq n}q^{m_{j}}\right)\\
& \psi_2\left(\left|\mathbf{m}\right\rangle \otimes\left|\be_{k}\right\rangle \right)  =\left|\mathbf{m}+\be_{k}\right\rangle  \prod_{k<j\le n}q^{-m_{j}}
\end{align*}
for $|\bf{m}\rangle$ and $1\le k\le n$.
Here we also understand $\left|\mathbf{m}\right\rangle=0$ whenever $\mathbf{m}\notin\mathbb{Z}_{+}^{n}(\epsilon)$ .
Note that when $\epsilon=(1^{N})$, $\psi_1$ and $\psi_2$ coincide with the maps in \cite[Lemma B.1]{AK} up to a constant multiple.

\begin{lem}
The maps $\psi_1$ and $\psi_2$ are $\mathcal{U}(\epsilon)$-linear.
\end{lem}
\pf Since the proof is rather straightforward, let us show that $\psi_1$ commutes with $e_i$ ($i\in I$), and leave the other details to the reader.\vskip 2mm

{\em Case 1}. Suppose that $i\in I\setminus\{0\}$.
First we have
\begin{equation}\label{eq:aux-1}
\begin{split}
e_{i}\psi_1\ket{\mathbf{m}}= & \sum\left[m_{k}\right]\prod_{j>k}q^{m_{j}}e_{i}\left(\ket{\be_{k}}\otimes\ket{\mathbf{m}-\be_{k}}\right)\\
= & \sum_{k\neq i,i+1}\left[m_{k}\right]\prod_{j>k}q^{m_{j}}\left[m_{i+1}\right]\ket{\be_{k}}\otimes\ket{\mathbf{m}-\be_{k}+\be_{i}-\be_{i+1}}\\
 & +\left[m_{i+1}\right]\prod_{j>i+1}q^{m_{j}}\left[m_{i+1}-1\right]\ket{\be_{i+1}}\otimes\ket{\mathbf{m}+\be_{i}-2\be_{i+1}}\\
 & +\left[m_{i+1}\right]\prod_{j>i+1}q^{m_{j}}\cdot q_{i}^{-m_{i}}q_{i+1}^{m_{i+1}-1}\ket{\be_{i}}\otimes\ket{\mathbf{m}-\be_{i+1}}\\
 & +\left[m_{i}\right]\prod_{j>i}q^{m_{j}}\left[m_{i+1}\right]\ket{\be_{i}}\otimes\ket{\mathbf{m}-\be_{i+1}}.
\end{split}
\end{equation}
Let $(\star)$ denote the sum of last two terms, that is,
\[
(\star)=\left[m_{i+1}\right]\prod_{j>i+1}q^{m_{j}}\cdot q_{i}^{-m_{i}}q_{i+1}^{m_{i+1}-1}\ket{\be_{i}}\otimes\ket{\mathbf{m}-\be_{i+1}}+\left[m_{i}\right]\prod_{j>i}q^{m_{j}}\left[m_{i+1}\right]\ket{\be_{i}}\otimes\ket{\mathbf{m}-\be_{i+1}}.
\]

Suppose first that $e_{i}\ket{\mathbf{m}}=0$.
Note that $e_{i}\ket{\mathbf{m}}=0$ if and only if $m_{i+1}=0$
or $m_{i+1}\neq0$, $m_{i}=1=\epsilon_{i}$. If $m_{i+1}=0$, then
$e_{i}\psi_1 \ket{\mathbf{m}}=0$. 
In the other case, $\ket{\mathbf{m}-\be_{k}+\be_{i}-\be_{i+1}}$
is nonzero if and only if $k=i$. So \eqref{eq:aux-1} is equal to 
\[
(\star)=\left[m_{i+1}\right]\prod_{j>i+1}q^{m_{j}}\left(\left[m_{i}\right]q^{m_{i+1}}+q_{i}^{-m_{i}}q_{i+1}^{m_{i+1}-1}\right)\ket{\be_{i}}\otimes\ket{\mathbf{m}-\be_{i+1}}.
\]
Since
\begin{align*}
\left[m_{i}\right]q^{m_{i+1}}+q_{i}^{-m_{i}}q_{i+1}^{m_{i+1}-1} & =\left[1\right]q^{m_{i+1}}+(-q)q_{i+1}^{m_{i+1}-1}\\
 & =
 \begin{cases}
q^{m_{i+1}}+(-q)q^{m_{i+1}-1} & \text{if }\epsilon_{i+1}=0,\\
q+(-q)(-q^{-1})^{0} & \text{if }\epsilon_{i+1}=1=m_{i+1},
\end{cases}\\
& =0,
\end{align*}
we have $e_{i}\psi_1 \ket{\mathbf{m}}=0=\psi_1 e_{i}\ket{\mathbf{m}}$ whenever
$e_{i}\ket{\mathbf{m}}=0$.

Next suppose that $e_{i}\ket{\mathbf{m}}\neq 0$ (necessarily $m_{i+1}\neq0$).
We have
\begin{align*}
\psi_1 e_{i}\ket{\mathbf{m}}= & \left[m_{i+1}\right]\psi_1 \ket{\mathbf{m}+\be_{i}-\be_{i+1}}\\
= & \left[m_{i+1}\right]\sum_{k\neq i,i+1}\left[m_{k}\right]\prod_{j>k}q^{m_{j}}\ket{\be_{k}}\otimes\ket{\mathbf{m}-\be_{k}+\be_{i}-\be_{i+1}}\\
 & +\left[m_{i+1}\right]\left[m_{i+1}-1\right]\prod_{j>i+1}q^{m_{j}}\ket{\be_{i+1}}\otimes\ket{\mathbf{m}-2\be_{i+1}+\be_{i}}\\
 & +\left[m_{i+1}\right]\left[m_{i}+1\right]q^{-1}\prod_{j>i}q^{m_{j}}\ket{\be_{i}}\otimes\ket{\mathbf{m}-\be_{i+1}}.
\end{align*}
It is equal to \eqref{eq:aux-1} if
\begin{equation}\label{eq:aux-2}
(\star)=\left[m_{i+1}\right]\left[m_{i}+1\right]q^{-1}\prod_{j>i}q^{m_{j}}\ket{\be_{i}}\otimes\ket{\mathbf{m}-\be_{i+1}}.
\end{equation}
Indeed, we have two possibilities: either $m_{i}=0$ or $m_{i}\neq0$ with $\epsilon_{i}=0$.
In the first case, we have
\begin{align*}
(\star) & =\left[m_{i+1}\right]\prod_{j>i+1}q^{m_{j}}\cdot q_{i+1}^{m_{i+1}-1}\ket{\be_{i}}\otimes\ket{\mathbf{m}-\be_{i+1}},
\end{align*}
and the product can be written as
\begin{align*}
\prod_{j>i+1}q^{m_{j}}\cdot q_{i+1}^{m_{i+1}-1} & =\begin{cases}
\prod_{j>i}q^{m_{j}}\cdot q^{-1} & \text{if }\epsilon_{i+1}=0\\
\prod_{j>i+1}q^{m_{j}} & \text{if }\epsilon_{i+1}=1=m_{i+1}.
\end{cases}\\
 & =\prod_{j>i}q^{m_{j}}\cdot q^{-1},
\end{align*}
which implies \eqref{eq:aux-2}.
In the other case, as $\prod_{j>i+1}q^{m_{j}}\cdot q_{i+1}^{m_{i+1}-1}=\prod_{j>i}q^{m_{j}}\cdot q^{-1}$
by the same reason, we have
\begin{align*}
(\star) & =\prod_{j>i}q^{m_{j}}\cdot q^{-1}\ket{\be_{i}}\otimes\ket{\mathbf{m}-\be_{i+1}}\left[m_{i+1}\right]\left(\left[m_{i}\right]q+q^{-m_{i}}\right)\\
 & =\prod_{j>i}q^{m_{j}}\cdot q^{-1}\ket{\be_{i}}\otimes\ket{\mathbf{m}-\be_{i+1}}\left[m_{i+1}\right]\left[m_{i}+1\right].
\end{align*}
Hence \eqref{eq:aux-2} holds.

{\em Case 2}. Suppose that $i=0$.
The proof is similar except that we should consider spectral parameters.
First, we have
{\allowdisplaybreaks
\begin{align*}
e_{0}\psi_1 \ket{\mathbf{m}}= & \sum\left[m_{k}\right]\prod_{j>k}q^{m_{j}}e_{0}\left(\ket{\be_{k}}\otimes\ket{\mathbf{m}-\be_{k}}\right)\\
= & \sum_{k\neq 1,n}\left[m_{k}\right]\prod_{j>k}q^{m_{j}}\left[m_{1}\right]\ket{\be_{k}}\otimes\ket{\mathbf{m}-\be_{k}+\be_{n}-\be_{1}}\cdot q\\
 & +\left[m_{1}\right]\prod_{j>1}q^{m_{j}}\left[m_{1}-1\right]\ket{\be_{1}}\otimes\ket{\mathbf{m}-2\be_{1}+\be_{n}}\cdot q\\
 & +\left[m_{1}\right]\prod_{j>1}q^{m_{j}}\cdot q_{n}^{-m_{n}}q_{1}^{m_{1}-1}\ket{\be_{n}}\otimes\ket{\mathbf{m}-\be_{1}}\cdot q^{1-\ell}\\
 & +\left[m_{n}\right]\left[m_{1}\right]\ket{\be_{n}}\otimes\ket{\mathbf{m}-\be_{1}}\cdot q.
\end{align*}}
Note that since $\ell=\sum m_{j}$, we have in the third term above
\[
\prod_{j>1}q^{m_{j}}\cdot q_{n}^{-m_{n}}q_{1}^{m_{1}-1}q^{1-\ell}=q^{-m_{1}}q_{n}^{-m_{n}}q_{1}^{m_{1}-1}q. 
\]
Similarly we have
\begin{align*}
\psi_1 e_{0}\ket{\mathbf{m}}= & \left[m_{1}\right]\psi_1 \ket{\mathbf{m}+\be_{n}-\be_{1}}\cdot1\\
= & \left[m_{1}\right]\sum_{k\neq 1,n}\left[m_{k}\right]\prod_{j>k}q^{m_{j}}\cdot q\ket{\be_{k}}\otimes\ket{\mathbf{m}-\be_{k}+\be_{n}-\be_{1}}\\
 & +\left[m_{1}\right]\left[m_{1}-1\right]\prod_{j>1}q^{m_{j}}\cdot q\ket{\be_{1}}\otimes\ket{\mathbf{m}-2\be_{1}+\be_{n}}\\
 & +\left[m_{1}\right]\left[m_{n}+1\right]\ket{\be_{n}}\otimes\ket{\mathbf{m}-\be_{1}}.
\end{align*}
Now the same argument applies as in {\em Case 1}. If $e_{0}\ket{\mathbf{m}}=0$, then
either $m_{1}=0$ or $m_{1}\neq0$ with $m_{n}=1=\epsilon_{n}$. 
In the first case, we clearly have $\psi_1 e_{0}\ket{\mathbf{m}}=e_{0}\psi_1 \ket{\mathbf{m}}=0$. 
In the latter case, we have
\[
e_{0}\psi_1 \ket{\mathbf{m}}=\left[m_{1}\right]\left(\left[1\right]q+q^{-m_{1}}(-q^{-1})^{-1}q_{1}^{m_{1}-1}q\right)\ket{\be_{n}}\otimes\ket{\mathbf{m}-\be_{1}}=0
\]
as $q-qq^{1-m_{1}}q_{1}^{m_{1}-1}$ vanishes regardless of $\epsilon_{1}$. 

Next, if $e_{0}\ket{\mathbf{m}}\neq0$ and $m_{1}\neq0$, then again
we have $\psi_1 e_{0}\ket{\mathbf{m}}=e_{0}\psi_1 \ket{\mathbf{m}}$ since
\begin{align*}
\left[m_{1}\right]\left(\left[m_{n}\right]q+q^{-m_{1}}q_{n}^{-m_{n}}q_{1}^{m_{1}-1}q\right) & =\begin{cases}
\left[m_{1}\right]\left(0+q^{-m_{1}}q_{1}^{m_{1}-1}q\right) & \text{if }m_{n}=0\\
\left[m_{1}\right]\left(\left[m_{n}\right]q+q^{-m_{n}}q^{-m_{1}}q_{1}^{m_{1}-1}q\right) & \text{if }m_{n}\neq0,\,\epsilon_{n}=0
\end{cases}\\
& =\left[m_{1}\right]\left[m_{n}+1\right].
\end{align*}
This completes the proof.
\qed

\begin{lem}
We have
\begin{itemize}
\item[(1)] $\psi_1$ is injective and $\mathcal{R}\circ \psi_1=0$,

\item[(2)] $\psi_2$ is surjective and $\psi_2\circ\mathcal{R}=0$.
\end{itemize}
\end{lem}

\pf (1) It is clear that $\psi_1$ is injective since $\psi_1$ is non-zero and $\mathcal{W}_{\ell,\e}(1)$ is irreducible.

By definition, we have $\psi_1(\ket{\ell\be_{1}})=C\ket{\be_{1}}\otimes\ket{(\ell-1)\be_{1}}=C v_1$ for a nonzero constant $C$ \eqref{eq:maximal vectors}. 
The $\ring{\mathcal{U}}(\e)$-highest weight vector $v_1$ is sent to zero by $\mathcal{R}$ since $\mc{R}=\mc{P}_1$ by \eqref{eq:R matrix for appendix}. 
This implies that $\mathcal{R}\circ \psi_1=0$.
 
(2) Since $\psi_2$ is non-zero and $\mathcal{W}_{\ell,\e}(1)$ is irreducible, it is surjective. 
Note that $v_2=\ket{(\ell-1)\be_{1}}\otimes\ket{\be_{2}}-q\ket{(\ell-2)\be_{1}+\be_{2}}\otimes\ket{\be_{1}}$ generates $\mathrm{Im}\mathcal{R}$, which is isomorphic to $V((\ell-1,1))$. Since $\psi_2(v_2)=0$, we have $\psi_2\circ\mathcal{R}=0$.
\qed
\begin{lem}
The sequence \eqref{eq:exact sequence} is exact.
\end{lem}

\pf By the previous lemmas and the universal mapping properties of ${\rm Ker}$ and ${\rm Coker}$, we have the following commutative diagram of $\mathcal{U}(\e)$-modules:
\begin{equation*}
\begin{tikzcd}[sep=1.9em]
& & & & \mathcal{W}_{\ell,\e}(1) \arrow[r] & 0 \\
0 \arrow[r] & {\rm Ker} \mathcal{R} \arrow[r] 
& \mathcal{W}_{1,\e}(x)\otimes\mathcal{W}_{\ell-1,\e}(y) \arrow[r, "\mathcal{R}"]
& \mathcal{W}_{\ell-1,\e}(y)\otimes\mathcal{W}_{1,\e}(x) \arrow[r] \arrow[ru,"\psi_2"]
& \mathrm{Coker} \mathcal{R} \arrow[r] \arrow[u, dashed, two heads] & 0 \\
0 \arrow[r] & \mathcal{W}_{\ell,\e}(1) \arrow[u, dashed, hook] \arrow[ru,"\psi_1"'] & & & & \\
\end{tikzcd}
\end{equation*}
Hence two vertical arrows are isomorphisms. This implies that \eqref{eq:exact sequence} is exact.
\qed

\subsection{Proof of Theorem \ref{thm:equiv-QASWD}}\label{appendix B}

We assume that $\ell < n$. Put $\mc{F}=\mc{F}^*_{\e,\ell}$.
We first show that
\begin{equation}\label{eq:CP equivalence}
\xymatrixcolsep{2pc}\xymatrixrowsep{0pc}\xymatrix{
\mc{F}_{\e,\ell}^{\ast} : \aha\text{{\rm -mod}} \ar@{->}[r] & \ \mc{C}^\ell(\e) \\
M \ar@{|->}[r] & \ \mathcal{V}^{\otimes\ell}\otimes_{\fha}M}
\end{equation}
is an equivalence of categories.
We almost follow the arguments in \cite[Section 4.3 - 4.6]{CP} except {a part of Lemma \ref{lem:B.2}}. 
The following easy lemma is essential for the later computation.
\begin{lem}[cf. {\cite[Lemma 4.3]{CP}}]\label{lem:pf-aff-SWD}\mbox{}
\begin{enumerate}
\item Let $M$ be a finite-dimensional $\fha$-module. If $v\in\mathcal{V}^{\otimes\ell}$
has non-zero components in each isotypical component of $\mathcal{J}_{\ell}(M)$,
then the $\mathbf{k}$-linear map
\[
\xymatrixcolsep{2pc}\xymatrixrowsep{0pc}\xymatrix{
 M \ar@{->}[r] & \ \mathcal{V}^{\otimes\ell}\otimes_{\fha}M=\mathcal{J}_{\ell}(M) \\
m \ar@{|->}[r] & \ v \ot m},
\]
is injective.

\item Let $\left\{\,v_{i}:=\ket{\mathbf{e}_{i}}\,|\,{i=1,\dots,n}\,\right\}$
be the standard basis of $\mathcal{V}$. If $i_{1},\dots,i_{\ell}\in\left\{ 1,\dots,n \right\} $
are distinct, then the $\mathring{\mathcal{U}}(\epsilon)$-module
$\mathcal{V}^{\otimes\ell}$ is generated by a single vector $v_{i_{1}}\otimes\cdots\otimes v_{i_{\ell}}$.
In particular, the vector satisfies the condition in (1).
\end{enumerate}
\end{lem}
 
We first prove that $\mathcal{F}$ is essentially surjective.
Suppose that $W\in\mc{C}^\ell(\e)$ is given. By Theorem \ref{thm:finite-SWD},
there exists a $\fha$-module $M$ for which $W\cong\mathcal{J}_{\ell}(M)=\mathcal{V}^{\otimes\ell}\otimes M$
as a $\mathring{\mathcal{U}}(\epsilon)$-module. We shall extend the $\fha$-action on $M$ to $\aha$ so that 
$W\cong\mathcal{V}^{\otimes\ell}\otimes_{\fha}M\cong V_{\mathbb{O}}^{\otimes\ell}\ot_{\mathbb{O}\aha}M$
as a $\mathcal{U}(\epsilon)$-module.

For $1\leq j\leq\ell$, set
$v^{(j)}=v_{2}\otimes\cdots\otimes v_{j}\otimes v_{n}\otimes v_{j+1}\otimes\cdots\otimes v_{\ell}$. 
Regarding $\mathcal{V}^{\otimes \ell}\otimes_{\fha}M$ as a $\mathcal{U}(\epsilon)$-module, the weight of $f_{0}\left(v^{(j)}\otimes m\right)$ is $\delta_{1}+\cdots+\delta_{\ell}\in P$.
As 
\begin{equation}\label{eq:basis-1}
\left\{ v_{i_{1}}\otimes\cdots\otimes v_{i_{\ell}}\,|\,1\le i_{1},\dots,i_{\ell}\le\ell  \text{ are distinct}\right\} 
\end{equation}
is a basis of $\left(\mathcal{V}^{\otimes\ell}\right)_{\delta_{1}+\cdots+\delta_{\ell}}$,
we can write as
\begin{equation}\label{eq:basis-2}
f_{0}\left(v^{(j)}\otimes m\right)=\sum_{\mb{i}}\left(v_{i_{1}}\otimes\cdots\otimes v_{i_{\ell}}\right)\otimes m_{\mb{i}},
\end{equation}
where the sum is over $\mb{i}=\left(i_{1},\dots,i_{\ell}\right)$ such that $v_{i_{1}}\otimes\cdots\otimes v_{i_{\ell}}$ belongs to \eqref{eq:basis-1}, and $m_{\mb{i}}\in M$.
In fact, considering the $\fha$-action by $\mathcal{R}$ in Proposition \ref{thm:finite-SWD}, for each $\mb{i}=\left(i_{1},\dots,i_{\ell}\right)$ in \eqref{eq:basis-2}, there exists $h_{\mb{i}}\in\fha$
such that 
\[
v_{i_{1}}\otimes\cdots\otimes v_{i_{\ell}}=
\left(v_{2}\otimes\cdots\otimes v_{j}\otimes v_{1}\otimes v_{j+1}\otimes\cdots\otimes v_{\ell}\right)h_{\mb{i}}.
\]
Hence \eqref{eq:basis-2} is reduced to 
\begin{equation}\label{eq:basis-3}
f_{0}(v^{(j)}\otimes m)=\left(v_{2}\otimes\cdots\otimes v_{j}\otimes v_{1}\otimes v_{j+1}\otimes\cdots\otimes v_{\ell}\right)\otimes m^{\prime},
\end{equation}
for some $m^{\prime}\in M$. 
By Lemma~\ref{lem:pf-aff-SWD}, such $m^{\prime}$ is unique.
Therefore we obtain a $\mathbf{k}$-linear endomorphism $\alpha_{j}^{-}\in\mathrm{End}(\fha)$
sending $m$ to $m^{\prime}$. Considering $e_{0}$-action instead
yields $\alpha_{j}^{+}$. 
So we have
\begin{equation}\label{eq:0 action}
\begin{split}
e_{0}\left(v^{(j)}\otimes m\right) & =\left(\Delta_{j}(e_{0})v^{(j)}\right)\otimes\alpha_{j}^{+}(m)=\sum_{1\le i\le \ell}\left(\Delta_{i}(e_{0})v^{(j)}\right)\otimes\alpha_{i}^{+}(m),\\
f_{0}\left(v^{(j)}\otimes m\right) & =\left(\Delta_{j}(f_{0})v^{(j)}\right)\otimes\alpha_{j}^{-}(m)=\sum_{1\le i\le \ell}\left(\Delta_{i}(f_{0})v^{(j)}\right)\otimes\alpha_{i}^{-}(m),
\end{split}
\end{equation}
where $\Delta_{i}(e_{0})$ and $\Delta_{i}(f_{0})$ are given by
\begin{align*}
\Delta_{i}(e_{0}) & = 1^{\otimes i-1}\otimes e_{0}\otimes\left(k_{0}^{-1}\right)^{\otimes\ell-i},\\
\Delta_{i}(f_{0}) & = k_{0}^{\otimes i-1}\otimes f_{0}\otimes1^{\otimes\ell-i},
\end{align*}
acting on $\mathcal{V}^{\otimes\ell}$. 
Note that $\Delta_{i}(e_{0})v^{(j)}=0$ unless $i=j$. 
Indeed, $v^{(j)}$ in \eqref{eq:0 action} can be replaced by arbitrary $v\in \mc{V}^{\ot \ell}$.
\begin{lem}\label{lem:B.2}
For $v\in\mathcal{V}^{\otimes\ell}$ and $m\in M$, we have
\begin{align*}
e_{0}(v\otimes m) & =\sum_{1\le j\le \ell}\left(\Delta_{j}(e_{0})v\right)\otimes\alpha_{j}^{+}(m),\\
f_{0}(v\otimes m) & =\sum_{1\le j\le \ell}\left(\Delta_{j}(f_{0})v\right)\otimes\alpha_{j}^{-}(m).
\end{align*}
\end{lem}

\pf
We only prove the case for $f_0$ since the other case is similar.
Take
$v=v_{i_{1}}\otimes\cdots\otimes v_{i_{\ell}}$. If none of $i_{j}$
is equal to $n$, then $\Delta_{j}(f_{0})v=0$ for any $j$.
On the other hand, we have $f_{0}(v\otimes m)=0$ since $\delta_{i_{1}}+\cdots+\delta_{i_{\ell}}-\delta_{n}+\delta_{1}$
is not a weight of $\mathcal{V}^{\otimes\ell}$. Hence the identity
holds.

For each pair of sequences
\begin{align*}
\mb{j} & =(j_{1}<j_{2}<\cdots<j_{r}),\quad
\mb{j}' =(j_{1}^{\prime}<j_{2}^{\prime}<\cdots<j_{s}^{\prime})
\end{align*}
in $\left\{ 1,\dots,\ell\right\}$, which are disjoint, let $\mc{V}^{(\mb{j},\mb{j}')}$
be the subspace of $\mathcal{V}^{\otimes\ell}$ spanned by vectors
of the form $v_{i_{1}}\otimes\cdots\otimes v_{i_{\ell}}$, where $i_{j_{t}}=1$
($1\leq t\leq r$), $i_{j_{t}^{\prime}}=n$ ($1\leq t\leq s$) and
$i_{j}\neq1,n$ for others. 
Clearly $\mc{V}^{\ot\ell}=\bigoplus\mc{V}^{(\mb{j},\mb{j}^{\prime})}$,
so that we may prove the identity for $v$ in each $\mathcal{V}^{(\mathbf{j},\mathbf{j}^{\prime})}$.

In addition, it is enough to check the identity for $v=v_{i_{1}}\otimes\cdots\otimes v_{i_{\ell}}\in\mathcal{V}^{(\mb{j},\mb{j}^{\prime})}$
with no $v_{2},\,\dots,\,v_{n-1}$ appearing more than once, because
of Lemma~\ref{lem:pf-aff-SWD}(2) (with respect to the subalgebra
of $\mathring{\mathcal{U}}(\epsilon)$ generated by $e_{i},\,f_{i}$
and $k_{i}^{\pm1}$ for $i=2,\dots,n-1$). There is always such a vector
since $\ell< n$.

We shall prove the identity by induction on $s$.  
We start with $s=1$, and use induction on $r$.
The case when $r=0$ and $s=1$ has already been done when we define $\alpha_{j}^{\pm}$ with $v=v^{(j)}$. 

Suppose that it is true for $r-1$. 
Choose $v=v_{i_{1}}\otimes\cdots\otimes v_{i_{\ell}}\in\mathcal{V}^{(\mb{j},\mb{j}')}$
such that only $v_{3},\,\dots,\,v_{n-1}$ appear as a factor of $v$
without repetition (which is possible as $s,r\geq 1$). 
Let $v^{\prime}$ be the vector obtained from $v$ by replacing the last $v_{1}$ (that
is, $v_{j_{r}}$) by $v_{2}$ so that $v^{\prime}$ has one less $v_{1}$
than $v$. By our choice of $v$, $e_{1}v^{\prime}=v$. Then we compute
as
\begin{align*}
f_{0}(v\otimes m) & =f_{0}e_{1}(v^{\prime}\otimes m) =e_{1}f_{0}(v^{\prime}\otimes m) =e_{1}\sum_{j}\left(\Delta_{j}(f_{0})v^{\prime}\right)\otimes\alpha_{j}^{-}(m)\\
 & =e_{1}\left(q_{1}^{-\left|\left\{ t\,|\,t<r,\, j_{t}<j_{1}^{\prime}\right\} \right|}v^{\prime\prime}\right)\otimes\alpha_{j_{1}'}^{-}(m)\\
 & =q_{1}^{-\left|\left\{ t\,|\,t<r,\,j_{t}<j_{1}^{\prime}\right\} \right|}(e_{1}v^{\prime\prime})\otimes\alpha_{j_{1}'}^{-}(m)\\
 & =q_{1}^{-\left|\left\{ t\,|\,t<r,\,j_{t}<j_{1}^{\prime}\right\} \right|}q_{1}^{-\delta(j_{r}<j_{1}^{\prime})}\left[(1^{\otimes j_{r}-1}\otimes e_{1}\otimes1^{\otimes\ell-j_{r}})v^{\prime\prime}\right]\otimes\alpha_{j_{1}'}^{-}(m)\\
 & =q_{1}^{-\left|\left\{ t\,|\,t\leq r,\,j_{t}<j_{1}^{\prime}\right\} \right|}\left[(1^{\otimes j_{r}-1}\otimes e_{1}\otimes1^{\otimes\ell-j_{r}})v^{\prime\prime}\right]\otimes\alpha_{j_{1}'}^{-}(m)\\
 & =\left(\Delta_{j_{1}'}(f_{0})v\right)\otimes\alpha_{j_{1}'}^{-}(m)=\sum_{j}\left(\Delta_{j}(f_{0})v\right)\otimes\alpha_{j}^{-}(m).
\end{align*}
Here the third equality follows from induction hypothesis on $r$,
$v^{\prime\prime}$ is the resulting vector of replacing (the unique)
$v_{n}$ factor of $v^{\prime}$ by $v_{1}$, the last equality
holds since $v$ has exactly one $v_{n}$ factor, and $\delta(P)$ is 1
if the statement $P$ is true and 0 otherwise.

Now assume the result for $s-1$ and let us prove it for $s\geq2$.
Choose $v=v_{i_{1}}\otimes\cdots\otimes v_{i_{\ell}}\in\mathcal{V}^{(\mb{j},\mb{j}')}$
such that $v_{n-1}$ does not appear as a factor of $v$ and for each $i=2,\,\dots,\,n-2$, $v_i$ occurs at most once (which is possible as $s\geq2$). 
We shall compute $\left[e_{n-1},f_{n-1}\right]f_{0}(v\otimes m)$ in two different ways.
 
We first have
\begin{equation}\label{eq:0 action-2}
\left[e_{n-1},f_{n-1}\right]f_{0}(v\otimes m)=\frac{q_{n}^{1-s}-q_{n}^{s-1}}{q-q^{-1}}f_{0}(v\otimes m),
\end{equation}
since $\left[e_{n-1},f_{n-1}\right]=\frac{k_{n-1}-k_{n-1}^{-1}}{q-q^{-1}}$ and
the weight of $f_{0}(v\otimes m)$ is 
\[
\sum_{i_{k}\neq n-1,n}\delta_{i_{k}}+s\delta_{n}+\delta_{1}-\delta_{n}.
\] 

Next, by similar arguments for \eqref{eq:basis-3}, $f_{0}\left(v\otimes m\right)$
can be written as a sum of $v_{\mb{k}}\otimes m_{\mb{k}}$ for some $m_{\mb{k}}\in M$
and $v_{\mb{k}}=v_{k_{1}}\otimes\cdots\otimes v_{k_{\ell}}$ with none
of $v_{k_{i}}$ is equal to $v_{n-1}$. Hence $f_{n-1}f_{0}\left(v\otimes m\right)=0$
and so 
\begin{equation}\label{eq:0 action-3}
\left[e_{n-1},f_{n-1}\right]f_{0}(v\otimes m)=-f_{n-1}e_{n-1}f_{0}(v\otimes m)=-f_{n-1}f_{0}e_{n-1}(v\otimes m).
\end{equation}
We first compute
\[
e_{n-1}(v\otimes m)=\left(e_{n-1}v\right)\otimes m=\left(q_{n}^{s-1}v^{\prime,1}+q_{n}^{s-2}v^{\prime,2}+\cdots+v^{\prime,s}\right)\otimes m,
\]
where $v^{\prime,p}$ is obtained from $v$ by replacing $j_{p}^{\prime}$-th
factor (which is $v_{n}$) by $v_{n-1}$. The vector $v^{\prime,p}$
has one less $v_{n}$'s than $v$, so that the induction hypothesis deduces
\[
f_{0}e_{n-1}\left(v\otimes m\right)
=\sum_{p=1}^{s}q_{n}^{s-p}f_{0}(v^{\prime,p}\otimes m)
=\sum_{p=1}^{s}q_{n}^{s-p}\sum_{t}\left(\Delta_{t}(f_{0})v^{\prime,p}\right)\otimes\alpha_{t}^{-}(m).
\]
By definition of $\Delta_{t}(f_{0})$, 
\[
\Delta_{t}(f_{0})v^{\prime,p}=
\begin{cases}
q_{n}^{u-1-\delta(u>p)}q_{1}^{-\left|\left\{ k|j_{k}<t\right\} \right|}v^{\prime\prime,p,u} & \text{if }t=j_{u}^{\prime}\text{ for some }u\neq p,\\
0 & \text{otherwise}.
\end{cases}
\]
where $v^{\prime\prime,p,u}$ is obtained from $v^{\prime,p}$ by
replacing $j_{u}^{\prime}$-th factor (which is $v_{n}$) by $v_{1}$.
Since any nonzero $v^{\prime\prime,p,u}$ has exactly one $v_{n-1}$,
\begin{equation*}
\begin{split}
f_{n-1}f_{0}e_{n-1}\left(v\otimes m\right) 
& =f_{n-1}\sum_{p=1}^{s}q_{n}^{s-p}\sum_{t}\left(\Delta_{t}(f_{0})v^{\prime,p}\right)\otimes\alpha_{t}^{-}(m)\\
& =\sum_{p=1}^{s}q_{n}^{s-p}\left[\sum_{u\neq p}q_{n}^{u-1-\delta(u>p)}q_{1}^{-\left|\left\{ k|j_{k}<j_{u}^{\prime}\right\} \right|}\left(f_{n-1}v^{\prime\prime,p,u}\right)\otimes\alpha_{j_{u}^{\prime}}^{-}(m)\right]\\
& =\sum_{p=1}^{s}q_{n}^{s-p}\left[\sum_{u\neq p}q_{n}^{u-1-\delta(u>p)}q_{1}^{-\left|\left\{ k|j_{k}<j_{u}^{\prime}\right\} \right|}q_{n}^{1+\delta(u<p)-p}v^{u}\otimes\alpha_{j_{u}^{\prime}}^{-}(m)\right]
\end{split}
\end{equation*}
where $v^{u}$ is obtained from $v$ by replacing $j_{u}^{\prime}$-th
factor (which is $v_{n}$) by $v_{1}$. Now for $1\leq u\leq s$,
the coefficient of $v^{u}\otimes\alpha_{j_{u}^{\prime}}^{-}(m)$ is
\begin{align*}
 & \sum_{p<u}q_{n}^{s-p}q_{n}^{u-2}q_{1}^{-\left|\left\{ k|j_{k}<j_{u}^{\prime}\right\} \right|}q_{n}^{1-p}+\sum_{p>u}q_{n}^{s-p}q_{n}^{u-1}q_{1}^{-\left|\left\{ k|j_{k}<j_{u}^{\prime}\right\} \right|}q_{n}^{2-p}\\
= & q_{1}^{-\left|\left\{ k|j_{k}<j_{u}^{\prime}\right\} \right|}\left[\sum_{p=1}^{u-1}q_{n}^{s-p}q_{n}^{u-2}q_{n}^{1-p}+\sum_{p=u}^{s-1}q_{n}^{s-p-1}q_{n}^{u-1}q_{n}^{1-p}\right]\\
= & q_{1}^{-\left|\left\{ k|j_{k}<j_{u}^{\prime}\right\} \right|}\sum_{p=1}^{s-1}q_{n}^{s+u-1-2p}=q_{n}^{u-1}q_{1}^{-\left|\left\{ k|j_{k}<j_{u}^{\prime}\right\} \right|}\frac{q_{n}^{s-1}-q_{n}^{1-s}}{q_{n}-q_{n}^{-1}}.
\end{align*}

Finally, combining the computation of \eqref{eq:0 action-2} and \eqref{eq:0 action-3}
we obtain
\begin{align*}
f_{0}\left(v\otimes m\right) &=\frac{q-q^{-1}}{q_{n}^{s-1}-q_{n}^{1-s}}\left[f_{n-1},e_{n-1}\right]f_{0}\left(v\otimes m\right)  =\frac{q-q^{-1}}{q_{n}^{s-1}-q_{n}^{1-s}}f_{n-1}f_{0}e_{n-1}\left(v\otimes m\right)\\
 & =\sum_{u=1}^{s}\left(\Delta_{j_{u}^{\prime}}(f_{0})v\right)\otimes\alpha_{j_{u}^{\prime}}^{-}(m) =\sum_{t=1}^{\ell}\left(\Delta_{t}(f_{0})v\right)\otimes\alpha_{t}^{-}(m).
\end{align*}
since $\Delta_{j_{u}^{\prime}}(f_{0})v=q_{n}^{u-1}q_{1}^{-\left|\left\{ k|j_{k}<j_{u}^{\prime}\right\} \right|}v^{u}$
and $\Delta_{t}(f_{0})v=0$ for $t\neq j_{u}^{\prime}$. 
This completes the induction.
\qed\vskip 2mm

Now, we define 
\begin{equation}\label{eq:X action}
X_{j}^{\pm1}m=\alpha_{j}^{\pm}(m) \quad (m\in M,\ 1\le j\le \ell).
\end{equation}

\begin{lem}\label{lem:B.3}
$M$ is a $\aha$-module with respect to \eqref{eq:X action}, and $W$ is isomorphic to ${V}^{\otimes\ell}\otimes_{\aha}M$ as a $\mathcal{U}(\e)$-module.
\end{lem}
\pf The proof is almost identical to the one in \cite{CP}. So we leave it to the reader. \qed\vskip 2mm

This completes the proof for essential surjectivity of $\mathcal{F}$.

\begin{lem}
The functor $\mathcal{F}$ is fully faithful.
\end{lem}

\pf
First, $\mathcal{F}$ is faithful since $\mathcal{J}_{\ell}$
is faithful. So it suffices to show that $\mathcal{F}_{\ell}$ is surjective on morphisms. 

Suppose that $F:\mathcal{F}_{\ell}(M)\rightarrow\mathcal{F}_{\ell}(M^{\prime})$ is a $\mathcal{U}(\epsilon)$-linear map
for $M,\,M^{\prime}\in\aha\text{-mod}$. Since $\mathcal{J}_{\ell}$
is an equivalence, there is a $\fha$-linear map $f:M\rightarrow M^{\prime}$
such that $\mathcal{J}_{\ell}(f)=F$. 
Since $F$ is $\mathcal{U}(\epsilon)$-linear,
$e_{0}F\left(v\otimes m\right)=F\left(e_{0}\left(v\otimes m\right)\right)$.
The left hand side is equal to
\[
e_{0}F\left(v\otimes m\right)=e_{0}\left(v\otimes f(m)\right)=\sum_{j}\Delta_{j}(e_{0})v\otimes X_{j}f(m),
\]
while the right hand side is
\[
F\left(e_{0}\left(v\otimes m\right)\right)=F\left(\sum_{j}\Delta_{j}(e_{0})v\otimes X_{j}m\right)=\sum_{j}\Delta_{j}(e_{0})v\otimes f(X_{j}m).
\]
Now for each $i$, we can choose a vector $v(i)$ so that $\Delta_{j}(e_{0})v=0$
unless $j=i$, and at the same time $\Delta_{i}(e_{0})v$ is of the form $v_{i_{1}}\otimes\cdots\otimes v_{i_{\ell}}$, whose factors are all distinct $v_{k}$'s. For example, we
may take $v(1)=v_{1}\otimes v_{2}\otimes\cdots\otimes v_{\ell}$.
Putting $v=v(i)$ in the above identities, we obtain 
$X_{i}f(m)=f(X_{i}m)$ by Lemma~\ref{lem:pf-aff-SWD}.
Hence $f$ is $\aha$-linear as well. 
\qed\vskip 2mm
 
Therefore, $\mc{F}$ in \eqref{eq:CP equivalence} is an equivalence of categories.
Since every simple object in $\aha\text{{\rm -mod}}$ is a quotient of $L(a_1)\circ \dots \circ L(a_\ell)$ for some $a_1,\dots,a_\ell\in {\bf k}$, where $\circ$ is a convolution product, and $\mc{F}(L(a_1)\circ \dots \circ L(a_\ell))=\mathcal{W}_1(a_1)\ot \cdots\ot \mathcal{W}_1(a_\ell)$, $\mc {F}$ induces the equivalence  
\begin{equation*}\label{eq:CP equivalence-2}
\xymatrixcolsep{2pc}\xymatrixrowsep{0pc}\xymatrix{
\mc{F}_{\e,\ell}^{\ast} : \aha\text{{\rm -mod}}_J \ar@{->}[r] & \ \mc{C}^\ell_J(\e).}
\end{equation*}


{\small
\begin{thebibliography}{HK}

\bibitem{AK}
T. Akasaka, M. Kashiwara, {\em Finite-dimensional representations of quantum affine algebras}, Publ. Res. Inst. Math. Sci. \textbf{33} (1997) 839--867.

\bibitem{BKK}
G. Benkart, S.-J. Kang, M. Kashiwara, {\em Crystal bases for the
quantum superalgebra $U_q(\frak{gl}(m,n))$}, J. Amer. Math. Soc. \textbf{13} (2000) 295-331.


\bibitem{B}
J. Brundan, {\em Kazhdan–Lusztig polynomials and character formulae for the lie superalgebra $\mf{gl}(m|n)$}, J. Amer. Math. Soc. \textbf{16} (2003) 185--231.

\bibitem{BK} 
J. Brundan, A. Kleshchev, {\em Blocks of cyclotomic Hecke algebras
and Khovanov-Lauda algebras}, Invent. Math. \textbf{178} (3) (2009) 451--484.

\bibitem{BLW}
J. Brundan, I. Losev, B. Webster, {\em Tensor product categorifications and the super Kazhdan-Lusztig conjecture}, Int. Math. Res. Not. (2017) 6329--6410.

\bibitem{BS}
J. Brundan, C. Stroppel, {\em Highest weight categories arising from Khovanov’s diagram algebra IV: the general linear supergroup}, J. Eur. Math. Soc. \textbf{14} (2012) 373--419.

\bibitem{CP}
V. Chari, A. Pressley, {\em Quantum affine algebras and affine Hecke algebras}, Pacific J. Math. \textbf{174} (2) (1996) 295--326.

\bibitem{CL}
S.-J. Cheng, N. Lam, {\em Irreducible characters of general linear superalgebra and super duality}, Comm. Math. Phys. \textbf{298} (2010) 645--672.

\bibitem{CLW}
S.-J. Cheng, N. Lam, W. Wang, {\em Super duality and irreducible characters of ortho-symplectic Lie superalgebras}, Invent. Math. \textbf{183} (2011) 189--224.

\bibitem{CLW15}
S.-J. Cheng, N. Lam, W. Wang, {\em The Brundan-Kazhdan-Lusztig conjecture for general linear Lie superalgebras}, Duke Math. J. \textbf{164} (2015) 617--695.

\bibitem{CW}
S.-J. Cheng, W. Wang, {\em Dualities and Representations of Lie Superalgebras}, Graduate Studies in Mathematics 144, Amer. Math. Soc., 2013.

\bibitem{CWZ}
S.-J. Cheng, W. Wang, R. B.  Zhang, {\em Super duality and Kazhdan–Lusztig polynomials}, Trans. Amer. Math. Soc. \textbf{360} (2008) 5883--5924.

\bibitem{DO}
E. Date, M. Okado, {\em Calculation of excitation spectra of the spin model related with the vector representation of the quantized affine algebra of type $A_n^{(1)}$}, Internat. J. Modern Phys. A  \textbf{9} (1994) 399--417.

\bibitem{EA}
I. Entova-Aizenbud, {\em Notes on restricted inverse limits of categories}, preprint (2015), arXiv:1504.01121.

\bibitem{FR}
E. Frenkel, N. Reshetikhin, {\em The $q$-characters of representations of quantum affine algebras and deformations of $W$-algebras}, Recent Developments in Quantum Affine Algebras and related topics, Contemp. Math. \textbf{248} (1999) 163--205.

\bibitem{Her06}
D. Hernandez, {\em The Kirillov-Reshetikhin conjecture and solutions of T-systems}, J. Reine Angew. Math. \textbf{596} (2006) 63--87.

\bibitem{Her11}
D. Hernandez, {\em The algebra $U_q(\widehat{sl}_\infty)$ and applications}, J. Algebra \textbf{329} (2011) 147--162.

\bibitem{HL10}
D. Hernandez, B. Leclerc, {\em Cluster algebras and quantum affine algebras}, Duke Math. J. \textbf{154} (2010) 265--341.

\bibitem{HL16}
D. Hernandez, B. Leclerc, {\em Cluster algebras and category $\mc{O}$ for representations of Borel subalgebras of quantum affine algebras}, Algebra Number Theory \textbf{10} (2016) 2015--2052.

\bibitem{KKK}
S.-J. Kang, M. Kashiwara, M.-H. Kim, {\em Symmetric quiver Hecke algebras and $R$-matrices of quantum affine algebras}, Invent. math. \textbf{211} (2018) 591–685.

\bibitem{KKK2}
S.-J. Kang, M. Kashiwara, M.-H. Kim, {\em Symmetric quiver Hecke algebras and $R$-matrices of quantum affine algebras, II}, Duke Math. J. \textbf{164} (2015) 1549--1602.

\bibitem{KKKO15}
S.-J. Kang, M. Kashiwara, M.-H. Kim, S.-J. Oh, {\em Simplicity of heads and socles of tensor products}, Compos. Math. \textbf{151} (2015) 377--396.

\bibitem{KKKO18}
S.-J. Kang, M. Kashiwara, M.-H. Kim, S.-J. Oh, {\em Monoidal categorification of cluster algebras}, J. Amer. Math. Soc. \textbf{31} (2018) 349--426.

\bibitem{KKOP}
M. Kashiwara, M.-H. Kim, S.-J. Oh, E. Park, {\em Monoidal categorification and quantum affine algebras}, Compos. Math. \textbf{156} (2020) 1039--1077.

\bibitem{KMN}
S-J. Kang, M. Kashiwara, K. C. Misra, T. Miwa, T. Nakashima and A. Nakayashiki, 
{\em Affine crystals and vertex models}, Internat. J. Modern Phys. A \textbf{7} (suppl. 1A) (1992) 449--484.

\bibitem{KP}
S-J. Kang, E. Park, {\em Irreducible modules over Khovanov-Lauda-Rouquier algebras of type $A_n$ and semistandard tableaux}, J. Algebra \textbf{339} (2011) 223--251.

\bibitem{Kas02}
M. Kashiwara, {\em On level zero representations of quantum affine algebras}, Duke Math. J. \textbf{112} (2002) 117--175.

\bibitem{KL}
M. Khovanov, A. Lauda, {\em A diagrammatic approach to categorification of quantum groups I}, Represent. Theory \textbf{13} (2011) 309--347.

\bibitem{Kim}
M.-H. Kim, {\em Khovanov-Lauda-Rouquier algebras and $R$-matrices}, Ph.D. thesis, Seoul National University (2012).

\bibitem{KNS}
A. Kuniba, T. Nakanishi, J. Suzuki, {\em Functional relations in solvable lattice models: I. Functional relations and representation theory}, Internat. J. Modern Phys. A \textbf{9} (1994) 5215--5266.

\bibitem{KOS}
A. Kuniba, M. Okado, S. Sergeev, {\em Tetrahedron equation and generalized quantum groups}, 
J. Phys. A \textbf{48} (2015) 304001 (38p).

\bibitem{KO}
J.-H. Kwon, M. Okado, {\em Kirillov-Reshetikhin modules over generalized quantum groups of type $A$}, preprint (2018) arXiv:1804.05456, to appear in Publ. Res. Inst. Math. Sci.

\bibitem{KY}
J.-H. Kwon, J. Yu, {\em $R$-matrix for generalized quantum groups of type $A$}, J. Algebra, \textbf{566} (2021) 309--341.

\bibitem{Lu93}
G. Lusztig, {\em Introduction to quantum groups}, Progr. Math. \textbf{110}, Birkh\"{a}user, 1993.

\bibitem{Na03}
H. Nakajima, {\em $t$-analogs of $q$-characters of Kirillov-Reshetikhin modules of quantum affine algebras}, Represent. Theory \textbf{7} (2003) 259--274.

\bibitem{Na04}
H. Nakajima, {\em  Quiver varieties and $t$-analogs of $q$-characters of quantum affine algebras}, Ann. of  Math. \textbf{160} (2004) 1057--1097.

\bibitem{Na04-2}
H. Nakajima, {\em Extremal weight modules of quantum affine algebras}, Adv. Stud. Pure Math. \textbf{40} (2004) 343--369.

\bibitem{R}
R. Rouquier, {\em 2-Kac-Moody algebras}, preprint (2008) arXiv:0812.5023v1.

\bibitem{Sc}
O. M. Schn\"{u}rer, {\em Equivariant sheaves on flag varieties}, Math. Z. \textbf{267} (2011) 27--80.

\bibitem{Ser}
V. Serganova, {\em Characters of irreducible representations of simple lie superalgebras}, Doc. Math. Extra Volume ICM II (1998) 583--593.


\bibitem{Ya99}
H. Yamane, {\em On defining relations of affine Lie superalgebras and affine quantized universal enveloping superalgebras}, Publ. Res. Inst. Math. Sci. \textbf{35} (1999) 321--390.

\bibitem{Z14}
H. Zhang, {\em Representations of quantum affine superalgebras}, Math. Z. \textbf{278} (2014) 663--703.

\bibitem{Z16}
H. Zhang, {\em RTT realization of quantum affine superalgebras and tensor products}, Int. Math. Res. Not. \textbf{2016} (2016) 1126--1157.

\bibitem{Z17}
H. Zhang, {\em Fundamental representations of quantum affine superalgebras and $R$-matrices}, Transform. Groups \textbf{22} (2017) 559--590.


\end{thebibliography}}
\end{document}